\documentclass[a4paper,12pt]{article}
\usepackage{amsfonts, amsmath, amssymb, amscd, latexsym, graphicx,color}
\usepackage{a4,amsfonts,amsthm,amsmath,amssymb,manfnt}
\newcommand{\N}{\mathbb{N}}
\newcommand{\Z}{\mathbb{Z}}

\newcommand{\T}{\mathcal{T}(\varLambda)}
\newcommand{\FL}{\mathcal{F}(\varLambda)}
\newcommand{\BFL}{\text{BF}(\varLambda)}
\newcommand{\BFS}{\text{BF}(\varLambda;\mathcal{S})}
\newcommand{\BFE}{\text{BF}(\mathcal E;\mathcal{S})}
\newcommand{\C}{{\text{\rm Compl}}}
\newcommand{\im}{\text{im}}
\newcommand{\Aut}{\text{Aut}}
\newtheorem{theorem}{Theorem}[subsection]
\newtheorem{lemma}[theorem]{Lemma}
\newtheorem{defn}[theorem]{Definition}
\newtheorem{cor}[theorem]{Corollary}
\newtheorem{rmk}[theorem]{Remark}

\newtheorem{prop}[theorem]{Proposition}
\newtheorem{notation}[theorem]{Notation}
\newtheorem{ex}[theorem]{Example}

\newcommand{\demo}{ {\it   Proof. }}

\title{Generalized presentations of infinite groups,\\ in particular of ${\text{\rm Aut}}(F_{\omega})$}
\author{\hspace*{-20mm}O. Bogopolski\hspace*{50mm} W. Singhof\\
\hspace*{0mm}\begin{minipage}[t]{150mm}
\small{Institute of Mathematics of Siberian}\hspace*{15mm} {\small D\"{u}sseldorf University, Germany}\\ \small{Branch of Russian Academy
of Sciences,}\hspace*{10mm}\small{e-mail:
singhof@math.uni-duesseldorf.de} \\ {\small Novosibirsk, Russia}\\ {\small and D\"{u}sseldorf University, Germany} \\ \small{e-mail:
Oleg$\_$Bogopolski@yahoo.com}
\end{minipage}
\\}

\begin{document}

\medskip

\maketitle

\begin{abstract}

We develop a theory of generalized presentations of groups.
We give generalized presentations of the symmetric group $\varSigma(X)$ for an arbitrary set $X$
and of the automorphism group of the free group of countable rank, $\text{\rm Aut}(F_{\omega})$.

\end{abstract}

\tableofcontents
\setcounter{tocdepth}{2}

\section{Introduction}

According to classical group theory, a group $G$ is generated by a subset $\varLambda$ if every element of $G$
is a {\it finite product} of elements of $\varLambda\cup \varLambda^{-1}$.
In particular, if $G$ is uncountable, it cannot be generated by
a countable set.
However, if we allow appropriate {\it infinite products} this
becomes possible for at least the two types of groups mentioned in the abstract.
We show that $\varSigma(X)$ and  ${\text{\rm Aut}}(F_{\omega})$ can be generated in this generalized sense by transpositions and by elementary Nielsen automorphisms of the first kind, respectively.
Moreover, we describe {\it generalized presentations} (see Section~\ref{Generalized presentations}) of these groups on these sets of generalized generators.

In classical group theory, the free groups serve as universal objects;
we have to replace them in our theory by what we call {\it generalized free groups}, which include the {\it big free groups}.

In Section~\ref{permut} we describe two generalized presentations of the group $\varSigma(X)$,
see Theorems~\ref{ThSigma1} and~\ref{ThSigma2}. In Section~\ref{GpAut} we describe a generalized presentation of the group ${\text{\rm Aut}}(F_{\omega})$, see Theorem~\ref{kernel2}. A more algebraic description is given in
Theorems~\ref{kernel1} and~\ref{kernel3}.

Whereas in the case of the symmetric  groups
$\varSigma(X)$, we admit arbitrary sets $X$, for the automorphism groups ${\text{\rm Aut}}(F(X))$ 
we have to assume that $X$ is countable.
For bigger sets $X$, the corresponding questions about ${\text{\rm Aut}}(F(X))$
remain open. In Appendix~B, we formulate a number of questions.

\section{Generalized free groups and generalized\\ presentations}\label{BF and GP}

Throughout the paper we use the following notations.
For a subset $R$ of a group $G$, we denote by $\langle R\rangle_G$ the subgroup generated by $R$ and
by $\langle\!\langle R\rangle\!\rangle_G$ the normal closure of $R$ in $G$.
We skip the subscript $G$ if the ambient group is clear from a context.
If $G$ is a topological group with topology $\frak{T}$, we denote by $\overline{R}^{\frak{T}}$ the (topological)
closure of the subset $R$ in $G$.

\subsection{Big free groups}~\label{Bfg}

We begin with a few recollections concerning big free groups~\cite{CC1}. Let $\varLambda$ be a set.
 By $\T$ we denote the set of all maps
\[f:S \rightarrow \varLambda\]

with the following properties:
\begin{itemize}
 \item $S$ is a totally ordered set.
 \item For each $\lambda \in \varLambda$, the set $f^{-1}(\lambda)$ is finite.
\end{itemize}

We identify two elements $f:S \rightarrow \varLambda$ and $f':S' \rightarrow \varLambda$ of $\T$ if there is
an order preserving bijection $\varphi:S\rightarrow S'$ with $f=f' \circ \varphi$.

\smallskip
For the remainder of Section~\ref{BF and GP},
we assume that we are given a free involution $\lambda\mapsto\lambda^{-1}$ on $\varLambda$.

\emph{Definition.} A subset $I$ of a linearly ordered set $S$ is called an \emph{interval} of $S$ if the following holds:
for $s,s',s'' \in S$ with $s< s'' < s'$ and $s, s' \in I$, we have $s'' \in I$.

\smallskip
We define
\[ [s,s']:=[s,s']_S:=\{x \in S \mid s \leq x \leq s'\}\,.\]

Other types of intervals such as $]\,s,s'\,[$ or $]\,s,\infty\,[$ are defined accordingly.

\medskip
Given $f,g \in \T$, we say that $g$ \emph{is obtained from} $f$ \emph{by cancellation} and write $f\searrow g$,
if the following holds:
\begin{itemize}
 \item If $f$ is of the form $f:S\rightarrow \varLambda$, there is $T\subseteq S$ such that $g =f_{\mid S \smallsetminus T}$.
 \item There is an involution $\ast$ on $T$ such that for all $t \in T$ we have

 \begin{eqnarray}  f(t^{\ast})=f(t)^{-1} \label{S1}\\[1ex]
   [t,t^{\ast}] \subseteq T\qquad \label{S2}\\[1ex]
  ([t,t^{\ast}])^{\ast} = [t,t^{\ast}]\, .\label{S3}
 \end{eqnarray}
\end{itemize}

Let $\approx$ be the equivalence relation on $\T$ generated by $f \searrow g$ and let
\[\text{BF}(\varLambda):= \T /\approx\, .\]

By $[f] \in \BFL$ we denote the class of $f \in \T$. Then $\BFL$ becomes a group by $[f]\,[f']:=[f f']$
where $ff'$ is the concatenation of $f$ and $f'$. We call $\BFL$ the \emph{big free group} over
$\varLambda$.
Note that BF is a functor from the category of sets with free
involution to the category of groups.
Recall the following definition from ~\cite{CC1}:

\medskip
\emph{Definition.} An element $f \in \T$ is called \emph{reduced} if no element of $\T$ except $f$ itself can be
obtained from $f$ by cancellation.

\medskip

As shown in~\cite[Theorem 3.9$'$]{CC1}, every element of $\BFL$ admits a unique reduced representing element in $\T$.
This important property (which is analogous to the uniqueness property for free groups) implies the following statements.

\medskip

1) Any injective map of sets with free involutions $\varLambda'\rightarrow \varLambda$ induces a monomorphism of big free groups
${\text{\rm BF}}(\varLambda')\rightarrow{\text{\rm BF}}(\varLambda)$.

2) The natural embedding $\varLambda\rightarrow \T$, given by the rule $\lambda\mapsto (f:\{\ast\}\rightarrow \varLambda)$, where $f(\ast)=\lambda$, induces the natural embedding $\varLambda\rightarrow \BFL$. So, we will identify $\varLambda$ with its image in $\BFL$.

\smallskip

Let $\varLambda^+$ be a subset of $\varLambda$ containing, for every $\lambda \in \varLambda$, exactly one of the two elements $\lambda,\lambda^{-1}$.
The subgroup of $\BFL$ generated (in the classical sense) by $\varLambda$ will be denoted by $\FL$.
Clearly, $\FL$ is isomorphic to the free group $F(\varLambda^+)$ with basis $\varLambda^+$.

It is obvious that if $\varLambda$ is finite,
then $\BFL$ coincides with $\FL$. However, for infinite $\varLambda$ the big free group $\BFL$ is not free. To explain this, let us recall  the definition of the Hawaiian Earring.

\medskip

{\it Definition. } {\it The Hawaiian Earring} $\mathcal{H}$ is the topological space which is the countable union of circles of radii $\frac{1}{n}$, $n=1,2,\dots $, embedded into the Euclidean plane in such a way that they have only
one common point $x$.
The topology of the Hawaiian Earring is induced by the topology of the plane.

\vspace*{-1mm}
\hspace*{-8.8mm}
\unitlength 1mm 
\linethickness{0.4pt}
\ifx\plotpoint\undefined\newsavebox{\plotpoint}\fi 
\begin{picture}(76.39,38)(-15,32)
\put(70,50){\circle{12.01}}
\put(76.39,50){\line(0,1){.57}}
\put(76.37,50.57){\line(0,1){.568}}
\put(76.33,51.14){\line(0,1){.565}}
\put(76.25,51.7){\line(0,1){.56}}
\multiput(76.14,52.26)(-.02721,.11066){5}{\line(0,1){.11066}}
\multiput(76.01,52.82)(-.0331,.10905){5}{\line(0,1){.10905}}
\multiput(75.84,53.36)(-.03241,.08926){6}{\line(0,1){.08926}}
\multiput(75.65,53.9)(-.03184,.07491){7}{\line(0,1){.07491}}
\multiput(75.42,54.42)(-.03133,.06396){8}{\line(0,1){.06396}}
\multiput(75.17,54.93)(-.03085,.05528){9}{\line(0,1){.05528}}
\multiput(74.9,55.43)(-.03039,.0482){10}{\line(0,1){.0482}}
\multiput(74.59,55.91)(-.03293,.0465){10}{\line(0,1){.0465}}
\multiput(74.26,56.38)(-.03216,.04061){11}{\line(0,1){.04061}}
\multiput(73.91,56.82)(-.03143,.03559){12}{\line(0,1){.03559}}
\multiput(73.53,57.25)(-.03329,.03386){12}{\line(0,1){.03386}}
\multiput(73.13,57.66)(-.03505,.03203){12}{\line(-1,0){.03505}}
\multiput(72.71,58.04)(-.04006,.03284){11}{\line(-1,0){.04006}}
\multiput(72.27,58.4)(-.04593,.03372){10}{\line(-1,0){.04593}}
\multiput(71.81,58.74)(-.04767,.03121){10}{\line(-1,0){.04767}}
\multiput(71.33,59.05)(-.05475,.03179){9}{\line(-1,0){.05475}}
\multiput(70.84,59.34)(-.06342,.03241){8}{\line(-1,0){.06342}}
\multiput(70.33,59.6)(-.07436,.03311){7}{\line(-1,0){.07436}}
\multiput(69.81,59.83)(-.07603,.02908){7}{\line(-1,0){.07603}}
\multiput(69.28,60.03)(-.09039,.02912){6}{\line(-1,0){.09039}}
\multiput(68.74,60.21)(-.11018,.02909){5}{\line(-1,0){.11018}}
\put(68.19,60.35){\line(-1,0){.558}}
\put(67.63,60.47){\line(-1,0){.563}}
\put(67.07,60.55){\line(-1,0){.567}}
\put(66.5,60.61){\line(-1,0){.569}}
\put(65.93,60.63){\line(-1,0){.57}}
\put(65.36,60.63){\line(-1,0){.569}}
\put(64.79,60.59){\line(-1,0){.566}}
\put(64.23,60.53){\line(-1,0){.562}}
\multiput(63.67,60.43)(-.1389,-.0317){4}{\line(-1,0){.1389}}
\multiput(63.11,60.3)(-.10959,-.03123){5}{\line(-1,0){.10959}}
\multiput(62.56,60.15)(-.0898,-.03088){6}{\line(-1,0){.0898}}
\multiput(62.02,59.96)(-.07545,-.03055){7}{\line(-1,0){.07545}}
\multiput(61.49,59.75)(-.06449,-.03023){8}{\line(-1,0){.06449}}
\multiput(60.98,59.51)(-.06278,-.03364){8}{\line(-1,0){.06278}}
\multiput(60.48,59.24)(-.05412,-.03285){9}{\line(-1,0){.05412}}
\multiput(59.99,58.94)(-.04706,-.03213){10}{\line(-1,0){.04706}}
\multiput(59.52,58.62)(-.04115,-.03146){11}{\line(-1,0){.04115}}
\multiput(59.07,58.27)(-.03941,-.03362){11}{\line(-1,0){.03941}}
\multiput(58.63,57.9)(-.03442,-.03271){12}{\line(-1,0){.03442}}
\multiput(58.22,57.51)(-.03262,-.0345){12}{\line(0,-1){.0345}}
\multiput(57.83,57.1)(-.03352,-.03949){11}{\line(0,-1){.03949}}
\multiput(57.46,56.66)(-.03136,-.04123){11}{\line(0,-1){.04123}}
\multiput(57.11,56.21)(-.03201,-.04713){10}{\line(0,-1){.04713}}
\multiput(56.79,55.74)(-.03272,-.0542){9}{\line(0,-1){.0542}}
\multiput(56.5,55.25)(-.03349,-.06286){8}{\line(0,-1){.06286}}
\multiput(56.23,54.75)(-.03007,-.06456){8}{\line(0,-1){.06456}}
\multiput(55.99,54.23)(-.03037,-.07552){7}{\line(0,-1){.07552}}
\multiput(55.78,53.7)(-.03066,-.08988){6}{\line(0,-1){.08988}}
\multiput(55.6,53.16)(-.03096,-.10967){5}{\line(0,-1){.10967}}
\multiput(55.44,52.61)(-.0313,-.139){4}{\line(0,-1){.139}}
\put(55.32,52.06){\line(0,-1){.562}}
\put(55.22,51.5){\line(0,-1){.566}}
\put(55.15,50.93){\line(0,-1){1.139}}
\put(55.12,49.79){\line(0,-1){.569}}
\put(55.14,49.22){\line(0,-1){.567}}
\put(55.2,48.66){\line(0,-1){.563}}
\put(55.29,48.09){\line(0,-1){.558}}
\multiput(55.4,47.54)(.02936,-.11011){5}{\line(0,-1){.11011}}
\multiput(55.55,46.99)(.02935,-.09032){6}{\line(0,-1){.09032}}
\multiput(55.73,46.44)(.02926,-.07596){7}{\line(0,-1){.07596}}
\multiput(55.93,45.91)(.03329,-.07428){7}{\line(0,-1){.07428}}
\multiput(56.16,45.39)(.03257,-.06334){8}{\line(0,-1){.06334}}
\multiput(56.42,44.88)(.03192,-.05467){9}{\line(0,-1){.05467}}
\multiput(56.71,44.39)(.03132,-.0476){10}{\line(0,-1){.0476}}
\multiput(57.03,43.92)(.03075,-.04168){11}{\line(0,-1){.04168}}
\multiput(57.36,43.46)(.03294,-.03998){11}{\line(0,-1){.03998}}
\multiput(57.73,43.02)(.03211,-.03497){12}{\line(0,-1){.03497}}
\multiput(58.11,42.6)(.03394,-.0332){12}{\line(1,0){.03394}}
\multiput(58.52,42.2)(.03567,-.03134){12}{\line(1,0){.03567}}
\multiput(58.95,41.82)(.04069,-.03206){11}{\line(1,0){.04069}}
\multiput(59.39,41.47)(.04658,-.03281){10}{\line(1,0){.04658}}
\multiput(59.86,41.14)(.05364,-.03363){9}{\line(1,0){.05364}}
\multiput(60.34,40.84)(.05536,-.03071){9}{\line(1,0){.05536}}
\multiput(60.84,40.56)(.06404,-.03117){8}{\line(1,0){.06404}}
\multiput(61.35,40.32)(.07499,-.03165){7}{\line(1,0){.07499}}
\multiput(61.88,40.09)(.08934,-.03219){6}{\line(1,0){.08934}}
\multiput(62.41,39.9)(.10913,-.03283){5}{\line(1,0){.10913}}
\multiput(62.96,39.74)(.1384,-.0337){4}{\line(1,0){.1384}}
\put(63.51,39.6){\line(1,0){.56}}
\put(64.07,39.5){\line(1,0){.565}}
\put(64.64,39.42){\line(1,0){.568}}
\put(65.21,39.38){\line(1,0){1.139}}
\put(66.35,39.38){\line(1,0){.568}}
\put(66.91,39.43){\line(1,0){.564}}
\put(67.48,39.51){\line(1,0){.56}}
\multiput(68.04,39.61)(.11059,.02748){5}{\line(1,0){.11059}}
\multiput(68.59,39.75)(.10896,.03337){5}{\line(1,0){.10896}}
\multiput(69.14,39.92)(.08918,.03263){6}{\line(1,0){.08918}}
\multiput(69.67,40.11)(.07484,.03202){7}{\line(1,0){.07484}}
\multiput(70.19,40.34)(.06389,.03148){8}{\line(1,0){.06389}}
\multiput(70.71,40.59)(.05521,.03099){9}{\line(1,0){.05521}}
\multiput(71.2,40.87)(.04812,.03051){10}{\line(1,0){.04812}}
\multiput(71.68,41.17)(.04642,.03304){10}{\line(1,0){.04642}}
\multiput(72.15,41.5)(.04053,.03226){11}{\line(1,0){.04053}}
\multiput(72.59,41.86)(.03552,.03151){12}{\line(1,0){.03552}}
\multiput(73.02,42.24)(.03378,.03337){12}{\line(1,0){.03378}}
\multiput(73.43,42.64)(.03194,.03513){12}{\line(0,1){.03513}}
\multiput(73.81,43.06)(.03274,.04014){11}{\line(0,1){.04014}}
\multiput(74.17,43.5)(.0336,.04602){10}{\line(0,1){.04602}}
\multiput(74.5,43.96)(.03109,.04775){10}{\line(0,1){.04775}}
\multiput(74.82,44.44)(.03165,.05483){9}{\line(0,1){.05483}}
\multiput(75.1,44.93)(.03226,.0635){8}{\line(0,1){.0635}}
\multiput(75.36,45.44)(.03292,.07444){7}{\line(0,1){.07444}}
\multiput(75.59,45.96)(.03371,.08878){6}{\line(0,1){.08878}}
\multiput(75.79,46.49)(.0289,.09046){6}{\line(0,1){.09046}}
\multiput(75.96,47.04)(.02882,.11025){5}{\line(0,1){.11025}}
\put(76.11,47.59){\line(0,1){.558}}
\put(76.22,48.14){\line(0,1){.564}}
\put(76.31,48.71){\line(0,1){.567}}
\put(76.36,49.28){\line(0,1){.724}}
\put(76.31,50){\line(0,1){.692}}
\put(76.29,50.69){\line(0,1){.69}}
\put(76.24,51.38){\line(0,1){.686}}
\multiput(76.15,52.07)(-.0307,.1702){4}{\line(0,1){.1702}}
\multiput(76.03,52.75)(-.0315,.13474){5}{\line(0,1){.13474}}
\multiput(75.87,53.42)(-.03194,.11079){6}{\line(0,1){.11079}}
\multiput(75.68,54.09)(-.03219,.09345){7}{\line(0,1){.09345}}
\multiput(75.45,54.74)(-.0323,.08022){8}{\line(0,1){.08022}}
\multiput(75.19,55.38)(-.03231,.06975){9}{\line(0,1){.06975}}
\multiput(74.9,56.01)(-.03224,.06121){10}{\line(0,1){.06121}}
\multiput(74.58,56.62)(-.03211,.05408){11}{\line(0,1){.05408}}
\multiput(74.23,57.22)(-.03193,.04801){12}{\line(0,1){.04801}}
\multiput(73.84,57.79)(-.031694,.042752){13}{\line(0,1){.042752}}
\multiput(73.43,58.35)(-.031417,.038145){14}{\line(0,1){.038145}}
\multiput(72.99,58.88)(-.033322,.036493){14}{\line(0,1){.036493}}
\multiput(72.53,59.39)(-.032797,.032429){15}{\line(-1,0){.032797}}
\multiput(72.03,59.88)(-.036867,.032908){14}{\line(-1,0){.036867}}
\multiput(71.52,60.34)(-.041459,.033368){13}{\line(-1,0){.041459}}
\multiput(70.98,60.78)(-.043107,.031209){13}{\line(-1,0){.043107}}
\multiput(70.42,61.18)(-.04836,.03138){12}{\line(-1,0){.04836}}
\multiput(69.84,61.56)(-.05444,.0315){11}{\line(-1,0){.05444}}
\multiput(69.24,61.9)(-.06157,.03155){10}{\line(-1,0){.06157}}
\multiput(68.62,62.22)(-.07011,.03152){9}{\line(-1,0){.07011}}
\multiput(67.99,62.5)(-.08058,.03139){8}{\line(-1,0){.08058}}
\multiput(67.35,62.75)(-.0938,.03113){7}{\line(-1,0){.0938}}
\multiput(66.69,62.97)(-.11115,.03069){6}{\line(-1,0){.11115}}
\multiput(66.02,63.16)(-.13508,.02998){5}{\line(-1,0){.13508}}
\put(65.35,63.31){\line(-1,0){.682}}
\put(64.67,63.42){\line(-1,0){.687}}
\put(63.98,63.5){\line(-1,0){.69}}
\put(63.29,63.55){\line(-1,0){.692}}
\put(62.6,63.56){\line(-1,0){.691}}
\put(61.91,63.53){\line(-1,0){.689}}
\put(61.22,63.47){\line(-1,0){.685}}
\multiput(60.53,63.37)(-.1699,-.0326){4}{\line(-1,0){.1699}}
\multiput(59.85,63.24)(-.13437,-.03302){5}{\line(-1,0){.13437}}
\multiput(59.18,63.08)(-.11043,-.03319){6}{\line(-1,0){.11043}}
\multiput(58.52,62.88)(-.09308,-.03324){7}{\line(-1,0){.09308}}
\multiput(57.87,62.65)(-.07985,-.0332){8}{\line(-1,0){.07985}}
\multiput(57.23,62.38)(-.06938,-.0331){9}{\line(-1,0){.06938}}
\multiput(56.6,62.08)(-.06084,-.03293){10}{\line(-1,0){.06084}}
\multiput(55.99,61.75)(-.05371,-.03272){11}{\line(-1,0){.05371}}
\multiput(55.4,61.39)(-.04764,-.03247){12}{\line(-1,0){.04764}}
\multiput(54.83,61.01)(-.042392,-.032174){13}{\line(-1,0){.042392}}
\multiput(54.28,60.59)(-.037788,-.031845){14}{\line(-1,0){.037788}}
\multiput(53.75,60.14)(-.036115,-.033731){14}{\line(-1,0){.036115}}
\multiput(53.25,59.67)(-.032057,-.033161){15}{\line(0,-1){.033161}}
\multiput(52.77,59.17)(-.03249,-.037236){14}{\line(0,-1){.037236}}
\multiput(52.31,58.65)(-.032898,-.041833){13}{\line(0,-1){.041833}}
\multiput(51.88,58.11)(-.03328,-.04708){12}{\line(0,-1){.04708}}
\multiput(51.48,57.54)(-.03364,-.05314){11}{\line(0,-1){.05314}}
\multiput(51.11,56.96)(-.03088,-.05479){11}{\line(0,-1){.05479}}
\multiput(50.77,56.35)(-.03085,-.06192){10}{\line(0,-1){.06192}}
\multiput(50.46,55.73)(-.03073,-.07046){9}{\line(0,-1){.07046}}
\multiput(50.19,55.1)(-.03048,-.08093){8}{\line(0,-1){.08093}}
\multiput(49.94,54.45)(-.03007,-.09415){7}{\line(0,-1){.09415}}
\multiput(49.73,53.79)(-.02943,-.11149){6}{\line(0,-1){.11149}}
\multiput(49.56,53.13)(-.02845,-.13541){5}{\line(0,-1){.13541}}
\put(49.42,52.45){\line(0,-1){.683}}
\put(49.31,51.76){\line(0,-1){.688}}
\put(49.24,51.08){\line(0,-1){2.762}}
\put(49.3,48.31){\line(0,-1){.684}}
\multiput(49.4,47.63)(.02765,-.13558){5}{\line(0,-1){.13558}}
\multiput(49.54,46.95)(.02878,-.11166){6}{\line(0,-1){.11166}}
\multiput(49.71,46.28)(.02952,-.09432){7}{\line(0,-1){.09432}}
\multiput(49.92,45.62)(.03,-.08111){8}{\line(0,-1){.08111}}
\multiput(50.16,44.97)(.03031,-.07064){9}{\line(0,-1){.07064}}
\multiput(50.43,44.34)(.03049,-.0621){10}{\line(0,-1){.0621}}
\multiput(50.74,43.72)(.03362,-.06047){10}{\line(0,-1){.06047}}
\multiput(51.07,43.11)(.03333,-.05334){11}{\line(0,-1){.05334}}
\multiput(51.44,42.53)(.033,-.04727){12}{\line(0,-1){.04727}}
\multiput(51.83,41.96)(.032651,-.042026){13}{\line(0,-1){.042026}}
\multiput(52.26,41.41)(.03227,-.037427){14}{\line(0,-1){.037427}}
\multiput(52.71,40.89)(.031861,-.03335){15}{\line(0,-1){.03335}}
\multiput(53.19,40.39)(.033521,-.031681){15}{\line(1,0){.033521}}
\multiput(53.69,39.91)(.0376,-.032067){14}{\line(1,0){.0376}}
\multiput(54.22,39.46)(.042201,-.032424){13}{\line(1,0){.042201}}
\multiput(54.77,39.04)(.04745,-.03275){12}{\line(1,0){.04745}}
\multiput(55.34,38.65)(.05352,-.03304){11}{\line(1,0){.05352}}
\multiput(55.93,38.29)(.06065,-.03329){10}{\line(1,0){.06065}}
\multiput(56.53,37.95)(.06919,-.0335){9}{\line(1,0){.06919}}
\multiput(57.15,37.65)(.07966,-.03367){8}{\line(1,0){.07966}}
\multiput(57.79,37.38)(.08127,-.02957){8}{\line(1,0){.08127}}
\multiput(58.44,37.15)(.09448,-.02901){7}{\line(1,0){.09448}}
\multiput(59.1,36.94)(.11181,-.02817){6}{\line(1,0){.11181}}
\multiput(59.77,36.77)(.1697,-.0336){4}{\line(1,0){.1697}}
\put(60.45,36.64){\line(1,0){.685}}
\put(61.14,36.54){\line(1,0){.689}}
\put(61.83,36.47){\line(1,0){.691}}
\put(62.52,36.44){\line(1,0){.692}}
\put(63.21,36.45){\line(1,0){.691}}
\put(63.9,36.49){\line(1,0){.688}}
\put(64.59,36.57){\line(1,0){.683}}
\multiput(65.27,36.68)(.13526,.02918){5}{\line(1,0){.13526}}
\multiput(65.95,36.82)(.11133,.03003){6}{\line(1,0){.11133}}
\multiput(66.61,37)(.09399,.03058){7}{\line(1,0){.09399}}
\multiput(67.27,37.22)(.08077,.03092){8}{\line(1,0){.08077}}
\multiput(67.92,37.47)(.0703,.03111){9}{\line(1,0){.0703}}
\multiput(68.55,37.75)(.06176,.03119){10}{\line(1,0){.06176}}
\multiput(69.17,38.06)(.05462,.03118){11}{\line(1,0){.05462}}
\multiput(69.77,38.4)(.04855,.0311){12}{\line(1,0){.04855}}
\multiput(70.35,38.77)(.0469,.03353){12}{\line(1,0){.0469}}
\multiput(70.91,39.18)(.041655,.033123){13}{\line(1,0){.041655}}
\multiput(71.46,39.61)(.03706,.03269){14}{\line(1,0){.03706}}
\multiput(71.98,40.06)(.032988,.032235){15}{\line(1,0){.032988}}
\multiput(72.47,40.55)(.033536,.036296){14}{\line(0,1){.036296}}
\multiput(72.94,41.06)(.031641,.037959){14}{\line(0,1){.037959}}
\multiput(73.38,41.59)(.031945,.042564){13}{\line(0,1){.042564}}
\multiput(73.8,42.14)(.03221,.04782){12}{\line(0,1){.04782}}
\multiput(74.18,42.71)(.03243,.05389){11}{\line(0,1){.05389}}
\multiput(74.54,43.31)(.0326,.06102){10}{\line(0,1){.06102}}
\multiput(74.87,43.92)(.03272,.06956){9}{\line(0,1){.06956}}
\multiput(75.16,44.54)(.03277,.08003){8}{\line(0,1){.08003}}
\multiput(75.42,45.18)(.03274,.09326){7}{\line(0,1){.09326}}
\multiput(75.65,45.84)(.03259,.1106){6}{\line(0,1){.1106}}
\multiput(75.85,46.5)(.03229,.13455){5}{\line(0,1){.13455}}
\multiput(76.01,47.17)(.0317,.17){4}{\line(0,1){.17}}
\put(76.14,47.85){\line(0,1){.686}}
\put(76.23,48.54){\line(0,1){.689}}
\put(76.29,49.23){\line(0,1){.771}}
\put(76.33,50){\line(0,1){.452}}
\put(76.32,50.45){\line(0,1){.451}}
\put(76.28,50.9){\line(0,1){.448}}
\put(76.21,51.35){\line(0,1){.444}}
\put(76.13,51.8){\line(0,1){.438}}
\multiput(76.01,52.23)(-.02744,.08623){5}{\line(0,1){.08623}}
\multiput(75.88,52.66)(-.03223,.08456){5}{\line(0,1){.08456}}
\multiput(75.71,53.09)(-.03076,.06885){6}{\line(0,1){.06885}}
\multiput(75.53,53.5)(-.02963,.05745){7}{\line(0,1){.05745}}
\multiput(75.32,53.9)(-.0328,.0557){7}{\line(0,1){.0557}}
\multiput(75.09,54.29)(-.03138,.04705){8}{\line(0,1){.04705}}
\multiput(74.84,54.67)(-.0302,.0402){9}{\line(0,1){.0402}}
\multiput(74.57,55.03)(-.0324,.03844){9}{\line(0,1){.03844}}
\multiput(74.28,55.38)(-.03105,.03291){10}{\line(0,1){.03291}}
\multiput(73.97,55.71)(-.03284,.03112){10}{\line(-1,0){.03284}}
\multiput(73.64,56.02)(-.03837,.03248){9}{\line(-1,0){.03837}}
\multiput(73.29,56.31)(-.04013,.03028){9}{\line(-1,0){.04013}}
\multiput(72.93,56.58)(-.04698,.03149){8}{\line(-1,0){.04698}}
\multiput(72.56,56.83)(-.05563,.03292){7}{\line(-1,0){.05563}}
\multiput(72.17,57.06)(-.05738,.02976){7}{\line(-1,0){.05738}}
\multiput(71.77,57.27)(-.06879,.03091){6}{\line(-1,0){.06879}}
\multiput(71.35,57.46)(-.08449,.03241){5}{\line(-1,0){.08449}}
\multiput(70.93,57.62)(-.08617,.02763){5}{\line(-1,0){.08617}}
\put(70.5,57.76){\line(-1,0){.438}}
\put(70.06,57.87){\line(-1,0){.444}}
\put(69.62,57.96){\line(-1,0){.448}}
\put(69.17,58.03){\line(-1,0){.451}}
\put(68.72,58.06){\line(-1,0){.905}}
\put(67.82,58.07){\line(-1,0){.451}}
\put(67.36,58.03){\line(-1,0){.448}}
\put(66.92,57.97){\line(-1,0){.444}}
\put(66.47,57.88){\line(-1,0){.438}}
\multiput(66.03,57.77)(-.08629,-.02726){5}{\line(-1,0){.08629}}
\multiput(65.6,57.63)(-.08463,-.03205){5}{\line(-1,0){.08463}}
\multiput(65.18,57.47)(-.06892,-.03061){6}{\line(-1,0){.06892}}
\multiput(64.77,57.29)(-.05751,-.0295){7}{\line(-1,0){.05751}}
\multiput(64.36,57.08)(-.05577,-.03268){7}{\line(-1,0){.05577}}
\multiput(63.97,56.85)(-.04712,-.03128){8}{\line(-1,0){.04712}}
\multiput(63.6,56.6)(-.04026,-.03011){9}{\line(-1,0){.04026}}
\multiput(63.23,56.33)(-.03851,-.03231){9}{\line(-1,0){.03851}}
\multiput(62.89,56.04)(-.03298,-.03098){10}{\line(-1,0){.03298}}
\multiput(62.56,55.73)(-.03119,-.03278){10}{\line(0,-1){.03278}}
\multiput(62.24,55.4)(-.03257,-.0383){9}{\line(0,-1){.0383}}
\multiput(61.95,55.06)(-.03037,-.04007){9}{\line(0,-1){.04007}}
\multiput(61.68,54.7)(-.03159,-.04692){8}{\line(0,-1){.04692}}
\multiput(61.43,54.32)(-.03304,-.05556){7}{\line(0,-1){.05556}}
\multiput(61.19,53.93)(-.02988,-.05732){7}{\line(0,-1){.05732}}
\multiput(60.99,53.53)(-.03106,-.06872){6}{\line(0,-1){.06872}}
\multiput(60.8,53.12)(-.0326,-.08442){5}{\line(0,-1){.08442}}
\multiput(60.64,52.7)(-.02782,-.08611){5}{\line(0,-1){.08611}}
\put(60.5,52.27){\line(0,-1){.438}}
\put(60.38,51.83){\line(0,-1){.443}}
\put(60.29,51.39){\line(0,-1){.448}}
\put(60.23,50.94){\line(0,-1){.451}}
\put(60.19,50.49){\line(0,-1){1.356}}
\put(60.22,49.13){\line(0,-1){.448}}
\put(60.28,48.68){\line(0,-1){.444}}
\put(60.37,48.24){\line(0,-1){.439}}
\multiput(60.48,47.8)(.02707,-.08635){5}{\line(0,-1){.08635}}
\multiput(60.61,47.37)(.03186,-.0847){5}{\line(0,-1){.0847}}
\multiput(60.77,46.95)(.03046,-.06899){6}{\line(0,-1){.06899}}
\multiput(60.95,46.53)(.02938,-.05758){7}{\line(0,-1){.05758}}
\multiput(61.16,46.13)(.03256,-.05584){7}{\line(0,-1){.05584}}
\multiput(61.39,45.74)(.03118,-.04719){8}{\line(0,-1){.04719}}
\multiput(61.64,45.36)(.03002,-.04033){9}{\line(0,-1){.04033}}
\multiput(61.91,45)(.03223,-.03858){9}{\line(0,-1){.03858}}
\multiput(62.2,44.65)(.03091,-.03305){10}{\line(0,-1){.03305}}
\multiput(62.51,44.32)(.03271,-.03126){10}{\line(1,0){.03271}}
\multiput(62.83,44.01)(.03823,-.03265){9}{\line(1,0){.03823}}
\multiput(63.18,43.71)(.04,-.03046){9}{\line(1,0){.04}}
\multiput(63.54,43.44)(.04685,-.03169){8}{\line(1,0){.04685}}
\multiput(63.91,43.19)(.05548,-.03316){7}{\line(1,0){.05548}}
\multiput(64.3,42.95)(.05725,-.03){7}{\line(1,0){.05725}}
\multiput(64.7,42.74)(.06865,-.03121){6}{\line(1,0){.06865}}
\multiput(65.11,42.56)(.08435,-.03278){5}{\line(1,0){.08435}}
\multiput(65.54,42.39)(.08605,-.02801){5}{\line(1,0){.08605}}
\put(65.97,42.25){\line(1,0){.437}}
\put(66.4,42.14){\line(1,0){.443}}
\put(66.85,42.05){\line(1,0){.448}}
\put(67.29,41.98){\line(1,0){.451}}
\put(67.75,41.94){\line(1,0){.905}}
\put(68.65,41.93){\line(1,0){.451}}
\put(69.1,41.97){\line(1,0){.448}}
\put(69.55,42.03){\line(1,0){.444}}
\put(69.99,42.11){\line(1,0){.439}}
\multiput(70.43,42.22)(.108,.0336){4}{\line(1,0){.108}}
\multiput(70.86,42.36)(.08477,.03168){5}{\line(1,0){.08477}}
\multiput(71.29,42.52)(.06905,.03031){6}{\line(1,0){.06905}}
\multiput(71.7,42.7)(.05764,.02925){7}{\line(1,0){.05764}}
\multiput(72.11,42.9)(.05591,.03244){7}{\line(1,0){.05591}}
\multiput(72.5,43.13)(.04726,.03108){8}{\line(1,0){.04726}}
\multiput(72.88,43.38)(.04544,.03367){8}{\line(1,0){.04544}}
\multiput(73.24,43.65)(.03865,.03215){9}{\line(1,0){.03865}}
\multiput(73.59,43.94)(.03311,.03083){10}{\line(1,0){.03311}}
\multiput(73.92,44.24)(.03134,.03264){10}{\line(0,1){.03264}}
\multiput(74.23,44.57)(.03273,.03816){9}{\line(0,1){.03816}}
\multiput(74.53,44.91)(.03054,.03993){9}{\line(0,1){.03993}}
\multiput(74.8,45.27)(.03179,.04678){8}{\line(0,1){.04678}}
\multiput(75.06,45.65)(.03328,.05541){7}{\line(0,1){.05541}}
\multiput(75.29,46.04)(.03013,.05719){7}{\line(0,1){.05719}}
\multiput(75.5,46.44)(.03136,.06858){6}{\line(0,1){.06858}}
\multiput(75.69,46.85)(.03296,.08428){5}{\line(0,1){.08428}}
\multiput(75.85,47.27)(.02819,.08599){5}{\line(0,1){.08599}}
\put(75.99,47.7){\line(0,1){.437}}
\put(76.11,48.14){\line(0,1){.443}}
\put(76.2,48.58){\line(0,1){.447}}
\put(76.27,49.03){\line(0,1){.451}}
\put(76.31,49.48){\line(0,1){.523}}
\put(68,50){\makebox(0,0)[cc]{$\dots$}}
\end{picture}

\begin{center}
Fig.~1
\end{center}

Let $\gamma_i$ be the closed path starting at $x$ and passing the $i$-th circle
in the clockwise direction. The infinite concatenation
of these paths $\gamma_1\gamma_2\dots $ determines an element of the fundamental group $\pi_1(\mathcal{H},x)$.
Actually any (infinite) concatenation of these paths determines an element of $\pi_1(\mathcal{H},x)$ as soon as
the number of occurrences of every $\gamma_i^{\pm 1}$ is finite.

\medskip

{\it Remark.} There is a general construction, which allows to consider the fundamental groups of topological spaces
as topological groups. Let $X$ be a path connected topological space and $x$ be a point of $X$.
For every path-connected neighborhood $U$ of $x$, we consider the homomorphism $\psi_U:\pi_1(U,x)\rightarrow \pi_1(X,x)$ induced by the embedding $U\hookrightarrow X$. As a basis of neighborhoods of $1$ in $\pi_1(X,x)$ we take the set of normal subgroups $\langle\!\langle {\text{\rm im}}(\psi_U)\rangle\!\rangle$, where $U$ runs over all path-connected neighborhoods of $X$ containing $x$. We call this topology on $\pi_1(X,x)$ {\it canonical}.

\medskip

In~\cite[Theorem 4.1]{MM}, Morgan and Morrison proved that the fundamental group of the Hawaiian Earring, $\pi_1(\mathcal{H},x)$, is canonically isomorphic to BF$(\mathbb{N})$.
Moreover, these groups are isomorphic as topological groups, where $\pi_1(\mathcal{H},x)$
is endowed by the canonical topology and
BF$(\mathbb{N})$ is endowed by the natural topology explained in Section~\ref{naturaltopology}.

The group $\pi_1(\mathcal{H},x)$ is not free.
This follows from~\cite[Theorem~4.1]{MM} and a remark in~\cite[page 80]{H}; see also a short proof in~\cite{S}.
By the above statement 1), this implies that $\BFL$ is not free for any infinite $\varLambda$.

\medskip

{\it Remark.} The word ``free'' for the big free group $\BFL$ is explained by the uniqueness of reduced representatives of elements of $\BFL$. The word ``big'' is explained by the fact $|\BFL|=2^{|\varLambda|}$ for infinite $\varLambda$.

\subsection{Admissible sets and generalized free groups}\label{Sect2.2}

The subgroups of $\BFL$ can be described by specifying certain subsets of $\T$ which we define next:

\medskip
\emph{Definition.} A subset $\mathcal{S}$ of $\T$ is called \emph{admissible} if it has the following properties:

\begin{itemize}
 \item $f, f' \in \mathcal{S} \Rightarrow f f' \in \mathcal{S}$.
 \item $f \in \mathcal{S} \Rightarrow \bar{f} \in \mathcal{S}$.\\
 (For a totally ordered set $S$, let $\bar{S}$ be the set $S$ with the reverse ordering. Given $f:S\rightarrow \Lambda$, we define
 $\bar{f}:\bar{S}\rightarrow \varLambda$ by $\bar{f}(s):=f(s)^{-1}$.)
 \item For $f \in \mathcal{S}$ and $g \in \T$ with $f\searrow g$, we have $g \in \mathcal{S}$.
\end{itemize}

Given an admissible subset $\mathcal{S}$ of $\T$, we obtain a group $\BFS:=\mathcal{S}/\approx$ where $\approx$
is the equivalence relation on $\mathcal{S}$ generated by $f \searrow g$.
Since every element of $\BFL$ admits a unique reduced representing element in $\T$,
it follows easily that $\BFS$ is a subgroup of $\BFL$ if $\mathcal{S}$ is admissible. Conversely, any subgroup of $\BFL$
is of the form $\BFS$ with an admissible set $\mathcal{S}$.

\smallskip

For certain questions concerning infinite groups, big free groups are a more appropriate tool than free groups. They present, however,
new problems: Whereas subgroups of free groups are again free groups, subgroups of big free groups need not be big free groups.

Indeed, the subgroup $\mathcal{F}(\mathbb{N})$ of BF$(\mathbb{N})$ is not a big free group: On the one hand,
it is not isomorphic to $\BFL=\mathcal{F}(\varLambda)$ for a finite set $\varLambda$; on the other hand,
it is not isomorphic to $\BFL$ for an infinite set $\varLambda$ since then $\BFL$ is uncountable.


\medskip
\emph{Definition.} Any group $G$ with $\FL\leqslant G \leqslant \BFL$ is called a \emph{generalized free group} over $\varLambda$.

\medskip
 For an admissible subset $\mathcal{S}$ of $\T$, the group
$\BFS$ is a generalized free group over $\varLambda$ if and only if $\varLambda \subseteq \mathcal{S}$.

\subsection{Big free groups as topological groups}\label{naturaltopology}

The groups $\BFL$, and hence also their subgroups $\BFS$, carry the structure of topological groups. This is what we must explain next.

\smallskip
For a subset $A$ of $\varLambda$ with $A^{-1}=A$ there is a map
\[ \varPhi_A :\T \rightarrow \mathcal{T}(A)\]
which sends $f:S \rightarrow \varLambda$ to its restriction to $f^{-1}(A)$. The map $\varPhi_A$ induces a homomorphism
\[ \varphi_A:\BFL\rightarrow \text{BF}(A)\; .\]

with $\ker \varphi_A=\langle\!\langle\, {\text{\rm BF}}(\Lambda\smallsetminus A)\rangle\!\rangle$.

Since $\text{BF}(A)=\mathcal{F}(A)$ for finite sets $A$, we obtain a homomorphism
\[\varphi:\BFL \longrightarrow \lim_{\stackrel{\longleftarrow}{A \;\text{finite}}}\mathcal{F}(A)\]

which is injective by~\cite[Theorem 3.10]{CC1}. We get a topology on $\BFL$ which has the subgroups $\ker \varphi_A$ for finite subsets $A$ of $\varLambda$
as a basis of the neighborhoods of the neutral element.

\smallskip
Let us call this topology on $\BFL$ and on its subgroups $\BFS$ the \emph{natural topology}. In the natural topology,
the free group $\FL$ is dense in any generalized free group over $\varLambda$.



Unfortunately, the natural topology is too coarse for many purposes. Depending on the situation, we have to consider topologies
which belong to the following class:

\medskip
\emph{Definition.} Given an admissible subset $\mathcal{S}$ of $\T$ with $\varLambda\subseteq \mathcal{S}$, a topology $\mathfrak{T}$ on $\BFS$ is called
\emph{admissible} if it has the following three properties:

\begin{itemize}
 \item With $\mathfrak{T}$, the group $\BFS$ becomes a topological group.
 \item The topology $\mathfrak{T}$ is finer than (i.e. contains at least as many open sets as) the natural topology.
 \item The free group $\FL$ is dense in $\BFS$ with respect to $\mathfrak{T}$.
\end{itemize}

\subsection{Infinite products in big free groups}\label{infinite products}

In the groups $\BFL$ and, more generally, $\BFS$, one can form certain infinite products: Suppose that $T$ is a totally ordered
set and that for every $t \in T$ there is given an element $x_t \in \BFL$. Let $f_t$ be the reduced representative of $x_t$,
and suppose that, for each $\lambda \in \varLambda$, the sets $f^{-1}_t(\lambda)$ are empty for all but finitely many values
of $t$. Then we can form, in an obvious manner, the element
\[\prod_{t \in T} f_t \in \T \]

and can define the infinite product
\[\prod_{t \in T} x_t:= \big[\prod_{t \in T} f_t \big] \in \BFL\; .\]



As a special case of these infinite products, consider any element $f:S\rightarrow \varLambda$ of $\mathcal{T}(\varLambda)$. We then have the elements $f(s)\in \varLambda\subseteq \BFL$ and obtain
$$[f]=\underset{s\in S}{\prod}f(s)\in \BFL\,.$$

In Section 4 we will use this ``word'' notation for the elements of $\BFL$.

\subsection{Generalized presentations}\label{Generalized presentations}

Now we come to the main definitions of the present paper.

\medskip


\begin{defn}\label{generalized generators}
{\rm Let $G$ be a group and $A$ a subset of $G$. We say that $G$ is \emph{generated by}
$A$ \emph{in the generalized sense}
if there exist
 a set $\varLambda$ with a free involution,
an admissible subset $\mathcal{S}$ of $\T$ with $\varLambda \subseteq \mathcal{S}$, and
an epimorphism $p:\BFS \twoheadrightarrow G$ with $p(\varLambda)= A\cup A^{-1}$.
}
\end{defn}

Clearly, if a group $G$ is generated by a subset $A$, it is also generated by $A$ in the generalized sense.
The converse is not always true: the group ${\text{\rm BF}}(\varLambda)$ is generated  by $\varLambda$ in the generalized sense, but it
is not generated by $\varLambda$ in the usual sense if $\varLambda$ is infinite. Indeed, if $\varLambda$ is infinite,
then $|\langle \varLambda\rangle|$=$|\varLambda|$ and $|{\text{\rm BF}}(\varLambda)|=2^{|\varLambda|}$.

\begin{prop}\label{generators} If a group $G$ is generated by a subset $A$ in the generalized sense,
then this holds also for any subset $A_1$ of $G$ containing $A$.
\end{prop}

We leave the proof to the reader as an exercise.


\begin{defn}\label{gp}
{\rm Let $G$ be a group. A \emph{generalized presentation} of $G$ is a tuple
$(\varLambda, \mathcal{S},\mathfrak{T},R)$ with the following properties:
\begin{itemize}
 \item $\varLambda$ is a set with a free involution.
 \item $\mathcal{S}$ is an admissible subset of $\T$ with $\varLambda \subseteq \mathcal{S}$.
 \item $\mathfrak{T}$ is an admissible topology on $\BFS$.
 \item $R$ is a subset of the generalized free group $\BFS$.
 \item There is an epimorphism $p: \BFS \twoheadrightarrow G$ such that $\ker p$ is the smallest normal subgroup which contains the set $R$
 and is closed with respect to $\mathfrak{T}$.
With other words, $\ker p=\overline{\langle\!\langle R\rangle\!\rangle}^{\frak{T}}$.
\end{itemize}

\medskip The sets $\varLambda$ and $R$ are called, respectively, the sets of {\it generalized generators} and {\it defining generalized relations} for this generalized presentation of $G$.

}
\end{defn}

\begin{rmk}
{\rm
a)
If $(\varLambda, \mathcal{S},\mathfrak{T},R)$  is a generalized presentation of a group $G$, then $G$ is
generated by $p(\varLambda)$ in the generalized sense, where $p$ is a map as in~\ref{gp}.

b) Every usual presentation of $G$ gives rise to a generalized one. Indeed, let
$(\varDelta,R)$ be a usual presentation of $G$, i.e. there exists an epimorphism $p: F(\varDelta) \twoheadrightarrow G$ such that $\ker p=\langle\!\langle R\rangle\!\rangle$.
We set $\varLambda=\varDelta^{\pm}$ and define the admissible subset $\mathcal{S}\subseteq \T$ consisting of all maps $f:S\rightarrow \varLambda$ with finite $S$. Then $(\varLambda, \mathcal{S},\mathfrak{T},R)$ is a generalized presentation of $G$, where $\mathfrak{T}$ is the discrete topology on $\BFS=\mathcal{F}(\varLambda)$.

c) Every big free group $\BFS$ has the generalized presentation $(\varLambda, \T,\frak{T}, \emptyset)$,
where $\frak{T}$ is an arbitrary admissible topology on $\BFS$.
}
\end{rmk}


\begin{defn}
{\rm Given a group $G$ and cardinal number $c$, we say that $G$ \emph{admits a generalized presentation of type} $c$ if there
is a generalized presentation
$(\varLambda, \mathcal{S},\mathfrak{T}, R)$ of $G$ with $|\varLambda| \leq c$ and $|R|\leq c$.}
\end{defn}

\begin{ex}
{\rm The additive group $\mathbb{R}$ is generated in the generalized sense by the subset $\{10^{-n}\mid n\in \mathbb{N}\}$. To see this, we define the corresponding $\varLambda$, the free involution on $\varLambda$ and an admissible subset $\mathcal{S}$ of $\mathcal{T}(\varLambda)$ as follows:
$$\varLambda:=\mathbb{Z}\smallsetminus \{0\}; \,\,n^{-1}:=-n,$$
$$\mathcal{S}=\{ (f:S\rightarrow \varLambda)\in \mathcal{T}(\varLambda)  \mid {\text {\rm there exists}}\hspace*{1mm}
M\in \mathbb{N}\hspace*{1mm} {\text {\rm with}}\hspace*{1mm} |f^{-1}(n)|\leqslant M \hspace*{1mm}{\text {\rm for all}}
\hspace*{1mm} n\in \varLambda
\}.$$

The epimorphism $p$ from Definition~\ref{generalized generators} is defined as the unique continuous homomorphism $p:\BFS \rightarrow \mathbb{R}$ with $p(n)=10^{-n}$ for $n\in \mathbb{N}$.

We leave it to the reader to show that $\mathbb{R}$ admits a generalized presentation of type $\aleph_0$.
}
\end{ex}

\begin{ex}
{\rm The group $\BFS$ from the previous example is, as Zastrow~\cite{Z} showed, a free group with $\aleph_1$ generators.
Hence each group with at most $\aleph_1$ elements is generated in the generalized sense by a countable subset.

}\end{ex}

\section{Generalized presentations for infinite\\ symmetric groups}\label{permut}

For a set $X$, let $\varSigma(X)$ be the symmetric group consisting of all bijections $\sigma:X\rightarrow X$. The element obtained
from $x$ by applying $\sigma$ will be denoted by $x\cdot\sigma$. In the present section, we will find two different generalized
presentations for $\varSigma(X)$. Both of them are of type $|X|$ if $X$ is an infinite set.

\medskip

For $x,y \in X$ with $x \neq y$, let $\tau_{x,y} \in \varSigma(X)$ be the transposition interchanging $x$ and $y$ and leaving
all other elements of $X$ fixed. We will show that $\varSigma(X)$ is, in the generalized sense introduced in Section~\ref{BF and GP}, generated by
these transpositions and that the usual relations between them are actually defining generalized relations.

\subsection{An example}\label{an example}
In $\varSigma(\Z)$, consider the shift $\sigma :n \mapsto n+1$. In a self-explaining way, we can write,
for instance,
\begin{eqnarray}
 \sigma =\ldots \tau_{2,3}\; \tau_{1,2}\; \tau_{0,1}\; \tau_{-1,0} \ldots
\end{eqnarray}
or
\begin{eqnarray}
 \sigma = \tau_{0,1}\; \tau_{0,-1}\; \tau_{-1,2} \;\tau_{-1,-2} \;\tau_{-2,3} \;\tau_{-2,-3} \ldots \; .
\end{eqnarray}

There is an important difference between these two ways of writing $\sigma$ as an infinite product: Consider e.g. the
``subword''
\[ \tau_{1,2}\; \tau_{0,1}\; \tau_{-1,0} \ldots\]
of (4). This does not represent an element of $\varSigma(\Z)$: Indeed, if it would represent $\rho \in \varSigma(\Z)$, what would
$\rho(2)$ be?

\smallskip
On the other hand, each subword (i.e. each finite or infinite string of consecutive letters) of (5) defines an element of $\varSigma(\Z)$.

\smallskip
In the first of our two generalized presentations of $\varSigma(\Z)$, both (4) and (5) will be legal ways of writing $\sigma$.
In the second presentation, (4) will be illegal, but (5) will remain legal.

\subsection{The admissible set $\mathcal{S}$}\label{admis set}

We return to the general case and will now describe the set $\varLambda$ and the admissible set $\mathcal{S}$ of a generalized
presentation $(\varLambda, \mathcal{S}, \frak{T},R)$ of $\varSigma(X)$.
Let
\[\varLambda:=\varLambda(X) := \{T_{x,y}\mid x,y\in X,\, x \neq y \}\; ;\]

the free involution on $\varLambda$ sends $T_{x,y}$ to $T_{y,x}$.

\medskip
The definition of $\mathcal{S}$ requires a certain amount of notation.


\smallskip
With every map $f:S \rightarrow \varLambda$, where $S$ is a totally ordered set,
we associate two maps $f_1:S\rightarrow X$
and $f_2:S\rightarrow X$ by the following rule: if $s\in S$
and $f(s)=T_{x,y}$, we set $f_1(s)=x$ and $f_2(s)=y$. For $z\in X$ we set

\[U(z, f):=\{s \in S \mid z\in \{f_1(s),f_2(s)\}\}\; .\]
Let us write
\[\mathcal{S}_0:= \{(f:S \rightarrow \varLambda) \mid  \;|U(x, f)|< \infty\; \forall\; x \in X\}\subseteq \mathcal{T}
 (\varLambda)\; .\]

We will describe a subset $\mathcal{S}$ of $\mathcal{S}_0$.

Let $f:S\rightarrow \varLambda$ be an element of $\mathcal{S}_0$ which is fixed for the moment; we have to define
what it means that $f$ belongs to $\mathcal{S}$. For each $x \in X$, we will define inductively four sequences
\begin{eqnarray}
 x_0^+ ,&  x_1^+ ,\qquad x_2^+ ,& \ldots \\[1ex]
 &  s_1^+(x) ,\quad s_2^+(x) ,&\ldots \\[1ex]
 x_0^-,& x_1^-,\qquad x_2^-,&\ldots \\[1ex]
  &s_1^-(x),\quad s_2^-(x),&\ldots
\end{eqnarray}

The sequences (6) and (8) will consist of elements of $X$, the sequences (7) and (9) of elements of $S$.

\medskip
These sequences may be
finite or infinite. The sequence (6) will be finite iff (7) is finite, and if this is the case, both will end with the term with
the same index, say $n^+(x)$. Here, we allow that $n^+(x)=0$; this is to mean that (6) is the 1-term sequence $x^+_0$ and (7) is the
empty sequence. If (6) and (7) are infinite sequences, let us put $n^+(x):=\infty$. So, for any $x \in X$, we will have
\[n^+(x) \in \N \cup \{0, \infty\}\; .\]

Similarly, for the sequences (8) and (9); so we will obtain also
\[n^-(x) \in \N \cup \{0, \infty\}\; .\]

Now we come to the actual definition of the four sequences:
\[x_0^+:=x_0^-:=x \; .\]

If $U(x,f)= \emptyset$, let $n^+(x):=n^-(x)=0$.

\smallskip
If $U(x,f)\neq \emptyset$, let
\begin{eqnarray*}
 s_1^+(x):= \min U(x,f)\; ,\\[1ex]
 s_1^-(x):= \max U(x,f)\; .
\end{eqnarray*}

Now suppose inductively that for some $n \in \N$, we have already defined
\[x_0^+,\ldots , x_{n-1}^+,\, s_1^+(x),\ldots , s_n^+(x)\]
and that $x_{k-1}^+ \in \{f_1(s_k^+(x)), f_2(s_k^+(x)) \}$ for all $k=1,\dots ,n$. Then we define $x_n^+$ by requiring
\[\{x_{n-1}^+,x_n^+\}=\{f_1(s_n^+(x)), f_2(s_n^+(x)) \}\; .\]

Let $s^+_{n+1}(x)$ be the smallest element of $S$ which is contained in $U(x^+_n, f)$ and is bigger than $s^+_n(x)$,
assuming that such an element exists. If there is no such element, let $n^+(x):=n$.

\medskip
This completes the definition of the sequences (6) and (7), and it should be obvious how, by symmetry, the sequences (8) and (9)
are defined. Observe that
\begin{eqnarray*}
 s_1^+(x) < s_2^+(x) < \ldots , \\[1ex]
 s_1^-(x) > s_2^-(x) > \ldots \, .
\end{eqnarray*}

Now we define the subset $\mathcal{S}$ of $\mathcal{S}_0$ by declaring that $f$ belongs to $\mathcal{S}$ if and only if, for all
$x \in X$, we have
\[n^+(x)< \infty \quad , \quad n^-(x) < \infty \,.\]

\begin{lemma}\label{Sigma1}
 The subset $\mathcal{S}$ of $\T$ is admissible.
\end{lemma}

\medskip {\it Proof.}
 We have to show that an element of $\T$ is contained in $\mathcal{S}$ if it is obtained from an element of $\mathcal{S}$ by cancellation.
 This amounts to the following: Suppose we are given an element $(f:S\rightarrow \varLambda) \in \mathcal{S}$ and a subset
 $T$ of $S$ with an involution $\ast$ such that, for all $t \in T$, we have, in addition to the conditions (2) and (3) of Section~\ref{BF and GP},
 that

 \parbox{11cm}{\begin{eqnarray*}
 f(t^{\ast})=T_{x,y}\Leftrightarrow f(t)=T_{y,x}\, . \end{eqnarray*}} \hfill
 \parbox{8mm}{$(1')$}

 \medskip
 Then we have to show that $g:=f|S \smallsetminus T$ belongs to $\mathcal{S}$. Let $x \in X$. As in the definition of $\mathcal{S}$,
 we have the four finite sequences $(x_n^{\pm})$ and $(s_n^{\pm}(x))$ associated with $f$. We have to consider the corresponding
 sequences associated with $g$. We denote them by $(\xi_n^{\pm})$ and $(\sigma_n^{\pm}(x))$ with $\xi_n^{\pm} \in X$ and
 $\sigma_n^{\pm}(x) \in S \smallsetminus T$. We have to show that they are finite sequences. We will show that  $(\xi_n^+)$ is a
 subsequence of $(x_n^+)$ and $(\sigma_n^+(x))$ is a subsequence of $(s_n^+(x))$. To abbreviate, let us write $s_n:=s_n^+(x)$
 and $\sigma_n:=\sigma_n^+(x)$. We have $\xi_0^+=x=x_0^+$. We will show that $\sigma_1$ is a term in the sequence $(s_n)$ and that
 $\xi_1^+$ is the corresponding term in the sequence $(x_n^+)$.
 We have

 \begin{eqnarray}
  s_1=\min \{s \in S \mid x\in \{f_1(s),f_2(s)\}\}\; ,\quad \\[1ex]
  \sigma_1= \min \{s \in S \smallsetminus T \mid x\in \{f_1(s),f_2(s)\}\} \; .
 \end{eqnarray}

If $s_1 \in S \smallsetminus T$, then obviously $\sigma_1=s_1$ and $\xi_1=x_1$. Therefore we can assume that $s_1 \in T$. Then we
conclude from $(1')$ and (2) that
\[s_1< s_1^ {\ast} <\sigma_1 \; .\]

By (3) and the definition of $s_2$, there are two possibilities:

\smallskip
Either $s_2=s_1^{\ast}$ or
\[s_1 < s_2 < s_2^{\ast} < s_1^{\ast} \; .\]

Iterating this argument, we see that there is a number $m$ such that the finite sequence $(s_n)$ begins with the terms
\[s_1 < s_2 < \ldots < s_m < s_m^{\ast} < \ldots < s_2^{\ast} < s_1^{\ast} \]
which all lie in $T$. From $(1')$, we conclude that the sequence $(x_n^+)$ begins with the terms
\[x_0^+,x_1^+, \ldots , x_{m-1}^+, x_m^+, x_{m-1}^+, \ldots , x_1^+, x_0^+ \; .\]

\vspace*{10mm}

\unitlength 1mm 
\linethickness{0.4pt}
\ifx\plotpoint\undefined\newsavebox{\plotpoint}\fi 
\begin{picture}(110,40)(-20,40)
\thicklines
\put(10,60){\line(1,0){80}}
\put(20,60){\circle*{2.5}}
\put(40,60){\circle*{2.5}}
\put(60,60){\circle*{2.5}}
\put(80,60){\circle*{2.5}}
\qbezier(40,60)(50,81.5)(60,60)
\qbezier(20,60)(50,101.625)(80,59.75)
\put(19,65){\makebox(0,0)[cc]{$s_1$}}
\put(39,65){\makebox(0,0)[cc]{$s_2$}}
\put(62,65){\makebox(0,0)[cc]{$s_2^{\ast}$}}
\put(82,65){\makebox(0,0)[cc]{$s_1^{\ast}$}}
\put(19.75,54.25){\makebox(0,0)[cc]{$T_{x_0^{+}x_1^{+}}$}}
\put(39.75,54.25){\makebox(0,0)[cc]{$T_{x_1^{+}x_2^{+}}$}}
\put(59.75,54.25){\makebox(0,0)[cc]{$T_{x_2^{+}x_1^{+}}$}}
\put(79.75,54.25){\makebox(0,0)[cc]{$T_{x_1^{+}x_0^{+}}$}}
\end{picture}
\vspace*{-10mm}
\begin{center} Figure 2. Case $m=2$.
\end{center}

Since $s_{2m}=s_1^{\ast} < \sigma_1$, we see that $n^+(x) > 2m$. If $s_{2m+1} \in S \smallsetminus T$, we have $\sigma_1=s_{2m+1}$ and
$\xi_1^+=x_{2m+1}^+$. If $s_{2m+1} \in T$, we can repeat the above argument; we find that the sequence $(s_n)$ begins with terms of
the form

\parbox{11cm}{
\begin{eqnarray*}
 s_{1,1}, \ldots , s_{1,m_1}\,,\, s_{1,m_1}^{\ast}, \ldots , s_{1,1}^{\ast},\qquad\;\;\,\\[1ex]
 s_{2,1}, \ldots , s_{2,m_2}\,,\, s_{2,m_2}^{\ast}, \ldots , s_{2,1}^{\ast}, \ldots , \sigma_1
\end{eqnarray*}}\hfill
\parbox{8mm}{\begin{eqnarray}\end{eqnarray}}

and that the sequence $(x_n^+)$ begins with terms of the form

\parbox{11cm}{
\begin{eqnarray*}
x, x_{1,1}, \ldots , x_{1,m_1-1}\,,\, x_{1,m_1}\,,\, x_{1,m_1-1},\ldots , x_{1,1},x,\qquad\quad\;\,\\[1ex]
 x_{2,1}, \ldots , x_{2,m_2-1}\,,\, x_{2,m_2}\,,\,x_{2,m_2-1}, \ldots , x_{2,1},x, \ldots , \xi_1^+\;.
\end{eqnarray*}}\hfill
\parbox{8mm}{\begin{eqnarray}\end{eqnarray}}

\unitlength 1mm 
\linethickness{0.4pt}
\ifx\plotpoint\undefined\newsavebox{\plotpoint}\fi 
\begin{picture}(110,50)(8,40)
\thicklines
\put(10,60){\line(1,0){140}}
\put(20,60){\circle*{2.5}}
\put(40,60){\circle*{2.5}}
\put(60,60){\circle*{2.5}}
\put(80,60){\circle*{2.5}}
\put(100,60){\circle*{2.5}}
\put(120,60){\circle*{2.5}}
\put(140,60){\circle*{2.5}}
\qbezier(40,60)(50,81.5)(60,60)
\qbezier(20,60)(50,101.625)(80,59.75)
\qbezier(100,60)(110,81.5)(120,60)
\put(19,65){\makebox(0,0)[cc]{$s_{1,1}$}}
\put(39,65){\makebox(0,0)[cc]{$s_{1,2}$}}
\put(62,65){\makebox(0,0)[cc]{$s_{1,2}^{\ast}$}}
\put(82,65){\makebox(0,0)[cc]{$s_{1,1}^{\ast}$}}
\put(99,65){\makebox(0,0)[cc]{$s_{2,1}$}}
\put(141,65){\makebox(0,0)[cc]{$\sigma_1$}}
\put(122,65){\makebox(0,0)[cc]{$s_{2,1}^{\ast}$}}
\put(19.75,54.25){\makebox(0,0)[cc]{$T_{x,x_{1,1}}$}}
\put(39.75,54.25){\makebox(0,0)[cc]{$T_{x_{1,1},x_{1,2}}$}}
\put(59.75,54.25){\makebox(0,0)[cc]{$T_{x_{1,2},x_{1,1}}$}}
\put(79.75,54.25){\makebox(0,0)[cc]{$T_{x_{1,1},x}$}}
\put(99.75,54.25){\makebox(0,0)[cc]{$T_{x,x_{2,1}}$}}
\put(119.75,54.25){\makebox(0,0)[cc]{$T_{x_{2,1},x}$}}
\put(140.75,54.25){\makebox(0,0)[cc]{$\xi_1^{+}$}}
\end{picture}

\vspace*{-10mm}
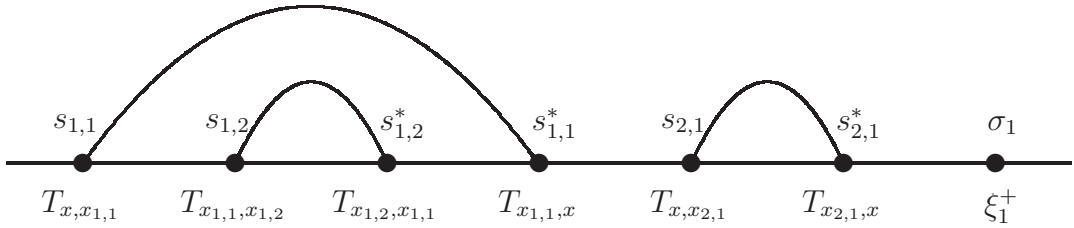
\begin{center} Figure 3. Case $m_1=2$, $m_2=1$.
\end{center}

Furthermore, the part of the sequence (12) lying between $\sigma_1$ and $\sigma_2$ has the same form as the part preceding $\sigma_1$, and so on.
$\Box$

\medskip
For later use, we state a fact which is clear from the proof of Lemma~\ref{Sigma1}.

\begin{lemma}\label{Sigma2}
 Given $f \in \mathcal{S}$ and $x \in X$, denote by $x_{\infty}(f)$ the last element in the finite sequence
 $x_0^+, x_1^+, \ldots \;$. If $f \searrow g$, then $x_{\infty}(f)=x_{\infty}(g)$. \hfill $\Box$
\end{lemma}

\subsection{The admissible topology $\mathfrak{T}$ and the homomorphism $p$}
Since $\mathcal{S}$ is an admissible subset of $\T$
containing $\varLambda$, the group $\BFS$ contains the free group $\mathcal F(\varLambda)$. In order to continue
with the description of a generalized presentation of $\varSigma(X)$, we will now define an admissible topology $\mathfrak{T}$ on
$\BFS$ and a homomorphism $p:\BFS \rightarrow \varSigma(X)$.

\medskip
On $\varSigma(X)$, there is a natural topology making $\varSigma(X)$ into a topological group. A basis for the neighborhoods of $1$
consists of the subgroups
\[ U_C:=\{\sigma \in \varSigma(X) \mid c \cdot \sigma=c \;\forall \; c \in C\},\]

where $C$ goes through the set of finite subsets of $X$.

\medskip
Using the notation explained in Section~\ref{naturaltopology}, we have for each finite subset $C$ of $X$ the homomorphism
\[\psi_C:=\varphi_{\varLambda(C)}\mid_{\BFS} :\BFS \rightarrow \text{BF}(\varLambda(C)) = \mathcal F(\varLambda(C))\hookrightarrow
\mathcal F(\varLambda)\; .\]
The kernels of these homomorphisms form a basis of neighborhoods of 1 for the natural topology on $\BFS$. We define a finer
topology $\mathfrak{T}$ by requiring that a basis of neighborhoods of 1 is given by the subgroups
\[ W_C:= \ker \psi_C \cap V_C\]
where
\[V_C:=\{ [\gamma] \in \BFS \mid x_{\infty}(\gamma)=x \quad \forall\; x \in C \}\; .\]

Here $x_{\infty}(\gamma)$ is the element introduced in Lemma~\ref{Sigma2}. It is easy to verify that $\mathfrak{T}$ is an admissible topology.

\begin{lemma}\label{Sigma3}
If we endow ${\text {\rm BF}}(\varLambda,\mathcal{S})$ with the topology $\mathfrak{T}$,
there is a unique continuous homomorphism $p:{\text {\rm BF}}(\varLambda,\mathcal{S}) \rightarrow \varSigma(X)$ with $p(T_{x,y})=\tau_{x,y}$ for all $x,y \in X, x\neq y$.

\smallskip
For $F=[f] \in {\text {\rm BF}}(\varLambda,\mathcal{S})$ with $f \in \mathcal{S}$ and $x \in X$, the element $x \cdot p(F) \in X$ is given by
\[ x \cdot p(F)= x_{\infty}(f)\; .\]
Equivalently, we can describe $p$ as follows:
There is a unique homomorphism $p:\mathcal F(\varLambda) \rightarrow \varSigma(X)$ sending $T_{x,y}$ to $ \tau_{x,y}$. We have to extend
 $p$ to all of ${\text {\rm BF}}(\varLambda,\mathcal{S})$. Given $F, f$ and $x$ as above, consider the finite sequence $(x^+_n)$ assigned to $x$ and
 $f$.
 Let $C$ be a finite subset of $X$ containing the elements $x^+_n$.
 Then $\psi_C(F) \in \mathcal F (\varLambda)$ is a finite word in $T_{y,z}$ with $y,z\in C$, and we have

 \[ x \cdot p(F)= x \cdot p(\psi_C(F)) \; .\]

 \hfill $\Box$

 \end{lemma}

\begin{lemma}\label{Sigma4}
 The homomorphism $p=p_X :{\text {\rm BF}}(\varLambda,\mathcal{S}) \rightarrow \varSigma(X)$ is surjective.
\end{lemma}

\medskip {\it Proof.}
 Let $\sigma \in \varSigma(X)$. Suppose that $X$ is the disjoint union of subsets $X_{\beta}$ such that $\sigma(X_{\beta})=X_{\beta}$
 for all $\beta$, and denote by $\sigma_{\beta} \in \varSigma(X_{\beta})$ the restriction of $\sigma$ to $X_{\beta}$. If $\sigma_{\beta} \in
 \im\, p_{X_{\beta}}$ for all $\beta$, then, obviously, $\sigma \in \im\, p_X$.

 \smallskip
 Therefore, to show that $\sigma \in \im\, p_X$, it suffices to assume that the group generated by $\sigma$ acts transitively on $X$.
 Then we are either in the trivial situation that $X$ is finite, or we are in the situation of the Example in Section~\ref{an example}, which is also obvious.
\hfill $\Box$

\subsection{A generalized presentation of $\varSigma(X)$}
Now we complete the description of the generalized presentation $(\varLambda,\mathcal{S},\mathfrak T,R)$ of $\varSigma(X)$.

\smallskip

Let $R$ be the subset of $\BFS$ consisting of the elements
\begin{itemize}
 \item $T^2_{x,y}$,
 \item $[T_{x,y},T_{z,w}]$ for $ |\{x,y,z,w\}|=4$,
 \item $T_{x,y}\, T_{x,z}\, T^{-1}_{x,y}\, T^{-1}_{y,z}$ for $|\{x,y,z\}| =3$.
\end{itemize}

For finite $X$, $\langle \varLambda\mid R\rangle$ is a presentation of $\varSigma(X)$. Indeed, this presentation can be easily obtained from the classical one (see~\cite[Theorem~7.1 in Chapter 2]{B}) with the help of Tietze transformations.

Let $N$ be the smallest closed (with respect to $\mathfrak T$) normal subgroup of $\BFS$ containing $R$.

\begin{lemma}\label{Sigma5}
 $\ker p =N$.
\end{lemma}

\medskip {\it Proof.} Clearly $N\subseteq \ker p$. We show that $\ker p\subseteq N$.
Given $g \in \ker p$, we have to show that every neighborhood $gW_B$ of $g$, where $W_B=\ker\psi_B \cap V_B$, contains an element of $\langle\!\langle R \rangle\!\rangle$.
It suffices to show that for each finite subset $B$ of $X$ there is an element
$h \in \langle\!\langle R \rangle\!\rangle$ with $g^{-1}h\in {\text{\rm ker}}\,\psi_B$.

(Indeed, we would have then $h\in \langle\!\langle R \rangle\!\rangle\subseteq \ker p\subseteq V_B$. Since
$g \in \ker p$, we obtain $g^{-1}h\in V_B$. This and $g^{-1}h\in {\text{\rm ker}}\,\psi_B$ would imply
$g^{-1}h\in W_B$, i.e. $h\in gW_B$.)

 \smallskip
 Let $g=[\gamma]$ with $\gamma \in \mathcal{S}$. There is a finite subset $C$ of $X$ which
 contains, for every $x \in B$, all elements $x^+_n$, formed with respect to $\gamma$. Then $p \circ \psi_C(g)$ is a permutation
 of $X$ which is the identity on $X \smallsetminus C$ and on $B$ since $g\in \ker p$.
 Hence there exists an element $a \in \mathcal F(\varLambda(C \smallsetminus B))$
 such that $p \circ \psi_C(g)=p(a)$, that is, $p \circ \psi_C(a^{-1}g)=1$. Since $C$ is finite, we conclude that
 \[h:= \psi_C(a^{-1}g)\in \langle\!\langle R \rangle\!\rangle\; ,\]

 hence $\psi_C(h^{-1}a^{-1}g)=1$, and therefore
 \[\psi_B(h^{-1}g)= \psi_B(h^{-1}a^{-1}g)=1 \; .\]

\hfill$\Box$

\medskip

The last two lemmas show

\begin{theorem}\label{ThSigma1} For any set $X$, the triple
 $(\varLambda,\mathcal{S},\mathfrak T, R)$ is a generalized presentation of $\varSigma(X)$. In particular, if $X$ is infinite, the group $\varSigma(X)$
 admits a generalized presentation of type $|X|$. \hfill $\Box$
\end{theorem}

\subsection{Another generalized presentation of $\varSigma(X)$}

As already indicated in Section~\ref{an example}, there is a second
generalized
presentation of $\varSigma(X)$ with an admissible set
$\mathcal{S}'$ which is smaller than $\mathcal{S}$.

\medskip
We continue to consider the admissible subset $\mathcal{S}$ of $\T$ introduced in Section~\ref{admis set}. Let $\mathcal{S}'$ be the subset
of $\T$ consisting of all maps $f:S\rightarrow \varLambda$ which satisfy the following condition:

\smallskip
{\it For each interval $I$ of $S$, we have $f \mid I \in \mathcal{S}$.}

\medskip
Note that the second representation of the shift $\sigma$ on $\Z$ in Example~\ref{an example} belongs to $\mathcal{S}'$.
We have the following lemma and theorem.

\begin{lemma}\label{Sigma6}
 $\mathcal{S}'$ is an admissible subset of $\mathcal{T}(\varLambda)$ and the relative topology $\mathfrak{T}'$ on ${\text{\rm BF}} (\varLambda,\mathcal{S}')$
 defined by $\mathfrak T$ is an admissible topology. \hfill $\Box$
\end{lemma}

\begin{theorem}\label{ThSigma2}
 $(\varLambda,\mathcal{S}',\mathfrak{T}', R)$ is a generalized presentation of $\varSigma(X)$. \hfill $\Box$
\end{theorem}

\section{A generalized presentation of $\Aut(F_{\omega})$}\label{GpAut}

We consider a set $X$ and the free group $F(X)$ with basis $X$. We introduce the set $X^\pm =X \cup X^{-1}$ with the free involution
$x \mapsto x^{-1}$. With the notation introduced in Section~\ref{Bfg}, we have $F(X) = \mathcal{F}(X^\pm)$. We will study the automorphism group
$\Aut \,F(X)$.

\smallskip
For $w \in F(X)$ and $\varphi \in \text{End}\, F(X)$, we write $w\varphi$ for the element obtained from $w$ by applying $\varphi$.
For $x,y \in X^\pm$ with $x \neq y^{\pm 1}$, let $E_{xy} \in \Aut(F(X))$ be the automorphism which sends $x$ to $xy$ and keeps all elements of $X^\pm \smallsetminus \{x,x^{-1}\}$ fixed.
Such automorphisms are called {\it elementary Nielsen of the first kind}. The automorphisms which
map $X^{\pm}$ onto itself are called {\it elementary Nielsen of the second kind} (or {\it monomial}).
Let $\mathcal{E}(X)$ be the set of all elementary Nielsen of the first kind and let $\mathcal{M}(X)$ be the group consisting of all monomial automorphisms of $F(X)$.

Nielsen showed that, for any nonempty finite set $X$, the elementary Nielsen automorphisms of the first kind generate a subgroup of $\Aut(F(X))$ of index 2.
He also gave a finite presentation of $\Aut(F(X))$.

We will show that, for a countably infinite set $X$, the group
$\Aut(F(X))$ is,
in our generalized sense, generated by the $E_{xy}$ and that the usual ``finite'' relations and some ``infinite''
relations are defining generalized relations for $\Aut(F(X))$.
In particular, we will find a generalized presentation $(\mathcal E,\mathcal S, \mathfrak T,R)$ for $\Aut(F(X))$ of type
$\aleph_0$ in the case where $X$ is a countably infinite set.

In Section~\ref{Appendix B} we will show that
in this case every proper normal subgroup of ${\text{\rm Aut}}(F(X))$ has index $2^{\aleph_0}$.

\smallskip
For the first steps of our argument, $X$ is allowed to be an arbitrary set.

\subsection{The admissible set $\mathcal S$}\label{AsS2}

The set $\mathcal{E}=\mathcal{E}(X)$ will play the role of the set called $\varLambda$ in the previous two sections.
We define the free involution $^{-1}:\mathcal{E}\rightarrow \mathcal{E}$ by the rule: $E_{xy}^{-1}=E_{xy^{-1}}$.
Of course, we have to consider certain infinite products of the $E_{xy}$.

\begin{ex}\rm {
 Let $X =\{x_1, x_2, \ldots \}$ be a countably infinite set.

 \smallskip
 1) The infinite product $E_{x_1 x_2} E_{x_3 x_4} E_{x_5 x_6} \ldots $ can be interpreted as the automorphism of $F(X)$ which
 fixes $x_i$ for even $i$ and sends $x_i$ to $x_i x_{i+1}$ for odd $i$.

 \smallskip
 2) The infinite product $\ldots E_{x_4 x_3} E_{x_3 x_2} E_{x_2 x_1}$ determines an automorphism of $F(X)$ which sends $x_i$
 to $x_i \ldots x_2 x_1$.

 \smallskip
 3) The infinite word $E_{x_1 x_2} E_{x_1 x_3} E_{x_1 x_4} \ldots$ does not determine an endomorphism of $F(X)$ since it would send
 $x_1$ to the infinite word $x_1 (\ldots x_4 x_3 x_2)$.}
\end{ex}

Given a totally ordered set $S$ and a map $f: S \rightarrow \mathcal E$, we obtain two maps $f_1, f_2: S \rightarrow X^{\pm}$ by
\[f(s) = E_{f_1(s)f_2(s)}\, .\]

The last example suggests that for a map $(f: S \rightarrow \mathcal E) \in \mathcal T(\mathcal E)$ to be admissible we should
at least require that $f$ belongs to
\[\mathcal S_0 := \{f \in \mathcal T (\mathcal E) | f_1^{-1}(x) \, \text{is finite for each} \, x \in X^\pm \}\, .\]

\begin{ex} \rm{
The infinite word $E_{x_1 x_2} E_{x_2 x_3} E_{x_3 x_4} \ldots $ is defined by an element of $\mathcal S_0$ but still does not
determine an endomorphism of $F(X)$ since $x_1$ would be sent to the infinite word $x_1 x_2 x_3 \ldots$.}
 \end{ex}

This example suggests that we should also require that $f$ belongs to
\[\begin{array}{lcl}\mathcal S_1 := \{(f: S\rightarrow \mathcal E) \in \mathcal T(\mathcal E) &\!\!| &\!\!\text{there is no injective and order preserving map} \\[1ex]
&& \varphi: \N \rightarrow S \,\, \text{with}\,\, f_1(\varphi(n+1))=(f_2(\varphi(n)))^{\pm 1}\}\,.\end{array}
\]

\begin{ex} \rm{
The infinite word $\ldots E_{x_3 x_4} E_{x_2 x_3} E_{x_1 x_2}$ is defined by an element of $\mathcal S_0 \cap \mathcal S_1$ and
it determines the endomorphism $\alpha$ of $F(X)$ which sends $x_i$ to $x_i x_{i+1}$ for all $i$. However, $\alpha$ is not
invertible since its image consists of words of even length.}
 \end{ex}

So we are finally led to the definition

\begin{center}
\fbox{$\mathcal S := \{f \in \mathcal T (\mathcal E) | f \in \mathcal S_0 \, \text{and} \, f, \bar{f} \in \mathcal S_1\}\,.$}
\end{center}

Recall that $\bar{f}$ was defined in Section~\ref{Sect2.2}. The condition $\bar{f} \in \mathcal S_1$ means that there is no injective
and order preserving map $\psi :(- \N) \rightarrow S$ with $f_1(\psi(n-1))=(f_2 (\psi(n)))^{\pm 1}$ for all $n \in (-\N)$.

\medskip
The set $\mathcal S$ will be the (obviously admissible) set used in our generalized presentation of $\Aut( F(X))$.

\subsection{Another description of the set $\mathcal{S}$}\label{Bas}

As a technical device, we have to introduce an admissible subset $\tilde{\mathcal S} \subseteq \mathcal S_0$ which is
very similar to the one used for $\varSigma (X)$. We will then show that $\tilde{\mathcal S}=\mathcal{S}$.

\smallskip
Let us fix an element $f \in \mathcal S_0$. We want to define what it means that $f$ belongs to $\tilde{\mathcal S}$.
Let us also fix for the moment an element $x \in X^\pm $. We define inductively two sequences

\begin{eqnarray}
 A_0 ,& \!\! A_1 ,\; A_2 ,& \ldots \label{S14} \\[1ex]
 &  s_1 ,\; s_2 ,&\ldots \label{S15}
\end{eqnarray}

The sequence (\ref{S14}) consists of subsets of $X^\pm$, the sequence (\ref{S15}) consists of elements of $S$. These sequences may be
finite or infinite. The sequence (\ref{S14}) will be finite iff (\ref{S15}) is finite, and if this is the case, both will end with the term
with the same index, say $\frak{n}^+ =\frak{n}^+(f,x)$. We allow $\frak{n}^+= 0$; this is to mean that (\ref{S14}) is the  1-term sequence $A_0$ and
(\ref{S15}) is the empty sequence. If (\ref{S14}) and (\ref{S15}) are infinite sequences, let us put $\frak{n}^+(f,x):= \infty$.

\smallskip
Now we come to the actual definition of the two sequences:

\[A_0 := \{x, x^{-1} \} \, .\]

If $f_1^{-1}(A_0)=\emptyset$, let $\frak{n}^+(f,x)=0$. Otherwise, let

\[\begin{array}{l}
 s_1:= \text{min} \, f_1^{-1}(A_0)\,,\\[1ex]
 A_1:= A_0 \cup \{f_2(s_1), f_2(s_1)^{-1}\}\,,\\[1ex]
 s_2:= \text{min} \,\{s \in S \,| \,s>s_1 \,\text{and}\, f_1(s) \in A_1\}\\
\end{array}\]

if the latter set is non-empty; otherwise put $\frak{n}^+ = 1$. Observe that the minimum exists since $f \in \mathcal S_0$. Put
\[  A_2 := A_1 \cup \{ f_2(s_2), f_2(s_2)^{-1}\} \qquad\qquad
\]
and so on. Obviously,

\[\tilde{\mathcal S} := \{ f \in \mathcal S_0 \;|\; \frak{n}^+(f,x) < \infty \,\text{ and}\,\, \frak{n}^+(\bar{f},x) < \infty \,\,\text{for all}\,\, x \in X^\pm \}\]

is an admissible subset of $\mathcal T(\mathcal E)$.

\begin{lemma}\label{use it}
 $\mathcal S = \mathcal{\tilde{S}}$.
\end{lemma}

\begin{demo} The inclusion $\mathcal {\tilde S} \subseteq \mathcal{S}$ is obvious. So, we will prove that $\mathcal S \subseteq \mathcal{\tilde{S}}$.
 Let $f \in \mathcal S$ and $ x \in X^\pm$. It suffices to show that $\frak{n}^+(f,x) < \infty$. Consider the corresponding sequence $(s_k)$ in $S$
 and the sequence $(A_k)$ of finite subsets of $X^{\pm}$.
 Let $V_k:= A_k /\sim$ where the equivalence relation $\sim $ is
 given by $a \sim a^{\pm 1}$. Then we have $V_0 \subseteq V_1 \subseteq V_2 \subseteq \ldots $,
 where the set $V_0$ consists of one element, say $v_0$.
 Let $V:= V_0 \cup V_1 \cup \ldots$.
 It suffices to show that $V$ is a finite set.

 \smallskip
 Let $K$ be the set of all natural numbers $k$ for which $V_k$ is bigger than $V_{k-1}$. For $k \in K$, denote by $v_k$ the
 element of $V_k \smallsetminus V_{k-1}$. Then $V= \{v_k \, | \, k \in K \cup \{0\}\}$.

\smallskip
We define a graph $\varGamma$ as follows: Let $V$ be the set of vertices of $\varGamma$. For each $k \in K$, we introduce
an edge $e_k$ beginning in the class of $f_1(s_k)$ and ending in the class of $f_2(s_k)$. Therefore $e_k$ begins in
$V_{k-1}$ and ends in $v_k$. Since $f \in \mathcal S_0$, there can start only finitely many edges in each vertex. Therefore
$\varGamma$ is a tree to which we can apply K{\"o}nig's lemma if $\varGamma$ is infinite. This leads to an immediate
contradiction to the condition $f \in \mathcal S_1$. So, $\Gamma$ and hence $V$ is finite. \hfill $\Box$
\end{demo}

\subsection{The homomorphism $\varPsi:\BFE \rightarrow \Aut(F(X))$} \label{Psi}

Suppose we are given a map $f:S \rightarrow \mathcal E$ belonging to the admissible set $\mathcal{S}=\mathcal {\tilde{S}}$. Let us write
\[ f(s)= E_{x_s y_s} \, .\]
It should be intuitively clear that we obtain an endomorphism $\alpha_f$ of $F(X)$ by
\[\alpha_f := \prod_{s \in S} E_{x_s y_s}\, .\]

The objective of the present subsection is to make this definition precise, to show that $\alpha_f$ is an automorphism and that
by assigning $\alpha_f$ to $f$ we obtain a homomorphism $\varPsi$ from $\BFE$ to $\Aut(F(X))$.

\medskip
Fix $f \in \mathcal{\tilde{S}}$ as above. For any finite subset $C:= \{s_1,\ldots ,s_n\}$ of $S$ with $s_1 < \ldots < s_n$, we
obtain the automorphism
\[E_C:=E_{x_{s_1} y_{s_1}} \ldots E_{x_{s_n} y_{s_n}} \, .\]

For $x \in X^\pm$, we have defined in Section~\ref{Bas} the finite subset
\[\{s_1, \ldots , s_{\,\frak{n}^+}\}=: C(f,x)\]

with $\frak{n}^+ = \frak{n}^+ (f,x)$. Observe that $C(f,x) = C(f,x^{-1})$. Hence we obtain an endomorphism $\alpha_f$ of $F(X)$ by
\[x \alpha_f:=x E_{C(f,x)} \; \text{for} \; x \in X^\pm \, .\]

\begin{lemma} \label{alpha1}
 Let $f,g \in \mathcal{\tilde{S}}$ with $f \searrow g$. Then $\alpha_f =\alpha_g$.
\end{lemma}

\begin{demo}
 There is a subset $T$ of $S$ with $g= f | S \smallsetminus T$ and there is an involution $\ast$ on $T$ satisfying (\ref{S1}),
 (\ref{S2}), (\ref{S3}). We can write the ordered set $C(f,x)$ in the form
 \[C(f,x) = T_1 S_1 \ldots T_r S_r  \]
 with subsets $T_i$ of $T$ and subsets $S_i$ of $S \smallsetminus T$; of these subsets, only $T_1$ and $S_r$ are allowed to be possibly
 empty. Then we have
 \[C(g,x) = S'_1 \ldots S'_r \]
 where $S'_i$ is a possibly empty subset of $S_i$. As in the proof of Lemma \ref{Sigma1}, we see that each $T_i$ is a
 concatenation of subsets of the form
 \[\sigma_1 < \sigma_2 < \ldots <\sigma_m < \sigma^{\ast}_m < \ldots < \sigma^{\ast}_2 < \sigma^{\ast}_1 \,.\]
 Therefore $E_{T_i} =$ id, and it is clear that the elements of $S_i$ which are not in $S'_i$ do not affect $x$.
 Hence we have

 \parbox{13cm}{\begin{eqnarray*} x E_{C(f,x)} = x E_{C(g,x)} \,.\end{eqnarray*}} \hfill
 \parbox{4mm}{$\Box$}

\end{demo}

\medskip
Observe that for any finite subset $C'$ of $S$ containing $C(f,x)$, we have
\[x \alpha_f =x E_{C'} \,.\]

Hence, for each $w \in F(X)$ and each finite subset $C'$ of $S$ containing the sets $C(f,x)$ for all the letters $x$ of
the word $w$, we have
\[w \alpha_f = w E_{C'}\,.\]
It is not difficult to conclude:
\begin{lemma} \label{alpha2}
 For $f,f' \in \mathcal{\tilde{S}}$, we have $\alpha_{ff'}= \alpha_f \alpha_{f'}$. \hfill $\Box$
\end{lemma}
As an immediate consequence of Lemmas \ref{alpha1} and \ref{alpha2}, we get:
\begin{prop}\label{alpha3}
 For each $f \in \mathcal{S}$, the endomorphism $\alpha_f$ is an automorphism. By $\varPsi [f]:=\alpha_f$ we obtain
 a homomorphism $\varPsi$ from {\text{\rm BF}}$(\mathcal E; \mathcal{S})$ to ${\text{\rm Aut}}(F(X))$. \hfill $\Box$
\end{prop}

\subsection{Refined Nielsen's method}\label{refined}

Here we introduce a {\it complexity of a word} and prove Theorem~\ref{Nielsen} which is a refinement of the classical
Nielsen method for simplifying tuples of elements in free groups (see~\cite[Chapter~1, Proposition~2.2]{LS}).
Condition (1) of this theorem is contained (in a similar form) in Nielsen's method, while
Condition (2) is new. Both conditions are present in Corollary~\ref{NielsenCor}, which will be used later in the proof of Theorem~\ref{general}.

\medskip

First we give some definitions.
Let $F$ be a free group with a basis $X$. We identify the elements of $F$ with reduced words in the alphabet $X^{\pm }$. For any element $w\in F$, we denote by $|w|$ the length of $w$ with respect to $X$.

We define four types of transformations on an arbitrary at most countable tuple of elements $U=(u_1,u_2,\dots )$ of $F$.

\medskip

(I$_i$)\hspace*{3mm} replace $u_i$ by $u_i^{-1}$ ({\it inversion});

\medskip

(R$_{ij}$) replace $u_i$ by $u_iu_j$ where $i\neq j$ ({\it right multiplication});

\medskip

(L$_{ij}$)\hspace*{0.5mm} replace $u_i$ by $u_ju_i$ where $i\neq j$ ({\it left multiplication});

\medskip

(D$_{i}$)\hspace*{1.5mm} delete $u_i$ if $u_i=1$ ({\it deletion}).

\medskip

In all cases it is understood that the $u_k$ for $k\neq i$ remain unchanged. These transformations and the inverse transformations R$_{ij}^{-1}$ and L$_{ij}^{-1}$ are called {\it elementary Nielsen transformations}.

A tuple $U=(u_1,u_2,\dots )$ of elements of $F$ is called {\it Nielsen reduced} if for any three elements $v_1,v_2,v_3$ of the form $u_i^{\pm 1}$, where $u_i\in U$, the following conditions hold:

\medskip

(N1) $v_1\neq 1$;

\medskip

(N2) if $v_1v_2\neq 1$, then $|v_1v_2|\geqslant |v_1|,|v_2|$;

\medskip

(N3) if $v_1v_2\neq 1$ and $v_2v_3\neq 1$, then $|v_1v_2v_3|>|v_1|-|v_2|+|v_3|$.

\medskip

Condition (N2) means that in the product $v_1v_2$ at most half of each factor cancels.
Condition (N3) means that in the product $v_1v_2v_3$ at least one letter of $v_2$ remains uncanceled.

Now we will introduce some notations. Suppose that the set $X\cup X^{-1}$ is totally ordered. This order induces
the {\it graded lexicographical order} $\preccurlyeq$ on the set of all reduced words in the alphabet
$X^{\pm }$ by the following rule.

Let $u$ and $v$ be two reduced words in the alphabet $X^{\pm}$. Denote by $w$ their maximal common initial segment.
We write $u\preccurlyeq v$ if either $|u|<|v|$ or $|u|=|v|$ and the letter of $u$ following $w$ (if it exists) occurs
earlier in the ordering than the letter of $v$ following $w$.

We write $u\prec v$ if $u\preccurlyeq v$ and $u\neq v$. Note that $u\prec v$ implies that $uw\prec vw$ for any $w\in F$, provided that the words $uw$ and $vw$ are reduced.
From now on and to the end of this subsection let $X$ be a finite set.
For any $w\in F$, let $\phi(w)$ denote the cardinality of the set $\{z\,|\, z \preccurlyeq w\}$.

Then $u\prec v\Longleftrightarrow \phi(u)<\phi(v)$ and
$\phi(u)<\phi(v)\Longleftrightarrow \phi(uw)<\phi(vw)$ provided the words $uw$ and $vw$ are reduced.

Let $v\in F$ be a reduced word. By $L(v)$ we denote the initial segment of $v$ of length $\lfloor(|v|+1)/2\rfloor$. The {\it weight} $W(v)$ of the word $v$ is defined to be $W(v)=\phi(L(v))+\phi(L(v^{-1}))$.
Obviously, $W(v)=W(v^{-1})$ and there exists only a finite number of words with weight not exceeding a given natural number. Note also, that $W(u)=W(v)$ does not imply $u=v^{\pm}$.

The {\it complexity} of $v$ is defined to be the pair of nonnegative integer numbers $(|v|,W(v))$ denoted $\C(v)$.
We will write $\C(v)<_{_{lex}}\C(u)$ if
either $|v|<|u|$ or $|v|=|u|$ and $W(v)<W(u)$. Obviously, $\C(u)=\C(u^{-1})$. Note that $\C(v)=\C(u)$ does not imply
$v=u^{\pm 1}$.





\begin{lemma}\label{w-lem} Let $u,v$ be two reduced words in $F$, such that $uv\neq 1$ and $|uv|<|u|=|v|$. Then $W(u)\neq W(v)$.
\end{lemma}

\demo The condition $|uv|<|u|=|v|$ implies that $L(u^{-1})=L(v)$. Suppose $W(u)=W(v)$. Then
$\phi(L(u))+\phi(L(u^{-1}))=\phi(L(v))+\phi(L(v^{-1}))$ and hence $L(u)=L(v^{-1})$.
The conditions $L(u^{-1})=L(v)$ and $L(u)=L(v^{-1})$ imply $uv=1$, a contradiction. \hfill $\Box$

\begin{theorem}\label{Nielsen} Let $F$ be a free group of finite rank. Let $U=(u_1,\dots,u_m)$ be a finite tuple of elements of $F$ with ${\text{\rm rank}}\langle U\rangle =m$, which is not Nielsen reduced.
Then there exists an elementary Nielsen transformation $\mathcal{N}$ of the form $({\text{\rm R}}_{ij})^{\pm 1}$ or $({\text{\rm L}}_{ij})^{\pm 1}$, which carries $U$ to another tuple $U'=(u_1',\dots ,u_m')$, such that the following holds:




\medskip

{\rm (1)} $\C(u_i')<_{_{lex}}\C(u_i)$ and $u_k'=u_k$ for all $k\in \{1,\dots,m\}\smallsetminus \{i\}$,



\medskip
{\rm (2)} $\C(u_j)<_{_{lex}}\C(u_i)$.



\end{theorem}

 \demo
Because of the assumption on the rank, Condition (N1) is satisfied.
We will use the following easy observation:
{\it it is sufficient to prove the theorem for a tuple $(u_1^{\epsilon_1},\dots,u_m^{\epsilon_m})$
with some choice of} $\epsilon_1,\dots,\epsilon_m\in \{1,-1\}$.\hfill ($\ast$)

\medskip

{\bf a)} Suppose that Condition (N2) is not satisfied. Then there are $v_i,v_j$
of the form $u_i^{\pm 1}$, where $u_i\in U$, and such that $v_1v_2\neq 1$ and additionally $|v_1v_2|<|v_1|$ or $|v_1v_2|<|v_2|$.

Using ($\ast$),
we may assume that $v_1,v_2\in U$, say $v_1=u_i$ and $v_2=u_j$.
Without loss of generality, we may also assume that $|u_iu_j|<|u_i|$.
Indeed, if $|u_iu_j|<|u_j|$, then $|u_j^{-1}u_i^{-1}|<|u_j^{-1}|$ and using ($\ast$) again,
we can reduce to the previous situation.

\medskip


Consider what happens if we apply (R$_{ij})$ or (L$_{ji})$ to the relevant part of $U$.
$$\begin{array}{ll}
(u_i,u_j)\overset{({\text{\rm R}}_{ij})}{\longrightarrow}(u_iu_j, u_j) & \vspace*{2mm}\\
(u_i,u_j)\overset{({\text{\rm L}}_{ji})}{\longrightarrow}(u_i, u_iu_j) &
\end{array}$$

If $|u_j|< |u_i|$, we choose $\mathcal{N}=({\text{\rm R}}_{ij})$. Then Condition (2) will be satisfied automatically and Condition (1) will be satisfied by our assumption $|u_iu_j|<|u_i|$.

If $|u_i|<|u_j|$, we choose $\mathcal{N}=({\text{\rm L}}_{ji})$. Then Condition (2) will be satisfied automatically and Condition (1) will be satisfied, since $|u_iu_j|<|u_i|<|u_j|$.

If $|u_j|=|u_i|$, then both $({\text{\rm R}}_{ij})$ and $({\text{\rm L}}_{ji})$ satisfy Condition (1). By Lemma~\ref{w-lem}, we can choose one of them as $\mathcal{N}$ so that Condition (2) will be satisfied.

\medskip

{\bf b)} Suppose that Condition (N2) is satisfied, but Condition (N3) is not.
Then there are $v_1,v_2,v_3$ of the form $u_i^{\pm 1}$, where $u_i\in U$, such that $v_1v_2\neq 1$, $v_2v_3\neq 1$, and
\begin{equation}\label{RefNiels1}
|v_1v_2v_3|\leqslant |v_1|-|v_2|+|v_3|.
\end{equation}

Let $v_1=ap^{-1}$ and $v_2=pb$, where $p$ is the maximal initial segment of $v_2$ canceling in the product $v_1v_2$. Similarly, we write $v_2=cq^{-1}$ and $v_3=qd$, where $q^{-1}$ is the maximal terminal segment
of $v_2$ cancelling in the product $v_2v_3$. By Condition (N2), we have $|p|,|q|\leqslant |v_2|/2$.
Then $v_2=prq^{-1}$ for some $r$. Assuming that $r\neq 1$ we would have
$$|v_1v_2v_3|=|v_1|-|v_2|+|v_3|+2|r|,$$
a contradiction to~(\ref{RefNiels1}). Therefore, $r=1$, $v_2=pq^{-1}$ and $|p|=|q|=|v_2|/2$.
Since $v_2\neq 1$, we obtain $p\neq q$.

\medskip

{\it Case} 1. Suppose that $\phi(p)<\phi(q)$.
Using $(\ast)$, we may assume that $v_2,v_3\in U$, say $v_2=u_i$, $v_3=u_j$. As above, we will look, what happens if we apply (R$_{ij})$ or (L$_{ji})$ to the relevant part of $U$:
$$\begin{array}{ll}
(u_i,u_j)=(pq^{-1},qd)\,\,\overset{({\text{\rm R}}_{ij})}{\longrightarrow}\,\,(pd, qd)=(u_iu_j,u_j) & \vspace*{2mm}\\
(u_i,u_j)=(pq^{-1},qd)\,\,\overset{({\text{\rm L}}_{ji})}{\longrightarrow}\,\,(pq^{-1}, pd)=(u_i,u_iu_j) & \\
\end{array}
$$
Note that $|u_j|\geqslant |u_i|$, since $|u_j|=|u_iu_j|\geqslant |u_i|$, where the last inequality follows from
Condition (N2).
Thus, we consider two subcases.

\medskip

{\it Case} 1.1. Suppose  $|u_i|<|u_j|$. In this case we choose $\mathcal{N}=({\text{\rm L}}_{ji})$.
Then Condition~(2) is trivially satisfied.
Condition (1) is also satisfied, since $|u_iu_j|=|u_j|$ and $W(u_iu_j)<W(u_j)$ because of $\phi(p)<\phi(q)$.

\medskip

{\it Case} 1.2. Suppose $|u_j|=|u_i|$. Then $|d|=|p|=|q|$ and so
\begin{equation}\label{RefNiels2}
|u_iu_j|=|u_i|=|u_j|.
\end{equation}
First we show that $W(u_i)\neq W(u_j)$.
Assume the contrary. Then from
$$\begin{array}{ll}
W(u_i) & = \, \phi(p)+\phi(q),\vspace*{2mm}\\
W(u_j) & = \, \phi(q)+\phi(d^{-1})$$
\end{array}$$
follows that $\phi(p)=\phi(d^{-1})$.
This
implies that $p=d^{-1}$, and hence $v_2v_3=1$, a contradiction.
Thus, it remains to consider two subcases.

\medskip

\hspace*{4mm}{\it Case} 1.2.1. Suppose that $W(u_j)<W(u_i)$. In this case we choose $\mathcal{N}=({\text{\rm R}}_{ij})$. Then Condition~(2) is trivially satisfied.
Now we show, that Condition~(1) is satisfied. In view of~(\ref{RefNiels2}), we shall prove that $W(u_iu_j)<W(u_i)$.


Since $W(u_j)<W(u_i)$, we have $\phi(q)+\phi(d^{-1})<\phi(p) +\phi(q)$ and hence $\phi(d^{-1})<\phi(p)$.
Recall, that we consider Case 1 where $\phi(p)<\phi(q)$.
Then $$W(u_iu_j)=W(pd)=\phi(p)+\phi(d^{-1})<\phi(p)+\phi(q)=W(u_i)$$
and we are done.

\medskip

\hspace*{4mm}{\it Case} 1.2.2. Suppose that $W(u_i)<W(u_j)$. In this case we choose $\mathcal{N}=({\text{\rm L}}_{ji})$. Then Condition~(2) is trivially satisfied. Now we show, that Condition~(1) is satisfied. In view of~(\ref{RefNiels2}), we shall prove, that $W(u_iu_j)< W(u_j)$. The later is valid:
$$W(u_iu_j)=W(pd)=\phi(p)+\phi(d^{-1})<\phi(q)+\phi(d^{-1})=W(u_j).$$

\medskip

{\it Case} 2. Suppose that $\phi(p)>\phi(q)$. Then we consider the pair $(v_2^{-1},v_1^{-1})$ instead of the pair $v_2,v_3$ and reduce to Case 1.\hfill $\Box$

\medskip

\begin{cor}\label{NielsenCor} Let $F$ be a free group of finite rank with a basis $X$. Every finite tuple $U=(u_1,\dots,u_m)$ of elements of $F$ with ${\text{\rm rank}}\langle U\rangle =m$ can be carried to a Nielsen reduced tuple $V=(v_1,\dots ,v_m)$ by
a finite number of Nielsen right and left multiplications, so that at each step Conditions {\rm (1)} and {\rm (2)} of Theorem~{\rm \ref{Nielsen}} are satisfied.

Furthermore, if some $u_i$ is contained in $X^{\pm}$, then the involved Nielsen transformations don't change the elements in place $i$.

Finally, $X^{\pm}\cap \langle U\rangle\subseteq V^{\pm }:=\{v_1^{\pm 1},\dots ,v_m^{\pm 1}\}$.


\end{cor}

\demo For every tuple $W=(w_1,\dots,w_m)$ we define its complexity by the rule $\C(U)=\sum_{i=1}^m \C(w_i)$,
where the sum is component-wise.
By Theorem~{\rm \ref{Nielsen}}, if the tuple $W$ is not Nielsen reduced, we can decrease its complexity
by applying an appropriate Nielsen transformation.
Since the complexity cannot decrease infinitely often, if we start from $U$, after a finite number of steps described in Theorem~{\rm \ref{Nielsen}} we get a Nielsen reduced tuple.

If $u_i\in X^{\pm}$, it will not be changed because of Condition (1).

Let $x\in X^{\pm}\cap \langle U\rangle$. Then $x\in X^{\pm}\cap \langle V\rangle$. Since $V$ is Nielsen reduced,
our claim follows from ~\cite[Corollary~2.4 in Chapter I]{LS} with $w=x$.
\hfill $\Box$

\medskip

In the next section we will use the following lemmas.

\begin{lemma}\label{hard2} Let $F$ be a free group of finite rank with a fixed basis $X$, and let $u,v\in F$.
Suppose that for some $\epsilon,\tau\in \{-1,1\}$ and $u'=(u^{\epsilon}v^{\tau})^{\epsilon}$ the following holds:

\medskip

{\rm (1)} $|u'|=|u|$,

\medskip

{\rm (2)} $\C(u')<_{_{lex}}\C(u)$,

\medskip
{\rm (3)} $\C(v)<_{_{lex}}\C(u)$.

\medskip

Then $L(u')\preccurlyeq L(u)$ and $L(u'^{-1})\preccurlyeq L(u^{-1})$.
\end{lemma}

\demo Condition (3) implies $|v|\leqslant |u|$. Together with $|u'|=|u|$ this gives $L(u^{\epsilon}v^{\tau})=L(u^{\epsilon})$, i.e. $L(u'^{\epsilon})=L(u^{\epsilon})$.
Then, by Condition (2), we have $L(u'^{-\epsilon})\prec L(u^{-\epsilon})$.\hfill $\Box$

\begin{lemma}\label{hard} Let $F$ be a free group of finite rank with a fixed basis $X$, and let $u,v\in F$. Suppose that for some $\epsilon,\tau\in \{-1,1\}$ and $u'=(u^{\epsilon}v^{\tau})^{\epsilon}$ the following holds:

\medskip

{\rm (1)} $|u|=|v|$,

\medskip

{\rm (2)} $\C(u')<_{_{lex}}\C(u)$,

\medskip
{\rm (3)} $\C(v)<_{_{lex}}\C(u)$.

\medskip

Then $L(v^{\tau})=L(u^{-\epsilon})\hspace{2mm} {\text{\rm and}}\hspace*{2mm} L(v^{-\tau})\prec L(u^{\epsilon})$.

\end{lemma}

\demo From Condition (2) we have $|u'|\leqslant |u|$. Together with $|u|=|v|$ this implies that $L(v^{\tau})=L(u^{-\epsilon})$.
From Conditions (3) and $|u|=|v|$ we have $W(v)<W(u)$, i.e. $\phi(L(v^{\tau}))+\phi(L(v^{-\tau}))<\phi(L(u^{\epsilon}))+\phi(L(u^{-\epsilon}))$. Hence $\phi(L(v^{-\tau}))< \phi(L(u^{\epsilon}))$ and so $L(v^{-\tau})\prec L(u^{\epsilon})$. \hfill $\Box$

\begin{notation} It is convenient to unify four different notations into one:\\
For $\epsilon,\tau\in \{-1,1\}$ and $i,j\in \mathbb{N}$ with $i\neq j$ we define
$$T_{i^{\epsilon}j^{\tau}}=
\begin{cases}
R_{ij}, & {\text{\rm if}}\hspace*{2mm} \epsilon=\tau=1\\
R_{ij}^{-1}, & {\text{\rm if}}\hspace*{2mm}  \epsilon=1,\tau=-1\\
L_{ij}, & {\text{\rm if}}\hspace*{2mm}  \epsilon=-1,\tau=-1\\
L_{ij}^{-1}, & {\text{\rm if}}\hspace*{2mm}  \epsilon=-1,\tau=1.\\
\end{cases}$$
\end{notation}

\subsection{The first step towards\\ the surjectivity of $\varPsi:{\text{\rm BF}}(\mathcal{E},\mathcal{S})\rightarrow {\text{\rm Aut}}(F_{\omega})$}\label{key}

For short, we denote by $F_n$ the free group with the finite basis $X_n=\{x_1,\dots,x_n\}$ and by $F_{\omega}$ the free group with the infinite countable basis $X_{\omega}=\{x_1,x_2,\dots \}$.

\medskip

In this section we prove Theorem~\ref{general} which states that, up to a monomial automorphism, every automorphism of $F_{\omega}$ lies in ${\text{\rm im}}\, \varPsi$. In the next section we show that every monomial automorphism of $F_{\omega}$ lies in ${\text{\rm im}}\, \varPsi$. With that the surjectivity of $\varPsi:{\text{\rm BF}}(\mathcal{E},\mathcal{S})\rightarrow {\text{\rm Aut}}(F_{\omega})$ will be proven.

\medskip

{\it Orders.} We will use the following order on $X_n^{\pm }$: $x_1^{-1}\prec x_1\prec \dots \prec x_n^{-1}\prec x_n$. The union of these orders for $n\in \mathbb{N}$ gives the order on $X_{\omega}^{\pm}$.
The corresponding graded lexicographical orders on $F_n$ and $F_{\omega}$ are denoted by $\preccurlyeq_n$ and $\preccurlyeq_{\omega}$. The complexity of a word $v\in F_n$ in the group $F_m$ with $m\geqslant n$ is denoted by
$\C_m(v)$. We emphasize that the complexity of the same word in different groups may be different.

\begin{rmk}\label{order} {\rm 1) If $u\prec_n v$ for some $u,v\in F_n$, then $u\prec_m  v$ for every $m>n$ and, moreover, $u\prec_{\omega}  v$.

2) For any $w\in F_n$, the set $\{z\in F_n\,|\, z \preccurlyeq_n w\}$ is finite. If $|w|\geqslant 2$, then
the set $\{z\in F_{\omega}\,|\, z \preccurlyeq_{\omega} w\}$ is infinite.

3) By nested induction, we see that every decreasing chain
$\dots \prec_{\omega} z_3\prec_{\omega} z_2\prec_{\omega} z_1$ in $F_{\omega}$ is finite.

}
\end{rmk}

In the proof of Theorem~\ref{general}
we will indirectly use the following density lemma (see~\cite[Proposition 4.1 in Chapter I]{LS}).

\begin{lemma}\label{density} For every $\alpha\in {\text{\rm Aut}}(F_{\omega})$ and for every $n$, there exists $m\geqslant n$ and an automorphism $\beta\in {\text{\rm Aut}}(F_m)$, such that $\alpha|_{X_n}=\beta|_{X_n}$.
\end{lemma}


Consider the set $\mathcal{E}=\{E_{xy}\,|\, x,y\in X_{\omega}^{\pm}, y\neq x,x^{-1}\}$
of the elementary Nielsen automorphisms of $F_{\omega}$ of the first kind.
Recall that with every map $f:S\rightarrow \mathcal{E}$ we associate two maps $f_1:S\rightarrow X_{\omega}^{\pm}$
and $f_2:S\rightarrow X_{\omega}^{\pm}$ by the following rule: if $s\in S$
and $f(s)=E_{xy}$, we set $f_1(s)=x$ and $f_2(s)=y$.

\begin{theorem}\label{general} For every automorphism $\alpha$ of $F_{\omega}$ there exists a monomial automorphism $\sigma\in \mathcal{M}(X_{\omega})$ such that $\alpha$ can be written in the form $\sigma\alpha_f$ for some $f:S\rightarrow \mathcal{E}$, where $S$ is either a finite (may be empty) ordered set or $S=(-\mathbb{N})$
and where the following properties are satisfied.

\medskip

{\rm (i)} For every $x\in X_{\omega}^{\pm}$ the set $\{n\in S\,|\, f(n)=E_{xy}\hspace*{2mm} {\text{\rm for some}}\hspace*{2mm}  y\in X_{\omega}^{\pm}\}$ is finite.

{\rm (ii)} There is no injective and order preserving map $\psi: (-\mathbb{N})\rightarrow S$
with conditions $f_1(\psi(n-1))=(f_2(\psi(n)))^{\pm 1}$ for all $n\in (-\mathbb{N})$.
\end{theorem}

\demo 
If there exist $\sigma\in \mathcal{M}(X_{\omega})$ and a finite subset $Y\subset X_{\omega}$ such that $\sigma^{-1}\alpha$ acts identically on $X_{\omega}\setminus Y$, we can take a finite $S$ by Nielsen's theorem.
This theorem says that
every automorphism
of a free group of finite rank is a finite product of automorphisms of type $E_{xy}$ and of a monomial automorphism.
So, we consider the opposite case.

Below in part {\bf a)} we define auxiliary elementary Nielsen transformations, which will be used in part {\bf b)} to construct the function $f$ satisfying this theorem.
Let $X=(x_1,x_2,\dots )$ and $w_i^{(0)}=x_i\alpha^{-1}$, and consider the infinite tuple $W^{(0)}=X\alpha^{-1}=(w_1^{(0)},w_2^{(0)},\dots )$.

\medskip

{\bf a)} We will inductively define Nielsen transformations $\mathcal{N}_t$ and tuples $W^{(t)}=(w_1^{(t)},w_2^{(t)},\dots )$, $t\geqslant 1$, such that they satisfy the following conditions:

\begin{itemize}
\item[\bf{a}1)] $\mathcal{N}_t$ is a right or a left multiplication, or the inverse to them,
\item[\bf{a}2)] $W^{(t)}=(\mathcal{N}_t\circ \dots \circ \mathcal{N}_1)\, W^{(0)}$,
\item[\bf{a}3)] if $\mathcal{N}_{t+1}=T_{i^{\epsilon}j^{\tau}}$, where $t\geqslant 0$, then $\bigl(w_i^{(t+1)}\bigr)^{\epsilon}=\bigl(w_i^{(t)}\bigr)^{\epsilon}\bigl(w_j^{(t)}\bigr)^{\tau}$ and
there exists a natural number $n(t)$, such that
%
%
%
%
%
\begin{eqnarray}\label{Cond.a31}
\C_{n(t)}(w_i^{(t+1)}) & <_{_{lex}}  \C_{n(t)}(w_i^{(t)})
\end{eqnarray}
\begin{eqnarray}\label{Cond.a32}
\C_{n(t)}(w_j^{(t)}) & <_{_{lex}}\C_{n(t)}(w_i^{(t)})
\end{eqnarray}

\item[\bf{a}4)] $\underset{t\rightarrow \infty}{\lim} W^{(t)}=(x_1,x_2,\dots )\sigma^{-1}$ for some $\sigma\in \mathcal{M}(X_{\omega})$.
\end{itemize}

Now we show how to define the transformations $\mathcal{N}_t$, $t\geqslant 1$. We set $\mathcal{N}_0=id$.
Suppose that $\mathcal{N}_0,\dots,\mathcal{N}_p$ are already defined and, for some $r\geqslant 0$, the tuple $W^{(p)}$ contains $x_1,\dots,x_r$, up to inversions, but not $x_{r+1}$.


Since $\alpha$ is an automorphism of $F_{\omega}$, there exists $m$ such that $x_{r+1}\in \langle w_1^{(p)},\dots,w_m^{(p)}\rangle$ and  $x_1,\dots,x_r$ are contained in
$(w_1^{(p)},\dots,w_m^{(p)})$ up to inversions. By Corollary~\ref{NielsenCor},
there exist elementary Nielsen transformations $\mathcal{N}_{p+1},\dots ,\mathcal{N}_q$, such that they satisfy Conditions a1)--a3) and the tuple $W^{(q)}=(\mathcal{N}_q\circ \dots \circ \mathcal{N}_{p+1})\, W^{(p)}$
contains $x_1,\dots,x_{r+1}$ up to inversions. Moreover, these transformations do not change the places of the involved tuples with $x_1^{\pm 1},\dots,x_r^{\pm 1}$. This recursive definition of $\mathcal{N}_t$ will give us $W^{(t)}$, $t\geqslant 1$, which obviously satisfy a4).

\medskip

{\it Remark.} Fix a natural number $n$. By a4), there exists $p$ such that $W^{(p)}$ contains a letter of $X^{\pm}$ in place~$n$.
Then the elementary Nielsen transformations $\mathcal{N}_{t}$, $t>p$, do not have the form $R_{nj}^{\pm 1}$ or $L_{nj}^{\pm 1}$.

\medskip

{\it Claim.} We fix $i\in \mathbb{N}$. If $|w_i^{(n)}|=|w_i^{(m)}|$ for some  $n>m$, then $$L\bigl((w_i^{(n)})^{\pm 1}\bigr)\preccurlyeq L\bigl((w_i^{(m)})^{\pm 1}\bigr).$$

Indeed, by formula~(\ref{Cond.a31}), we have $|w_i^{(t+1)}|\leqslant |w_i^{(t)}|$ for every $t\in \mathbb{N}$. Therefore $|w_i^{(t)}|=|w_i^{(m)}|$ for every $m\leqslant t\leqslant n$. Now the claim follows from a3) and Lemma~\ref{hard2}.

\medskip

{\bf b)} Now we show that $\alpha$ can be written in the form $\sigma \alpha_f$ for some $f:(-\mathbb{N})\rightarrow \mathcal{E}$ satisfying Conditions (i) and (ii). We will use $\sigma$ and the sequence of Nielsen transformations $\mathcal{N}_t$, $t\in \mathbb{N}$, which was defined in {\bf a)}.
Condition~a4) can be written as $\bigl(\dots \mathcal{N}_2\circ\mathcal{N}_1\bigr)\, W^{(0)}=X\sigma^{-1}$.
Since $W^{(0)}=X\alpha^{-1}$, we have
$$X\sigma^{-1}\alpha=\bigl(\dots \mathcal{N}_2\circ\mathcal{N}_1\bigr)X=X\bigl(\dots \mathcal{A}_{-2}\, \mathcal{A}_{-1}\bigr),$$
where $\mathcal{A}_{-t}$ are elementary Nielsen automorphisms of $F_{\omega}$, defined by the rule: $\mathcal{A}_{-t}=E_{x_i^{\epsilon}x_j^{\tau}}$ if $\mathcal{N}_t=T_{i^{\epsilon}j^{\tau}}$ for $i,j\in \mathbb{N}$ and $\epsilon,\tau\in \{-1,1\}$.
Then $$\sigma^{-1}\alpha= \dots \mathcal{A}_{-2}\, \mathcal{A}_{-1}$$ and we will show that the theorem is satisfied for the function $f:(-\mathbb{N})\rightarrow \mathcal{E}$, defined by the rule $f(-t)=\mathcal{A}_{-t}$, $t\in \mathbb{N}$.

Condition (i) follows from Remark in {\bf a)}.



Now we verify Condition (ii). If it is not satisfied, then there exists an infinite sequence of increasing natural numbers $t_1<t_2<t_3<\dots$, such that $\dots\mathcal{A}_{-t_3}, \mathcal{A}_{-t_2},\mathcal{A}_{-t_1}$ have the following form:
$$\dots \mathcal{A}_{-t_3}=E_{x_k^{\rho}x_l^{\theta}}, \hspace*{5mm} \mathcal{A}_{-t_2}=E_{x_j^{\delta}x_k^{\mu}},\hspace*{5mm} \mathcal{A}_{-t_1}=E_{x_i^{\epsilon}x_j^{\tau}}.$$
Then  $$\dots \mathcal{N}_{t_3}=T_{k^{\rho}l^{\theta}}, \hspace*{5mm}\mathcal{N}_{t_2}=T_{j^{\delta}k^{\mu}},\hspace*{5mm} \mathcal{N}_{t_1}=T_{i^{\epsilon}j^{\tau}}.$$







By Condition a3) we have
$$\dots \leqslant |w_k^{(t_2-1)}|\leqslant |w_j^{(t_2-1)}| \leqslant |w_j^{(t_1-1)}|\leqslant|w_i^{(t_1-1)}|.$$
Indeed, the odd inequalities from the right follow from~(\ref{Cond.a32}) and the even ones from~(\ref{Cond.a31}).

Since the length cannot decrease infinitely often, we may assume that all these lengths coincide.
Then we can apply Lemma~\ref{hard} to the following pairs of words and exponents:
$$(w_i^{(t_1-1)},w_j^{(t_1-1)}),\hspace*{2mm} (\epsilon, \tau),$$
$$(w_j^{(t_2-1)},w_k^{(t_2-1)}),\hspace*{2mm} (\delta, \mu),$$
$$(w_k^{(t_3-1)},w_l^{(t_3-1)}),\hspace*{2mm} (\rho, \theta),$$
$$\dots $$

We simplify and unify notations:

$$(u_1,u_2),\hspace*{2mm} (\sigma_1, \sigma_2),$$
$$(u_3,u_4),\hspace*{2mm} (\sigma_3, \sigma_4),$$
$$(u_5,u_6),\hspace*{2mm} (\sigma_5, \sigma_6),$$
$$\dots $$

By applying Lemma~\ref{hard} and Remark~\ref{order}.1) to pairs
$(u_{i},u_{i+1}),\hspace*{2mm} (\sigma_{i}, \sigma_{i+1})$, we obtain
\begin{eqnarray}\label{red}L(u_{i+1}^{\sigma_{i+1}})=L(u_i^{-\sigma_i})\hspace{2mm} {\text{\rm and}}\hspace*{2mm}
L(u_{i+1}^{-\sigma_{i+1}})\prec_{\omega} L(u_i^{\sigma_i})
\end{eqnarray}
for all odd $i\in \mathbb{N}$. By Claim before b), we have also
\begin{eqnarray}\label{blue}L(u_{i+2}^{\pm 1})\preccurlyeq_{\omega} L(u_{i+1}^{\pm 1})
\end{eqnarray}
for all odd $i\in \mathbb{N}$.
\medskip

{\it Claim.} The set $\underset{i\in \mathbb{N}}{\bigcup} \{ L(u_i),L(u_i^{-1})\} $ contains
an infinite subset $\{z_j\,|\, j\in \mathbb{N}\}$, such that
$$\dots \prec_{\omega} z_3\prec_{\omega} z_2\prec_{\omega} z_1.$$

{\it Proof.} We construct a graph $\mathcal{G}$ with edges colored in red and blue.
The vertex set of $\mathcal{G}$ is $\underset{i\in \mathbb{N}}{\bigcup} \{ L(u_i),L(u_i^{-1})\} $,
where $ L(u_i),L(u_i^{-1})$ are considered as formal symbols.
The edges of $\mathcal{G}$ are defined according to formulas~(\ref{red}) and~(\ref{blue}):

For every odd $i\in \mathbb{N}$, the vertices $L(u_i^{-\sigma_i})$ and $L(u_{i+1}^{\sigma_{i+1}})$ are connected by a red edge and the vertices $L(u_i^{\sigma_i})$ and $L(u_{i+1}^{-\sigma_{i+1}})$ are connected by a blue edge.
For every odd $i\in \mathbb{N}$, we connect the vertices $L(u_{i+1})$ and $L(u_{i+2})$ by a red edge and also
the vertices $L(u_{i+1}^{-1})$ and $L(u_{i+2}^{-1})$. An example is given in Figure~4.

\unitlength 1mm 
\linethickness{4.pt}
\ifx\plotpoint\undefined\newsavebox{\plotpoint}\fi 
\begin{picture}(145.275,70)(9,10)
\put(16,66){\circle*{2.5}}
\put(32,66){\circle*{2.5}}
\put(48,66){\circle*{2.5}}
\put(64,66){\circle*{2.5}}
\put(80,66){\circle*{2.5}}
\put(96,66){\circle*{2.5}}
\put(112,66){\circle*{2.5}}
\put(128,66){\circle*{2.5}}
\put(144,66){\circle*{2.5}}
\put(16,50){\circle*{2.5}}
\put(32,50){\circle*{2.5}}
\put(48,50){\circle*{2.5}}
\put(64,50){\circle*{2.5}}
\put(80,50){\circle*{2.5}}
\put(96,50){\circle*{2.5}}
\put(112,50){\circle*{2.5}}
\put(128,50){\circle*{2.5}}
\put(144,50){\circle*{2.5}}
\thicklines
\put(16,66.2){\line(1,-1){16}}
\put(16,66){\line(1,-1){16}}
\put(16,65.8){\line(1,-1){16}}
\thinlines
\put(16,50){\line(1,1){7}}
\put(32,66){\line(-1,-1){7}}
\thicklines
\put(32,66.15){\line(1,0){16}}
\put(32,66){\line(1,0){16}}
\put(32,65.85){\line(1,0){16}}
\put(48,66.15){\line(1,0){16}}
\put(48,66){\line(1,0){16}}
\put(48,65.85){\line(1,0){16}}
\put(64,66.15){\line(1,0){16}}
\put(64,66){\line(1,0){16}}
\put(64,65.85){\line(1,0){16}}
\put(128,66.15){\line(1,0){16}}
\put(128,66){\line(1,0){16}}
\put(128,65.85){\line(1,0){16}}
\put(96,66.15){\line(1,0){16}}
\put(96,66){\line(1,0){16}}
\put(96,65.85){\line(1,0){16}}
\put(64,50.15){\line(1,0){16}}
\put(64,50){\line(1,0){16}}
\put(64,49.85){\line(1,0){16}}
\put(128,50.15){\line(1,0){16}}
\put(128,50){\line(1,0){16}}
\put(128,49.85){\line(1,0){16}}
\put(96,50.15){\line(1,0){16}}
\put(96,50){\line(1,0){16}}
\put(96,49.85){\line(1,0){16}}
\put(32,50.15){\line(1,0){16}}
\put(32,50){\line(1,0){16}}
\put(32,49.85){\line(1,0){16}}
\thinlines
\put(48,50){\line(1,0){16}}
\thicklines
\put(96,66.2){\line(-1,-1){16}}
\put(96,66){\line(-1,-1){16}}
\put(96,65.8){\line(-1,-1){16}}
\put(128,66.2){\line(-1,-1){16}}
\put(128,66){\line(-1,-1){16}}
\put(128,65.8){\line(-1,-1){16}}
\thinlines
\put(80,66){\line(1,-1){7}}
\put(112,66){\line(1,-1){7}}
\put(96,50){\line(-1,1){7}}
\put(128,50){\line(-1,1){7}}
\put(15.25,70){\makebox(0,0)[cc]{$L(u_1)$}}
\put(31.25,70){\makebox(0,0)[cc]{$L(u_2)$}}
\put(47.25,70){\makebox(0,0)[cc]{$L(u_3)$}}
\put(63.25,70){\makebox(0,0)[cc]{$L(u_4)$}}
\put(79.25,70){\makebox(0,0)[cc]{$L(u_5)$}}
\put(95.25,70){\makebox(0,0)[cc]{$L(u_6)$}}
\put(111.25,70){\makebox(0,0)[cc]{$L(u_7)$}}
\put(127.25,70){\makebox(0,0)[cc]{$L(u_8)$}}
\put(143.25,70){\makebox(0,0)[cc]{$L(u_9)$}}
\put(15.25,44){\makebox(0,0)[cc]{$L(u_1^{-1})$}}
\put(31.25,44){\makebox(0,0)[cc]{$L(u_2^{-1})$}}
\put(47.25,44){\makebox(0,0)[cc]{$L(u_3^{-1})$}}
\put(63.25,44){\makebox(0,0)[cc]{$L(u_4^{-1})$}}
\put(79.25,44){\makebox(0,0)[cc]{$L(u_5^{-1})$}}
\put(95.25,44){\makebox(0,0)[cc]{$L(u_6^{-1})$}}
\put(111.25,44){\makebox(0,0)[cc]{$L(u_7^{-1})$}}
\put(127.25,44){\makebox(0,0)[cc]{$L(u_8^{-1})$}}
\put(143.25,44){\makebox(0,0)[cc]{$L(u_9^{-1})$}}
\put(150,50){\makebox(0,0)[cc]{$\dots$}}
\put(150,66){\makebox(0,0)[cc]{$\dots$}}
\end{picture}

\vspace*{-30mm}

\begin{center} Figure 4.
\end{center}

There are two infinite paths in this graph starting at $L(u_1)$ and at $L(u_1^{-1})$ respectively, say $p$ and $q$.
We consider these paths as subgraphs of $\mathcal{G}$. Since $\mathcal{G}=p\cup q$ contains infinitely many blue edges, we may assume that $p$ contains infinitely many blue edges. Let
$(L(u_i^{\epsilon_i}))_{i\in \mathbb{N}}$ be the sequence of vertices of $p$, where $\epsilon_i\in \{-1,1\}$.
Then the corresponding sequence $(L(u_i^{\epsilon_i}))_{i\in \mathbb{N}}$ of elements of $F(X_{\omega})$ contains an infinite subsequence
$(z_j)_{j\in \mathbb{N}}$, such that
$\dots \prec_{\omega} z_3\prec_{\omega} z_2\prec_{\omega} z_1.$ \hfill$\Box$

\medskip

Since the words $u_i$ have the same length, the words $z_j$ have the same length too and we have a contradiction
to Remark~\ref{order}.3). \hfill $\Box$


\subsection{\!The second step towards\\ the surjectivity of $\varPsi:{\text{\rm BF}}(\mathcal{E},\mathcal{S})\rightarrow {\text{\rm Aut}}(F_{\omega})$}
\label{permutations}

Recall that $X_{\omega}=\{x_1,x_2,\dots \}$ is the fixed basis of $F_{\omega}$ and $\mathcal{M}(X_{\omega})$ is the group of monomial automorphisms of $F_{\omega}$.
Let $\mathcal{I}(X_{\omega})$ be the subgroup of $\mathcal{M}(X_{\omega})$ consisting of the automorphisms which send each $x\in X_{\omega}$ to $x$ or $x^{-1}$.

\begin{theorem}\label{monom} Every monomial automorphism of $F_{\omega}$ lies in the image of the homomorphism $\varPsi:{\text{\rm BF}}(\mathcal{E},\mathcal{S})\rightarrow {\text{\rm Aut}}(F_{\omega})$.
\end{theorem}

\demo
Every monomial auto´morphism $\alpha\in {\text{\rm Aut}}(F_{\omega})$ can be expressed in the form $\alpha=\beta\gamma$ for some $\beta\in \mathcal{I}(X_{\omega})$, $\gamma\in \varSigma(X_{\omega})$.

We can express $\gamma$ as a product of independent countable cycles: $\gamma=\prod_{j\in J} \gamma_j$.
By Nielsen, every finite cycle $\gamma_j$ can be expressed in the form $\gamma_j=\nu_j\gamma_j'$, where
$\nu_j$ is the identity or the inversion of an element of $X_{\omega}$, and $\gamma_j'$ is a finite product of $\mathcal{E}$-automorphisms.
If $\gamma_j$ is an infinite countable cycle, it is similar to $\sigma$ in (5). Recall that
\begin{eqnarray}
\sigma = \prod_{i=0}^{\infty} (\tau_{-i,i+1}\; \tau_{-i,-(i+1)}).
\end{eqnarray}
One can check that
$$\tau_{-i,i+1}\; \tau_{-i,-(i+1)}=$$
$$=E_{x_{-(i+1)},x_{-i}}
E_{x_{-(i+1)}^{-1},x_{-i}}
E_{x_{-(i+1)}^{-1},x_{i+1}}
E_{x_{i+1}^{-1},x_{-(i+1)}^{-1}}
\cdot$$
$$E_{x_{-(i+1)},x_{i+1}^{-1}}
E_{x_{i+1},x_{-i}^{-1}}
E_{x_{-i}^{-1},x_{i+1}^{-1}}
E_{x_{i+1},x_{-i}^{-1}}.
$$

So, $\sigma$ lies in the $\varPsi$-image of an infinite product of $\mathcal{E}$-letters. It is not hard to check, that this product is admissible, i.e. lies in $\mathcal{S}$.

Thus, we have to consider an element $\rho\in \mathcal{I}(X_{\omega})$.
If the support of $\rho$ is infinite, then $\rho$ is an infinite product of some $\rho_j\in \mathcal{I}(X_{\omega})$
with disjoint supports of cardi\-nality 2. Each $\rho_j$ is a finite product of $\mathcal{E}$-automorphisms.
So, $\rho$ lies in the $\varPsi$-image of an admissible infinite product of $\mathcal{E}$-letters.
If the support of $\rho$ is finite, then we can write $\rho=\rho'\rho''$ for some $\rho'\rho''\in \mathcal{I}(X_{\omega})$
with infinite supports and so we have reduced to the previous case.\hfill $\Box$

\begin{theorem}\label{final} The homomorphism
$\varPsi:{\text{\rm BF}}(\mathcal{E},\mathcal{S})\rightarrow {\text{\rm Aut}}(F_{\omega})$
is surjective.
\end{theorem}

{\it Proof.} The proof follows immediately from Theorems~\ref{general} and~\ref{monom}.

\subsection{Three topologies on ${\text{\rm BF}}(\mathcal{E},\mathcal{S})$}\label{top}

We consider ${\text{\rm Aut}}(F(X))$ as a topological group with a basis for the neighborhoods of $1$ consisting of the subgroups
$${\text{\rm St}}(Y):=\{\alpha \in {\text{\rm Aut}}(F(X)) \mid y\alpha=y \;\forall \; y \in Y\},$$
where $Y$ runs through the finite subsets of $X$.
The preimage of this topology with respect to the homomorphism $\varPsi: {\text{\rm BF}}(\mathcal{E},\mathcal{S})\rightarrow {\text{\rm Aut}}(F(X))$
gives us a topology on ${\text{\rm BF}}(\mathcal{E},\mathcal{S})$; we denote it by $\frak{T}_{\text{\rm stab}}$ and call the {\it stabilizer topology}.

Recall that the group ${\text{\rm BF}}(\mathcal{E},\mathcal{S})$ also possesses the natural topology (see Section~\ref{naturaltopology});  we denote it by $\frak{T}_{\text{\rm nat}}$.

Let $\frak{T}$ be the topology on ${\text{\rm BF}}(\mathcal{E},\mathcal{S})$ generated by $\frak{T}_{\text{\rm nat}}$ and $\frak{T}_{\text{\rm stab}}$.
To be more precise, a basis for the neighborhoods of $1$ in $\frak{T}$ consists of the sets

$$\mathcal{U}_{\,Y,A}=\varPsi^{-1}({\text{\rm St}}(Y))\,\,\cap \,\,{\text{\rm ker}} \varphi_A,$$
where $Y$ runs through the finite subsets of $X$ and $A$ runs through the finite subsets of~$\mathcal{E}$.
Recall that the map $\varphi_A:{\text{\rm BF}}(\mathcal{E})\rightarrow {\text{\rm BF}}(A)$
sends every map $f:S\rightarrow \mathcal{E}$ from ${\text{\rm BF}}(\mathcal{E})$ to its restriction to $f^{-1}(A)$.

Note that $\mathcal{U}_{\,Y,A}$ are subgroups of ${\text{\rm BF}}(\mathcal{E},\mathcal{S})$ and that
${\text{\rm BF}}(\mathcal{E},\mathcal{S})$ is a topological group with respect to $\frak{T}$.

\subsection{A generalized presentation of ${\text{\rm Aut}}(F_{\omega})$}\label{formulations}

In this section we formulate our main Theorems~\ref{kernel1} and~\ref{kernel2}.
The first theorem describes the kernel of the epimorphism $\varPsi:{\text{\rm BF}}(\mathcal{E},\mathcal{S})\rightarrow {\text{\rm Aut}}(F_{\omega})$\break algebraically, while the second one topologically. We deduce Theorem~\ref{kernel2} from Theorem~\ref{kernel1}.
The proof of Theorem~\ref{kernel1} is based on technical Sections~\ref{inequality} and~\ref{transform}, and it will be given in Section~\ref{proof}. We conclude this section with Theorem~\ref{kernel3} which is a compact reformulation
of Theorem~\ref{kernel1}.

\medskip

We fix the basis $X=\{x_1,x_2,
\dots \}$ of $F_{\omega}$.

\medskip

Now we will work with the group $G={\text{\rm BF}}(\mathcal{E},\mathcal{S})$, where
$\mathcal{E}=\{E_{xy}\,|\, x,y\in X_{\omega}^{\pm }, y\neq x,x^{-1}\}$
and $\mathcal{S}$ is the set of classes of admissible maps defined in Section~\ref{AsS2}.
Let $G_n$ be the subset of $G$, consisting
of all classes of maps $f: S\rightarrow \mathcal{E}$ from $G$, such that $f_1(s)\in \{x_n,x_n^{-1},x_{n+1},x_{n+1}^{-1},\dots \}$ for every $s\in S$.
Note that $G_n$ is a subgroup of $G$ and we have $G=G_1\supset G_2\supset\dots  $, and $\overset{\infty}{\underset{n=1}{\cap}}\, G_n=1$. The sequence of subgroups $(G_n)_{n\in \mathbb{N}}$ is called
the {\it filtration on} $G$.

\begin{rmk}\label{neighborhood} Every neighborhood $\mathcal{U}_{\,Y,A}$ contains some $G_n$. Indeed, we can take
$$n=\max\{k \mid x_k\in Y\}+\max\{k \mid E_{x_k^{\pm 1}\ast}\in A\}+1.$$

\end{rmk}

\begin{defn} {\rm

Let $I$ be a totally ordered set and suppose we are given elements $g_i\in G$ for $i\in I$ such that,
for every $n$, all but finitely many $g_i$ lie in $G_n$. We can form the product $\underset{i\in I}\prod g_i=:g$
according to Section~\ref{infinite products}.

In this situation we call the product $\underset{i\in I}\prod g_i$ {\it admissible with respect to the filtration on} $G$ (or shortly {\it $(G_n)$-admis\-sible}) and will write
$$g=\underset{i\in I}{\underset{(G_n)-\text{\rm adm.}}{\prod}}g_i.$$}
\end{defn}


Observe that the product $g$ is an element of ${\text{\rm BF}}(\mathcal{E})$ and actually, with the notation of Section~\ref{AsS2}, of ${\text{\rm BF}}(\mathcal{E},\mathcal{S}_0)$, but not necessarily of $G_0$.

\medskip

By Theorem~\ref{final},
there is an epimorphism $$\varPsi: {\text{\rm BF}}(\mathcal{E},\mathcal{S})\rightarrow {\text{\rm Aut}}(F_{\omega}).$$

Our aim is to describe a subset $\mathcal{R}\subset {\text{\rm ker}}\,\varPsi$, which in some sense generates ${\text{\rm ker}}\,\varPsi$ (see Theorem~\ref{kernel2}). For that it is convenient to consider the
elements of ${\text{\rm BF}}(\mathcal{E},\mathcal{S})$ as (infinite) words in the alphabet $\mathcal{E}$.

Before we describe some infinite words of $\mathcal{R}$, we present here their finite analogon.
Let $x,y,z\in X^{\pm }$ , such that $y\notin \{x,x^{-1}\}$ and $z\notin \{x,x^{-1},y,y^{-1}\}$. Then the following
relations hold in ${\text{\rm Aut}}(F_{\omega})$:
$$E_{xy}E_{zx}=\bigl(E_{zy^{-1}}E_{zx}\bigr)E_{xy},$$
$$E_{xy}E_{zx^{-1}}=\bigl(E_{zx^{-1}}E_{zy}\bigr)E_{xy},$$
$$E_{xy}E_{zt}=E_{zt}E_{xy},\,\,\, t\notin\{x,x^{-1}\}.$$
Now we describe an ``infinite'' relation.

\begin{defn}\label{admiss}{\rm Let $\alpha=E_{xy}$ be a word consisting of one $\mathcal{E}$-letter (clearly $\alpha\in \mathcal{S}$).
A (possibly infinite) word $\beta\in \mathcal{S}$ is called {\it $\alpha$-admissible} if for each of its letters $E_{zt}$ we have $z\notin \{x,x^{-1},y,y^{-1}\}$.}
\end{defn}

Let $\beta$ be an $\alpha$-admissible word, where $\alpha=E_{xy}$.
We denote by $\beta_{\alpha}$ the word, obtained from $\beta$
by the following replacements of each occurrence of $E_{zx}$ and $E_{zx^{-1}}$:

$$\begin{array}{lll}
E_{zx} & \rightsquigarrow &  E_{zy^{-1}}E_{zx},\\
E_{zx^{-1}} & \rightsquigarrow & E_{zx^{-1}}E_{zy}.\\
\end{array}$$

\medskip

It is easy to check, that $\beta_{\alpha}\in \mathcal{S}$ and that $\varPsi(\alpha\beta)=\varPsi(\beta_{\alpha}\alpha)$
holds in ${\text{\rm Aut}}(F_{\omega})$.
Thus the following words lie in ${\text{\rm ker}}\,\varPsi$:

\medskip

{\bf (}$\mathcal{R}_{\alpha\beta}${\bf )} \hspace*{1mm} $\alpha\beta\alpha^{-1}\beta_{\alpha}^{-1}$, where $\alpha\in \mathcal{E}$ and $\beta$ is an $\alpha$-admissible (infinite) word from $G$.

\medskip

Now, let $F_n$ be the free group of finite rank $n$ with basis $\{x_1,\dots ,x_n\}$.
Consider the canonical epimorphism ${\text{\rm Aut}}(F_n)\rightarrow {\text{\rm GL}}_n(\mathbb{Z})$.
The full preimage of ${\text{\rm SL}}_n(\mathbb{Z})$ with respect to this epimorphism is a subgroup
of index 2 in ${\text{\rm Aut}}(F_n)$, which is denoted by ${\text{\rm SAut}}(F_n)$.
Gersten obtained the following presentation for this subgroup.

\begin{theorem}~\label{Gersten}{\rm \cite{G}}. A presentation for ${\text{\rm SAut}}(F_n)$ is given by generators $\{E_{ab}\,|\, a,b\in X_n^{\pm }, a\neq b,b^{-1}\}$ subject to relations:

\begin{itemize}

\item[\bf ($\mathcal{R}1$)] $E_{ab}^{-1}=E_{ab^{-1}}$,

\item[\bf ($\mathcal{R}2$)] $[E_{ab},E_{cd}]=1$ for $a\neq c,d,d^{-1}$ and $b\neq c,c^{-1}$,

\item[\bf ($\mathcal{R}3$)] $[E_{ab},E_{bc}]=E_{ac}$ for $a\neq c,c^{-1}$,

\item[\bf ($\mathcal{R}4$)] ${\text{\rm w}}_{ab}^{-1}E_{cd}{\text{\rm w}}_{ab}=E_{\sigma(c)\sigma(d)}$,\\
where ${\text{\rm w}}_{ab}$ is defined to be
$E_{ba}E_{a^{-1}b}E_{b^{-1}a}^{-1}$, and $\sigma$ is the monomial map, determined by ${\text{\rm w}}_{ab}$, i.e.
$a\mapsto b^{-1}, b\mapsto a$,

\item[\bf ($\mathcal{R}5$)] ${\text{\rm w}}_{ab}^4=1$.
\end{itemize}
\end{theorem}

\begin{defn}
{\rm We say that a word $W\in {\text{\rm BF}}(\mathcal{E},\mathcal{S})$
is of {\it type} $(\mathcal{R}i)$, if it can be written in the form $UV^{-1}$, where $U=V$ is the relation $(\mathcal{R}i)$
for some $n\in \mathbb{N}$. Let $\mathcal{R}_{\text{\rm fin}}$ be the set of all words of types $(\mathcal{R}1)$ -- $(\mathcal{R}5)$.

We say that a countable word $W\in {\text{\rm BF}}(\mathcal{E},\mathcal{S})$ is of {\it type} $(\mathcal{R}_{\alpha \beta})$ if
it can be written as $\alpha\beta\alpha^{-1}\beta_{\alpha}^{-1}$, where $\beta$ is $\alpha$-admissible.

 Let $\mathcal{R}$ be the set of all countable words of $G$ of types $(\mathcal{R}1)$ -- $(\mathcal{R}5)$ and $(\mathcal{R}_{\alpha\beta})$.}
\end{defn}




\begin{theorem}\label{kernel1} For $G={\text{\rm BF}}(\mathcal{S},\mathcal{E})$ the kernel of the epimorphism
$$\varPsi: G\rightarrow {\text{\rm Aut}}(F_{\omega})$$
coincides with the set of all products $UV$,
where $U,V\in G$ have the form

\begin{equation}\label{UandV}U=\underset{i\in I_1}{\underset{(G_n)-{\text{\rm adm.}}}{\prod}} f_i^{-1}r_i^{\epsilon_i}f_i,\hspace*{5mm}
V=\underset{i\in I_2}{\underset{(G_n)-{\text{\rm adm.}}}{\prod}} f_i^{-1}r_i^{\epsilon_i}f_i,
\end{equation}
where
$I_1$ is $\mathbb{N}$ or a finite initial segment of $\mathbb{N}$,
$I_2$ is $(-\mathbb{N})$ or a finite final segment of $(-\mathbb{N})$,
$f_i\in G$, $r_i\in \mathcal{R}$, $\epsilon_i\in\{-1,1\}$, and $\mathcal{R}$ is the set of all countable words of $G$ of types $(\mathcal{R}1)$ -- $(\mathcal{R}5)$ and $(\mathcal{R}_{\alpha\beta})$.
\end{theorem}

\begin{cor}\label{corollary1} Let $g\in {\text{\rm ker}}\, \varPsi$. Then for every natural $m$ there exist
elements $h_m\in \langle\!\langle \mathcal{R}\rangle\!\rangle_G$ and $g_m\in G_m$, such that $g=h_mg_m$.
\end{cor}

{\it Proof.} We write $g=UV$ for $U,V$ as in~(\ref{UandV}). For every natural $k$ we write $U=U_kU_k'$ and $V=V_k'V_k$, where
$$U_k=\underset{i=1}{\overset{k}{\prod}} f_i^{-1}r_i^{\epsilon_i}f_i,\hspace*{5mm}
V_k=\underset{i=-k}{\overset{-1}{\prod}} f_i^{-1}r_i^{\epsilon_i}f_i.$$

Let $m$ be a given natural number. Since the products in~(\ref{UandV})
are $(G_n)$-admissible, there exists $k$ such that $U_k',V_k'\in G_m$.
For this $k$ we set $h_m=U_k\cdot (U_k'V_k')V_k(U_k'V_k')^{-1}$ and $g_m=U_k'V_k'$.
Clealy, $h_m$ and $g_m$ satisfy the above condition.
\hfill $\Box$

\begin{lemma}\label{R} Every word $g\in G$ of type $(\mathcal{R}_{\alpha\beta})$ lies in the closure of  $\langle\!\langle\mathcal{R}_{\text{\rm fin}}\rangle\!\rangle_G$ in the topology $\frak{T}$.
\end{lemma}

{\it Proof.} Let $g=\alpha\beta\alpha^{-1}\beta_{\alpha}^{-1}$, where $\alpha=E_{xy}$ and $\beta\in G$ is an $\alpha$-admissible (infinite) word. We show that an arbitrary neighborhood $g\mathcal{U}$ of $g$ in the topology $\frak{T}$ contains an element $gu\in\langle\!\langle\mathcal{R}_{\text{\rm fin}}\rangle\!\rangle_G$.
We may assume that $\mathcal{U}=\mathcal{U}_{\,Y,A}$ for some finite subsets $Y\subset X_{\omega}$ and $A\subset \mathcal{E}$.
Let $A'$ be the minimal subset of $\mathcal{E}$, such that

\medskip

1) $A\cup \{E_{xy}\}\subseteq A'$,

2) if $E_{zx^{\epsilon}}\in A'$ for some $z\in X^{\pm 1}$ and $\epsilon\in \{-1,1\}$, then $E_{zy^{-\epsilon}}\in A'$.

\medskip

Clearly, $A'$ is finite. By the choice of $A'$, we have $\varphi_{A'}(g)=\alpha\beta'\alpha^{-1}(\beta')_{\alpha}^{-1}$ for $\beta'=\varphi_{A'}(\beta)$.
We set $u=g^{-1}\varphi_{A'}(g)$. Obviously, $gu\in\langle\!\langle\mathcal{R}_{\text{\rm fin}}\rangle\!\rangle_G$.
Since $g\in {\text{\rm ker}}\,\varPsi$, this implies $u\in {\text{\rm ker}}\,\varPsi$, and since
$\varphi_A\circ \varphi_{A'}=\varphi_A$, we have $u\in {\text{\rm ker}}\varphi_A$.
Hence $u\in \mathcal{U}$.
\hfill $\Box$

\begin{lemma}\label{RR} $\langle\!\langle\mathcal{R}\rangle\!\rangle_G$
lies in the closure of  $\langle\!\langle\mathcal{R}_{\text{\rm fin}}\rangle\!\rangle_G$ in the topology $\frak{T}$.
\end{lemma}

{\it Proof.} Let $g\in \langle\!\langle\mathcal{R}\rangle\!\rangle_G$ and $\mathcal{U}=\mathcal{U}_{\,Y,A}$ be a neighborhood of 1 in the group $G$. Then $g=\overset{k}{\underset{i=1}{\prod}} f_i^{-1}r_if_i$ for some natural $k$, some $r_i\in \mathcal{R}$ and $f_i\in G$. By Lemma~\ref{R}, $r_i=r_i'u_i$ for some $r_i'\in \langle\!\langle\mathcal{R}_{\text{\rm fin}}\rangle\!\rangle_G$ and $u_i\in \mathcal{U}$. Note that $u_i=(r_i')^{-1}r_i\in \mathcal{U}\cap{\text{\rm ker}}\,\varPsi$.
Since $\mathcal{U}\cap {\text{\rm ker}}\,\varPsi={\text{\rm ker}}\varphi_A\cap {\text{\rm ker}}\,\varPsi$ is normal in $G$, we have $$g\in g'
(\mathcal{U}\cap {\text{\rm ker}}\,\varPsi )\subseteq \langle\!\langle\mathcal{R}_{\text{\rm fin}}\rangle\!\rangle_G\, \mathcal{U}$$
for $g'=\overset{k}{\underset{i=1}{\prod}} f_i^{-1}r_i'f_i$.
Since this holds for every $\mathcal{U}$, the proof is completed.\hfill $\Box$

\begin{lemma}\label{RRR} ${\text{\rm ker}}\,\varPsi$
lies in the closure of  $\langle\!\langle\mathcal{R}_{\text{\rm fin}}\rangle\!\rangle_G$ in the topology $\frak{T}$.
\end{lemma}

{\it Proof.} Let $g\in {\text{\rm ker}}\,\varPsi$ and $\mathcal{U}=\mathcal{U}_{\,Y,A}$ be a neighborhood of 1 in $G$. By Remark~\ref{neighborhood}, there exists $m$, such that $G_m\subseteq \mathcal{U}$. By Corollary~\ref{corollary1}
and Lemma~\ref{RR} we have $$g\in \langle\!\langle \mathcal{R}\rangle\!\rangle_G\,G_m\subseteq
\langle\!\langle \mathcal{R}_{\text{\rm fin}}\rangle\!\rangle_G\, \mathcal{U}\,G_m=\langle\!\langle \mathcal{R}_{\text{\rm fin}}\rangle\!\rangle_G\, \mathcal{U}.
$$
\hfill $\Box$

\begin{theorem}\label{kernel2} For $G={\text{\rm BF}}(\mathcal{E},\mathcal{S})$ the kernel of the epimorphism
$$\varPsi:G\rightarrow {\text{\rm Aut}}(F_{\omega})$$
coincides with the closure of $\langle\!\langle\mathcal{R}_{\text{\rm fin}}\rangle\!\rangle_G$
in the topology $\frak{T}$ on $G$:
$${\text{\rm ker}}\, \varPsi=\overline{\langle\!\langle\mathcal{R}_{\text{\rm fin}}\rangle\!\rangle}_G^{\frak{T}}.$$
Here $\mathcal{R}_{\text{\rm fin}}$ is the set of all words of $G$ of types $(\mathcal{R}1)$ -- $(\mathcal{R}5)$.

In particular, $(\mathcal{E},\mathcal{S},\frak{T},\mathcal{R}_{\text{\rm fin}})$
is a generalized presentation of ${\text{\rm Aut}}(F_{\omega })$ of type $\aleph_0$.

\end{theorem}

{\it Proof.} In view of Lemma~\ref{RRR}, it is sufficient to prove that every $g\in \overline{\langle\!\langle\mathcal{R}\rangle\!\rangle}_G^{\frak{T}}$ lies in ${\text{\rm ker}}\,\varPsi$.
Fix a natural $n$.  Since $g\mathcal{U}_{\{x_n\},\emptyset}$ is a neighborhood of $g$, there exists $u\in \mathcal{U}_{\{x_n\},\emptyset}$, such that $gu\in \langle\!\langle\mathcal{R}\rangle\!\rangle_G$.
Then $\varPsi(g)=\varPsi(u^{-1})$ stabilizes $x_n$, and since this holds for every $n$, we have $g\in {\text{\rm ker}}\, \varPsi$.
\hfill $\Box$

\begin{rmk}\label{oberwolfach} ${\text{\rm ker}}\,\varPsi$ is properly contained in the topological closure of $\langle\!\langle\mathcal{R}\rangle\!\rangle_G$ in ${\text{\rm BF}}(\mathcal{S},\mathcal{E})$, with respect to the topology $\frak{T}_{\text{\rm nat}}$ defined in Section~\ref{top}. Indeed, the sequence of elements of ${\text{\rm ker}}\,\varPsi$ $$E_{x_1x_2}E_{x_1x_t}E_{x_tx_2^{-1}}E_{x_1x_t^{-1}}E_{x_tx_2}\,\, (t\in \mathbb{N})$$ converges to
$E_{x_1x_2}\notin {\text{\rm ker}}\,\varPsi$, when $t\rightarrow\infty$.
\end{rmk}

Finally we reformulate Theorem~\ref{kernel1} in a compact form.

\begin{theorem}\label{kernel3}
For $G={\text{\rm BF}}(\mathcal{E},\mathcal{S})$ the kernel of the epimorphism $$\varPsi: G\rightarrow {\text{\rm Aut}}(F_{\omega})$$ consists of  certain conjugates of elements of the set
\begin{equation}\label{set}\{ \underset{i\in I}{\prod_{(G_n)-\text{\rm adm.}}} f_i^{-1}r_i^{\pm 1}f_i\,\mid f_i\in G, r_i\in \mathcal{R}\}\cap G,
\end{equation}
where $I$ is a finite or infinite interval of $\mathbb{Z}$, and $\mathcal{R}$ is the set of all countable words of $G$ of types $(\mathcal{R}1)-(\mathcal{R}5)$ and $(\mathcal{R}_{\alpha\beta})$.
\end{theorem}

{\it Proof.} Let $g\in {\text {\rm ker}}\,\varPsi$. Then $g=UV$ for $U,V$ as in Theorem~\ref{kernel1}, and we can write $g=U(VU)U^{-1}$, where $VU$ lies in the set~(\ref{set}).
The inverse inclusion follows from the next lemma. \hfill $\Box$

\begin{lemma}\label{help} Let $I$ be a linearly ordered set and $g_i\in  {\text{\rm ker}}\,\varPsi$ for any $i\in I$. If
\begin{equation}\label{I}
g=\underset{i\in I}{\prod_{(G_n)-\text{\rm adm.}}} g_i
\end{equation}
and $g\in G$, then $g\in {\text{\rm ker}}\,\varPsi$.
\end{lemma}

{\it Proof.} Fix $x_r\in X_{\omega}$. Since the product~(\ref{I}) is $(G_n)$-admissible, the set $\{i\in I\mid g_i\notin G_{r+1}\}$ is finite. Let $i_1<i_2<\dots <i_k$ be its elements.
Denote
$$h_0=\underset{\{i\in I\mid i<i_1\}}{\underset{(G_n)-{\text{\rm adm.}}}{\prod}}g_i,\hspace*{10mm}
h_k=\underset{\{i\in I\mid i_k<i\}}{\underset{(G_n)-{\text{\rm adm.}}}{\prod}}g_i,\hspace*{10mm}
h_l=\underset{\{i\in I\mid i_l<i<i_{l+1}\}}{\underset{(G_n)-{\text{\rm adm.}}}{\prod}}g_i,$$
where $l=1,\dots,k-1$. Then $g=h_0g_{i_1}h_1g_{i_2}\dots g_{i_k}h_k$.
Since $g\in G$ and every (infinite) subword of a word from $G$ also lies in $G$, we have $h_j\in G$ for $j=0,\dots,k$. Then all $h_i$ belong to $G_{r+1}$, and hence $\varPsi (h_i)$ stabilize $x_r$. Moreover, $\varPsi(g_i)$ is the identity for any $i\in I$ by assumption. Therefore $\varPsi(g)$ stabilizes $x_r$. Since this holds
for every $r$, we have $g\in {\text{\rm ker}}\,\varPsi$. \hfill $\Box$

\begin{rmk} {\rm The following example shows that the assumption $g\in G$ in Lemma~\ref{help} cannot be omitted:
$$g=\underset{i\in \mathbb{N}}{\prod} [E_{x_{i-1}x_{i}},E_{x_{i+1}x_{i+2}}].$$

Indeed, $[E_{x_{i-1}x_{i}},E_{x_{i+1}x_{i+2}}]\in {\text{\rm ker}}\,\varPsi$ for every $i\in \mathbb{N}$, but $g\notin G$, and so $g\notin \!{\text{\rm ker}}\,\varPsi$.}
\end{rmk}

\subsection{An inequality for the length of chains}\label{inequality}

 By a slight abuse of notation, we will consider maps from $\mathcal{T}(\mathcal{E})$ as elements of ${\text{\rm BF}}(\mathcal{E})$.

\begin{defn}\label{chain}
{\rm Let $f:S\rightarrow \mathcal{E}$ be an element of ${\text{\rm BF}}(\mathcal{E})$ and let $x\in X^{\pm}_{\omega}$.
A sequence of elements of $S$, $(\dots, s_2, s_1)$, is called a {\it backward $x$-chain for~$f$}, if
the following conditions hold:

1) $\dots \prec s_2\prec s_1$,

2) $f_1(s_1)\in \{x,x^{-1}\}$,

3) $f_1(s_{i+1})\in \{f_2(s_i),f_2(s_i)^{-1}\}$ for $i\geqslant 1$.

\medskip

We denote by $M(f,x)$ the supremum of lengths of backward $x$-chains for $f$.
Note that $M(f,x) \leqslant \frak{n}^{+}(\bar{f},x)$ (see notation in Section~\ref{Bas}). Hence, by Lemma~\ref{use it}, $M(f,x)$ is finite for every
$f$ in $\mathcal{S}$.
}
\end{defn}

\begin{prop}\label{M-inequality} Let $f=\delta\alpha\beta\gamma$ be a word from $\mathcal{S}$, such that $\beta$ is $\alpha$-admissible; let $f'=\delta\beta_\alpha\alpha\gamma$. Then $M(f',x)\leqslant M(f,x)$ for every $x\in X^{\pm}_{\omega}$.
\end{prop}

\medskip To facilitate the understanding of the proof look at the following example.

\begin{ex}{\rm
$$f=\bigl(\dots \underset{d}{E_{x_2x_5}}\dots\bigr) \overset{\alpha}{\overbrace{\underset{a}{E_{x_1x_2}}}}
\overset{\beta}{\bigl(\overbrace{\,\,\dots \underset{b_1}{E_{x_3x_1}}\dots \underset{b_2}{E_{x_3x_1^{-1}}}\dots}\,\bigr)}\bigl(\dots \underset{c_1}{E_{x_4x_1^{-1}}}\dots \underset{c_2}{E_{x_4x_3^{-1}}}\dots \bigr).$$
$$f'=\bigl(\dots \underset{d}{E_{x_2x_5}\dots }\bigr)\overset{\beta_{\alpha}}{\bigl(\,\overbrace{\dots \underset{b_1'}{E_{x_3x_2^{-1}}}\underset{b_1}{E_{x_3x_1}}\dots \underset{b_2}{E_{x_3x_1^{-1}}}\underset{b_2'}{E_{x_3x_2}}\dots\,\,}\,\bigr)}\overset{\alpha}
{\overbrace{\underset{a'}{E_{x_1x_2}}}}\bigr(\dots \underset{c_1}{E_{x_4x_1^{-1}}}\dots \underset{c_2}{E_{x_4x_3^{-1}}}\dots \bigr).$$

We write down 3 backward $x_4$-chains in $f$:

$$C:\hspace*{10mm}(a,b_2,c_2),\, (d,a,b_1,c_2),\, (d,a,c_1),$$

and 3 backward $x_4$-chains in $f'$:

$$C':\hspace*{10mm}(b_2,c_2),\, (d,b_1',c_2),\, (d,a',c_1).$$
}
\end{ex}

\medskip

{\it Proof of Proposition~\ref{M-inequality}.} Let $f:(D\sqcup \{a\}\sqcup B\sqcup C)\rightarrow \mathcal{E}$
be the map corresponding to the product $f=\delta\alpha\beta\gamma$
and let $f':(D\sqcup B'\sqcup \{a'\}\sqcup C)\rightarrow \mathcal{E}$
be the map corresponding to the product $f'=\delta\beta_{\alpha}\alpha\gamma$. As in the example above, we assume,
that $B'$ is obtained from $B$ by inserting some $b'$-letters.

Let $\mathcal{C}'$ be a backward $x$-chain for $f'$.
It is sufficient to indicate a backward $x$-chain $\mathcal{C}$ for $f$, which is not shorter than $\mathcal{C}'$.
Note that if $\mathcal{C}'$ contains a $b_i'$-letter,
then all letters of this chain with smaller $s$-indices lie in $D$ (since $\beta$ is $\alpha$-admissible).

\medskip

{\it Case 1.} Suppose that $\mathcal{C}'$ does not contain $a'$ and does not contain any $b_i'$-letter.

Then we set $\mathcal{C}=\mathcal{C}'$.

{\it Case 2.} Suppose that $\mathcal{C}'$ does not contain $a'$, but contains some $b_i'$-letter.

Then we define $\mathcal{C}$ to be the chain obtained from $\mathcal{C}'$ by replacing
$b_i'$ by two letters $a, b_i$.

{\it Case 3.} Suppose that $\mathcal{C}'$ contains $a'$.

In this case $\mathcal{C}'$ does not contain elements of $B'$.
Then we define $\mathcal{C}$ to be the chain obtained from $\mathcal{C}'$ by replacing $a'$ by $a$.
\hfill $\Box$

\subsection{Transformation of elements of ${\text{\rm BF}}(\mathcal{E},\mathcal{S})$\\ to elements of
 ${\text{\rm BF}}_{-\mathbb{N}}(\mathcal{E},\mathcal{S})$  modulo ${\text{\rm ker}}\,\varPsi$}\label{transform}

We denote by ${\text{\rm BF}}_{-\mathbb{N}}(\mathcal{E},\mathcal{S})$
the subset of $G={\text{\rm BF}}(\mathcal{E},\mathcal{S})$
consisting of classes of all maps $g: I\rightarrow \mathcal{E}$, where $I$ equals $(-\mathbb{N})$
or a finite final segment of $(-\mathbb{N})$.
The main aim of this section is to prove that any element $f\in {\text{\rm BF}}(\mathcal{E},\mathcal{S})$
can be represented in the form $f=RA$ for some
$R\in {\text {\rm ker}}\, \varPsi$ and
$A\in {\text{\rm BF}}_{-\mathbb{N}}(\mathcal{E},\mathcal{S})$
(see Proposition~\ref{5claims}).

\medskip
In Section~\ref{formulations} we introduced the chain of subgroups $G=G_1\geqslant G_2\geqslant \dots$, such that
$\overset{\infty}{\underset{n=1}{\cap}}\, G_n=1$.
Now we decompose each difference $G_n\setminus G_{n+1}$ into smaller subsets.
By Lemma~\ref{use it}, the number $\frak{n}^{+}(f,x_n)$ is finite for every $f\in G$ and every $n\geqslant 1$.
Now, for every $m\geqslant 1$ we set $$G_{n,m}=\{f\in G_n\setminus G_{n+1}\,|\, \frak{n}^{+}(f,x_n)=m\}.$$

Clearly $G_n\setminus G_{n+1}=\underset{m\geqslant 1}{\bigsqcup} G_{n,m}$.

\medskip

Note that the last term in the sequence~(\ref{S15}) constructed for $f$ and $x=x_n$ is $s_{\frak{n}^{+}(f,x_n)}$. We denote it by $s^{+}(f,x_n)$.

\medskip

\begin{defn}(split and derived forms of $f$)\label{derived}
{\rm Suppose that $(f:S\rightarrow \mathcal{E})\in G_{n}\setminus G_{n+1}$ for some $n$.
Let $s^{+}=s^{+}(f,x_n)$ and
set $$S_1=\{s\in S\,|\, s\prec s^{+}\},\hspace*{2mm} S_2=\{s\in S\,|\, s^{+}\prec s\}.$$ Then we can write $f$ in the form $f=\delta\alpha\beta$, where $\delta=f_{|S_1}$, $\alpha=f_{|\{s^{+}\}}$ and $\beta=f_{|S_2}$. This form will be called the {\it split form} of $f$.

Note that $\beta$ is $\alpha$-admissible (see Definition~\ref{admiss}).
Using the definition of $\beta_{\alpha}$ given after this definition, we rewrite $f$ in the form
$f=r\overline{f}\alpha$, where $r=\delta(\alpha\beta\alpha^{-1}\beta_{\alpha}^{-1})\delta^{-1}$ and $\overline{f}=\delta\beta_{\alpha}$.
The expression $f=r\overline{f}\alpha$ will be called the {\it derived form} of~$f$.}
\end{defn}




\begin{lemma}\label{first} { Let $f\in G_{n,m}$ and let $f=r\overline{f}\alpha$ be the derived form of $f$.
Then the following properties are satisfied.

\medskip

{\bf P1.} $r\in \langle\!\langle G_{n}\cap \mathcal{R}\rangle\!\rangle_{G_{n}}$.

{\bf P2.} $\overline{f}\in G_{n,m-1}$ if $m>1$ and  $\overline{f}\in G_{n+1}$ if $m=1$.

{\bf P3.} $\alpha\in G_{n}$.}
\end{lemma}

{\it Proof.} The proof is straightforward. \hfill $\Box$

\medskip

We can continue and rewrite the element $\overline{f}$ in the derived form.
We want to do that infinitely countably many times.
So, we put $f^{(1)}=f$ and define three sequences of elements of $G$:
$(f^{(i)})_{i\geqslant 1}$, $(r^{(i)})_{i\geqslant 1}$ and $(\alpha^{(i)})_{i\geqslant 1}$, such that
\begin{eqnarray}\label{splitfi}
f^{(i)}=r^{(i)}\overline{f^{(i)}}\alpha^{(i)}
\end{eqnarray}
is the derived form of $f^{(i)}$ and $$f^{(i+1)}=\overline{f^{(i)}}.$$
If $f^{(i)}=1$ for some $i$, we will assume that $f^{(j)}=r^{(j)}=\alpha^{(j)}$ for all natural $j\geqslant i$.

We have
\begin{eqnarray}\label{splitf}
f=\bigl(r^{(1)}r^{(2)}\dots r^{(i)}\bigr)f^{(i+1)}\bigl(\alpha^{(i)}\dots \alpha^{(2)}\alpha^{(1)}\bigr).
\end{eqnarray}

\begin{lemma}\label{second} { Let $f=f^{(1)}\in G_{n,m}$. Then the following properties are satisfied
for every $1\leqslant i\leqslant m$.

\medskip

{\bf P4.}  $f^{(i)}\in G_{n,m-i+1}$; moreover $f^{(m+1)}\in G_{n+1}$.

\medskip

{\bf P5.} $r^{(i)}\in \langle\!\langle G_{n}\cap  \mathcal{R}\rangle\!\rangle_{G_{n}}$; moreover $r^{(m+1)}\in  \langle\!\langle G_{n+1}\cap \mathcal{R}\rangle\!\rangle_{G_{n+1}}$.

{\bf P6.} $\alpha^{(i)}\in G_{n}$; moreover $\alpha^{(m+1)}\in G_{n+1}$.
}
\end{lemma}

{\it Proof.} Property P4 follows by induction from P2. Properties P5 and P6 follow from P4 and~(\ref{splitfi}).
\hfill $\Box$

\begin{prop}\label{5claims} For every element $f\in {\text{\rm BF}}(\mathcal{E},\mathcal{S})$
the following claims hold.

\medskip

{\rm 1)} The sequence $(f^{(i)})_{i\geqslant 1}$ converges to 1 in the topology $\frak{T}_{\text{\rm nat}}$ on ${\text{\rm BF}}(\mathcal{E},\mathcal{S})$, when $i\rightarrow \infty$.

{\rm 2)} The product $R=\underset{i\in \mathbb{N}}{\prod} r^{(i)}$ is well defined, i.e. it belongs to $\text{\rm BF}(\mathcal{E})$.

{\rm 3)} The product $A=\underset{i\in -\mathbb{N}}{\prod} \alpha^{(-i)}$
belongs to ${\text{\rm BF}}_{-\mathbb{N}}(\mathcal{E},\mathcal{S})$.

{\rm 4)} $f=RA$.

{\rm 5)} The product $R=\overset{\infty}{\underset{i=1}{\prod}} r^{(i)}$ is $(G_n)$-admissible
and it belongs to ${\text{\rm BF}}(\mathcal{E},\mathcal{S})$.

{\rm 6)} $R\in {\text{\rm ker}}\, \varPsi$.

\end{prop}

{\it Proof.} Claims 1) and 2)  follow from P4 and P5. Now we prove the most difficult Claim 3).
Clearly $A$ is a map from $(-\mathbb{N})$ to $\mathcal{E}$.
By Property P6 we have $A\in {\text{\rm BF}}(\mathcal{E})$.
We have to verify, that for every $x\in X$, the cardinal number $M(A,x)$ is finite (see Definition~\ref{chain}).

\medskip

Let $f^{(k)}=\delta^{(k)}\alpha^{(k)}\beta^{(k)}$ be the split form of $f^{(k)}$. By definition, we have
$$f^{(k+1)}=\overline{f^{(k)}}=\delta^{(k)}(\beta^{(k)})_{\alpha^{(k)}}.$$ Denote $\gamma^{(k)}=\alpha^{(k-1)}\dots \alpha^{(1)}$. Then we have
$$f^{(k)}\gamma^{(k)}=\delta^{(k)}\alpha^{(k)}\beta^{(k)}\gamma^{(k)}$$
and
$$f^{(k+1)}\gamma^{(k+1)}=\delta^{(k)}(\beta^{(k)})_{\alpha^{(k)}}\,\alpha^{(k)}\gamma^{(k)}.$$

Applying Proposition~\ref{M-inequality} to these two words we get $$M(f^{(k+1)}\gamma^{(k+1)},x)\leqslant M(f^{(k)}\gamma^{(k)},x).$$

for every $x\in X^{\pm 1}_{\omega}$. Since $f=f^{(1)}\gamma^{(1)}$, we have by induction $$M(f^{(k)}\gamma^{(k)},x)\leqslant M(f,x)$$
and so $M(\gamma^{(k)},x)\leqslant M(f,x)$ for every $k\geqslant 1$. This implies $M(A,x)\leqslant M(f,x)$.

\medskip

Claim 4) follows from Equation~(\ref{splitf}) by taking the limit with respect to the natural topology on ${\text{\rm BF}}(\mathcal{E})$ and using Claims 1)-3).

\medskip

5) By the assumption on $f$ and by Claim 3), we have $R=fA^{-1}\in {\text{\rm BF}}(\mathcal{E},\mathcal{S})$.
The $(G_n)$-admissibility of the product $\overset{\infty}{\underset{i=1}{\prod}} r^{(i)}$ follows from P5.

6) This follows from 5) and Lemma~\ref{help}. \hfill $\Box$

\subsection{Proof of Theorem~\ref{kernel1} about ${\text{\rm ker}}\,\varPsi$}\label{proof}


The following technical lemma is used in Lemma~\ref{A}.

\begin{lemma}\label{two} Let $W:(-\mathbb{N})\rightarrow \mathcal{E}$ be an infinite word from $G_{k}\cap {\text{\rm ker}}\,\varPsi$. Then there exists an infinite word $W_2:(-\mathbb{N})\rightarrow \mathcal{E}$ from $G_{k+1}\cap {\text{\rm ker}}\,\varPsi$, and a finite word $W_1\in G_{k}\cap {\text{\rm ker}}\,\varPsi$, such that $W=W_2W_1$.
\end{lemma}

\medskip

{\it Proof.}
We write $W=BC$, where $C$ is a finite subword of $W$ of minimal length which contains all letters $E_{xy}$ with $x\in \{x_{k},x_{k}^{-1}\}$. Clearly $B\in G_{k+1}$, in particular $\varPsi(B)$
stabilizes $F_k$. Let
$$n=\max\{i,j\mid E_{x_ix_j}\hspace*{2mm}{\text{\rm is a letter of}}\hspace*{2mm} C\}+k.$$
Obviously, $n\geqslant k$.
Then $\varPsi(C)_{| F_n}$ is an automorphism of $F_n$, which stabilizes $F_k$ (this follows from $W=BC$, $\varPsi(W)=1$, and the fact that $\varPsi(B)$
stabilizes $F_k$). By Proposition~\ref{stabilizer}, there exists a finite word $C'$ in letters $E_{xy}$ with $x\in \{x_{k+1},\dots ,x_n\}^{\pm 1}$ and $y\in \{x_1,\dots ,x_n\}^{\pm 1}$, such that $\varPsi(C)_{| F_n}=\varPsi(C')_{| F_n}$ in ${\text{\rm Aut}}(F_n)$. We set $W_2=BC'$ and $W_1=(C')^{-1}C$. Then $W_1\in G_k\cap{\text{\rm ker}}\,\varPsi$.
Since $W=W_2W_1\in {\text{\rm ker}}\,\varPsi$, we conclude that $W_2\in G_{k+1}\cap {\text{\rm ker}}\,\varPsi$.
\hfill $\Box$

\begin{lemma}\label{A} Let $A:(-\mathbb{N})\rightarrow \mathcal{E}$ be an infinite word from ${\text{\rm ker}}\,\varPsi$. Then it can be written as
\begin{equation}\label{AU}A=\underset{i\in (-\mathbb{N})}{\underset{(G_n)-\text{\rm adm.}}{\prod}}f_i^{-1}r_i^{\pm 1}f_i,
\end{equation}
for some finite words $f_i$ and $r_i$ in the alphabet $\mathcal{E}$, where $r_i$'s are words of types $(\mathcal{R}1)-(\mathcal{R}5)$.
\end{lemma}

{\it Proof.} By Lemma~\ref{two}, we can write $A$ as $A=(\dots A_3A_2A_1]$ for some finite words
$A_k\in G_{k}\cap {\text{\rm ker}}\,\varPsi$. By Proposition~\ref{stabilizer}, $A_k$ can be written in the form
$$A_k=\underset{j\in J(k)}{\prod}f_j^{-1}r_j^{\pm 1}f_j,$$
where $J(k)$ is a finite set, $f_j$ and $r_j$ are finite words over $\mathcal{E}$, which belong to $G_k$,
and $r_j$'s are words of types $(\mathcal{R}1)-(\mathcal{R}5)$. This completes the proof. \hfill $\Box$

\bigskip

{\it Proof of Theorem~\ref{kernel1}}.
Let $f\in {\text{\rm ker}}\,\varPsi$. If $f$ is represented by a finite $\mathcal{E}$-word, the statement follows from Theorem~\ref{Gersten}.
So, suppose that $f$ is represented by an infinite countable $\mathcal{E}$-word. We show that $f$ can be written in the form $f=UV$ with $U,V$ satisfying the conclusion of Theorem~\ref{kernel1}. First, we write $f=RA$, where $R$ and $A$ as in Proposition~\ref{5claims}.
In particular, $A:(-\mathbb{N})\rightarrow \mathcal{E}$ is an infinite word from $G$ and

$$R\in\{ \underset{i\in \mathbb{N}}{\prod_{(G_n)-\text{\rm adm.}}} f_i^{-1}r_i^{\pm 1}f_i\,\mid f_i\in G, r_i\in \mathcal{R}\cup \{1\}\}\cap G.$$

Since $R\in {\text{\rm ker}}\,\varPsi$ by Proposition~\ref{5claims}, we have $A\in {\text{\rm ker}}\,\varPsi $. Then we can apply Lemma~\ref{A} and write $A$ in the form~(\ref{AU}). So, we can set $U=R$,~$V=A$.

\medskip

Conversely, if $f=UV$ with $U,V$ satisfying the conclusion of Theorem~\ref{kernel1},
then $f\in {\text{\rm ker}}\,\varPsi$ by Lemma~\ref{help}.  \hfill $\Box$

\section{Appendix A: Finite presentations for some stabilizers in ${\text{\rm SAut}}(F_n)$}\label{Stab}


Let $F_n$ be the free group with the basis $X_n=\{x_1,\dots,x_n\}$, where $n$ is finite.
If $a,b\in X_n^{\pm}$ and $a\neq b,b^{-1}$, the {\it Nielsen map} $E_{ba}$ is defined by the rule:
$$
\hspace*{-20.5mm}
\begin{array}{ll}
b\mapsto ba, & \\
x\mapsto x & {\text{\rm if}}\hspace*{2mm} x\in X_n^{\pm }\setminus \{b,b^{-1}\}.
\end{array}
$$

An automorphism of $F_n$ is called {\it monomial}, if it permutes the set $X_n\cup X_n^{-1}$.
Let $\mathcal{M}_n$ be the group of monomial automorphisms of $F_n$. If $a,b\in X_n^{\pm}$ and $a\neq b,b^{-1}$, the monomial automorphism $w_{ab}$ is defined by the rules:
$$
\hspace*{-11.5mm}
\begin{array}{ll}
a\mapsto b^{-1} & \\
b\mapsto a & \\
x\mapsto x & {\text{\rm if}}\hspace*{2mm} x\in X_n^{\pm}\setminus \{a,b,a^{-1},b^{-1}\}.
\end{array}
$$

If $A\subseteq X_n^{\pm}$, and $a\in A, a^{-1}\notin A$, the {\it Whitehead automorphism} $(A,a)$ is defined by
\begin{eqnarray}\label{WN}
(A,a)=\underset{b\neq a}{\underset{b\in A}{\prod}}\,E_{ba}.
\end{eqnarray}

Let $\mathcal{W}_n$ denote the set of all Whitehead automorphisms $(A,a)$.


We now recall McCool's presentation for ${\text{\rm Aut}}(F_n)$.

\begin{theorem}~\label{McCool}{\rm \cite{Mc1}} A finite presentation for ${\text{\rm Aut}}(F_n)$ is given by generators $\mathcal{M}_n\cup \mathcal{W}_n$ and relations:

\begin{itemize}
\item[{\rm \bf (M0)}] A defining set of relations for $\mathcal{M}_n$.

\item[{\rm \bf (M1)}] $(A,a)^{-1}=(A-a+a^{-1},a^{-1})$.

\item[{\rm \bf (M2)}] $(A,a)(B,a)=(A\cup B,a)$ where $A\cap B=\{a\}$.

\item[{\rm \bf (M3)}] $(A,a)(B,b)=(B,b)(A,a)$ where $A\cap B=\emptyset$, $a^{-1}\notin B$, $b^{-1}\notin A$.

\item[{\rm \bf (M4)}] $(A,a)(B,b)=(B,b)(A+B-b,a)$ where $A\cap B=\emptyset$, $a^{-1}\notin B$, $b^{-1}\in A$.

\item[{\rm \bf (M5)}] $(A,a)(A-a+a^{-1},b)=w_{ab}(A-b+b^{-1},a)$ with $w_{ab}\in \mathcal{M}_n$ as above, $b\in A$, $b^{-1}\notin A$, $a\neq b$.

\item[{\rm \bf (M6)}] $\sigma^{-1}(A,a)\sigma=(A\sigma,a\sigma)$ for $\sigma\in \mathcal{M}_n$.
\end{itemize}

\end{theorem}

Using Expression~(\ref{WN}) and rewriting McCool's relations in terms of Nielsen automorphisms $E_{ba}$, Gersten deduced in~\cite[Theorem~1.2]{G}
the following finite presentations of ${\text{\rm Aut}}(F_n)$.

\begin{theorem}~\label{Gersten2}{\rm \cite{G}}. A presentation for ${\text{\rm Aut}}(F_n)$ is given by generators $\{E_{ab}\,|\, a,b\in X_n^{\pm}, a\neq b,b^{-1}\}\cup \mathcal{M}_n$ subject to relations:

\begin{itemize}

\item[\bf ($\mathcal{S}0$)] A defining set of relations for $\mathcal{M}_n$,

\item[\bf ($\mathcal{S}1$)] $E_{ab}^{-1}=E_{ab^{-1}}$,

\item[\bf ($\mathcal{S}2$)] $[E_{ab},E_{cd}]=1$ for $a\neq c,d,d^{-1}$ and $b\neq c,c^{-1}$,

\item[\bf ($\mathcal{S}3$)] $[E_{ab},E_{bc}]=E_{ac}$ for $a\neq c,c^{-1}$,

\item[\bf ($\mathcal{S}4$)] $w_{ab}=E_{ba}E_{a^{-1}b}E_{b^{-1}a^{-1}}$

\item[\bf ($\mathcal{S}5$)] $\sigma^{-1}E_{ab}\sigma=E_{\sigma(a)\sigma(b)}$, $\sigma\in \mathcal{M}_n$.
\end{itemize}
\end{theorem}

By applying the Reidemeister-Schreier method, Gersten obtained in~\cite[Theorem~1.4]{G} the following
presentation for ${\text{\rm SAut}}(F_n)$ (in $\mathcal{R}4$ we correct a misprint in this paper).

\begin{theorem}~\label{Gersten1}{\rm \cite{G}}. A presentation for ${\text{\rm SAut}}(F_n)$ is given by generators $\{E_{ab}\,|\, a,b\in X_n^{\pm 1}, a\neq b,b^{-1}\}$ subject to relations:

\begin{itemize}

\item[\bf ($\mathcal{R}1$)] $E_{ab}^{-1}=E_{ab^{-1}}$,

\item[\bf ($\mathcal{R}2$)] $[E_{ab},E_{cd}]=1$ for $a\neq c,d,d^{-1}$ and $b\neq c,c^{-1}$,

\item[\bf ($\mathcal{R}3$)] $[E_{ab},E_{bc}]=E_{ac}$ for $a\neq c,c^{-1}$,

\item[\bf ($\mathcal{R}4$)] ${\text{\rm w}}_{ab}^{-1}E_{cd}{\text{\rm w}}_{ab}=E_{\sigma(c)\sigma(d)}$,\\
where ${\text{\rm w}}_{ab}$ is defined to be
$E_{ba}E_{a^{-1}b}E_{b^{-1}a}^{-1}$, and $\sigma$ is the monomial map, determined by ${\text{\rm w}}_{ab}$, i.e.
$a\mapsto b^{-1}, b\mapsto a$.

\item[\bf ($\mathcal{R}5$)] ${\text{\rm w}}_{ab}^4=1$.
\end{itemize}
\end{theorem}

\begin{notation} {\rm Let $F_n$ be the free group with basis $X_n=\{x_1,\dots,x_n\}$. For a subset $Y\subseteq F_n$ we denote by ${\text{\rm St}}_{{\text{\rm Aut}}(F_n)}(Y)$
the pointwise stabilizer of $Y$ in ${\text{\rm Aut}}(F_n)$ and by ${\text{\rm St}}_{{\text{\rm SAut}}(F_n)}(Y)$
the pointwise stabilizer of $Y$ in ${\text{\rm SAut}}(F_n)$.}
\end{notation}

\medskip

In~\cite{Mc2},  McCool proved that the group ${\text{\rm St}}_{{\text{\rm Aut}}(F_n)}(Y)$ is finitely presented for every finite $Y\subseteq F_n$ and gave an algorithm for finding such a presentation (see also Proposition~5.7 in Chapter~1 of ~\cite{LS} and a remark after it).

\medskip

We are specially interested in finding presentations of ${\text{\rm St}}_{{\text{\rm Aut}}(F_n)}(X_k)$ and
${\text{\rm St}}_{{\text{\rm SAut}}(F_n)}(X_k)$ for
$X_k=\{x_1,x_2,\dots,x_k\}$, where $1\leqslant k\leqslant n$.
Denote $$\mathcal{M}_{n,k}=\mathcal{M}_n\cap {\text{\rm St}}_{{\text{\rm Aut}}(F_n)}(X_k)$$ and
$$\mathcal{W}_{n,k}=\mathcal{W}_n\cap {\text{\rm St}}_{{\text{\rm Aut}}(F_n)}(X_k).$$
So, $\mathcal{W}_{n,k}$ consists of those $(A,a)\in \mathcal{W}_n$, for which
\begin{eqnarray}\label{restrW}
A\setminus\{a\}\subseteq \{x_{k+1},\dots,x_n\}^{\pm}.
\end{eqnarray}


\begin{prop}\label{stabilizer1} Let $X_n=\{x_1,\dots ,x_n\}$. For every $1\leqslant k\leqslant n$, a finite presentation for ${\text{\rm St}}_{{\text{\rm Aut}}(F_n)}(X_k)$ is given by generators $\mathcal{M}_{n,k}\cup \mathcal{W}_{n,k}$
subject to the relations of ${\rm \bf (M0)}-{\rm \bf (M6)}$ which contain only these generators.
\end{prop}

{\it Proof.} By the result of McCool, ${\text{\rm St}}_{{\text{\rm Aut}}(F_n)}(X_k)$ is isomorphic to the fundamental group of the following 2-dimensional simplicial complex $\mathcal{K}$.

The vertices of $\mathcal{K}$ are $k$-element ordered subsets of $X_n^{\pm}$, i.e. they have the form
$(x_{i_1}^{\epsilon_{i_1}},\dots ,x_{i_k}^{\epsilon_{i_k}})$, where $x_{i_j}\in X_n$ are different and $\epsilon_{i_j}\in \{1,-1\}$.
We distinguish the vertex ${\bf x}=(x_1,\dots,x_k)$.
Two vertices $(u_1,\dots ,u_k)$ and
$(v_1,\dots ,v_k)$ are joint by an edge with label $\alpha$ if $\alpha$ is a Whitehead automorphism of type $(A,a)$
or a monomial automorphism from $\mathcal{M}_n$, such that $v_j=u_j\alpha$ for $j=1,\dots,k$. The inverse edge is labeled  by $\alpha^{-1}$.
The initial vertex of an edge $e$ is denoted by $i(e)$ and the terminal one by $t(e)$.

\medskip

{\it Remark}. Clearly, if an edge $e$ of $\mathcal{K}$ is labelled by a Whitehead automorphism $(A,a)$, then $i(e)=t(e)$

\medskip

We have to define 2-cells of $K$.
Let $\phi(e)$ denote the label of an edge $e$.
This labeling can be obviously extended to paths in the 1-skeleton of $\mathcal{K}$: if $p=e_1e_2\dots e_m$ is a path there,
then we set $\phi(p)=\phi(e_1)\phi(e_2)\dots \phi(e_m)$. By a loop we understand a cyclically ordered sequence of edges
$(e_1,e_2,\dots ,e_m)$, such that $i(e_{j+1})=t(e_j)$, $j=1,\dots,m$, where the indexes are added modulo
$m$.

We glue a 2-cell along a loop $p=(e_1,e_2,\dots ,e_m)$ in the 1-skeleton of $\mathcal{K}$ if there exists a relation
$r=s$ from ${\rm \bf (M0)}-{\rm \bf (M6)}$, such that the
cyclic sequence $(\phi(e_1),\phi(e_2),\dots ,\phi(e_m))$ coincides with the cyclic
word $rs^{-1}$.

Clearly, if $p$ and $p'$ are homotopic paths in $\mathcal{K}$, then $\phi(p)=\phi(p')$. Moreover,
if $p$ is a closed path based at $\bf{x}$ in $\mathcal{K}$, then $\phi(p)$ is an element of ${\text{\rm St}}_{{\text{\rm Aut}}(F_n)}(X_k)$. Thus, $\phi$ induces a homomorphism $\Phi:\pi_1(\mathcal{K},{\bf x})\rightarrow {\text{\rm St}}_{{\text{\rm Aut}}(F_n)}(X_k)$. By the cited general result of McCool, $\Phi$ is an isomorphism.

To describe $\pi_1(\mathcal{K},x)$, we choose a maximal subtree $T$ in $\mathcal{K}$ in the following way.
For every vertex $y\neq \bf{x}$ of $\mathcal{K}$, we choose an edge $e_y$ from $\bf{x}$ to $y$, such that $\phi(e_y)\in \mathcal{M}_n$.
Let $T$ be the maximal subtree in the 1-skeleton of $\mathcal{K}$, consisting of all vertices and all these edges.
For convenience we introduce the formal symbol $e_{\bf x}$, which we identify with $\emptyset$.

For every edge $f$ in $\mathcal{K}$, we denote by $\gamma(f)$ the homotopy class of the path $e_{i(f)}fe_{t(f)}^{-1}$. The elements $\gamma(f)$ generate $\pi_1(\mathcal{K},\bf{x})$.


For every 2-cell in $\mathcal{K}$ with the boundary $f_1f_2\dots f_m$ we write the corresponding relation
$\gamma(f_1)\gamma(f_2)\dots \gamma(f_m)$. These relations together with the trivial relations $\gamma(f)^{-1}=\gamma(f^{-1})$ form a complete set of defining relations for the chosen set of generators
of $\pi_1(\mathcal{K},\bf{x})$.

By applying $\Phi$, we get generators and defining relations for ${\text{\rm St}}_{{\text{\rm Aut}}(F_n)}(X_k)$.
Now we describe precisely these generators and relations.

\medskip

{\it Generators.} Let $f$ be an edge. Denote $\sigma=\phi(e_{i(f)})$ and $\tau=\phi(e_{t(f)})$. By definition of $T$ we have $\sigma,\tau\in \mathcal{M}_n$.

\hspace*{5mm} a) Suppose that $\phi(f)=(A,a)\in \mathcal{W}_n$. By the above remark, the initial and the terminal vertices of $f$ coincide. In particular $\sigma=\tau$.
  It follows that $\Phi(\gamma(f))=\sigma\cdot (A,a)\cdot \sigma^{-1}=(A\sigma^{-1},a\sigma^{-1})$ is an automorphism from $\mathcal{W}_{n,k}$.

\hspace*{5mm} b) Suppose that $\phi(f)=\nu\in \mathcal{M}_n$. Then $\Phi(\gamma(f))=\sigma\nu\tau^{-1}\in \mathcal{M}_{n,k}$.

Thus, $\mathcal{M}_{n,k}\cup \mathcal{W}_{n,k}$ is a generator set for ${\text{\rm St}}_{{\text{\rm Aut}}(F_n)}(X_k)$.


\medskip

{\it Relations.} Let $\Delta$ be a 2-cell in $\mathcal{K}$ with the boundary $f_1f_2\dots f_m$.
By construction, $\phi(f_1)\phi(f_2)\dots \phi(f_m)$ is one of the relations ${\rm \bf (M0)}-{\rm \bf (M6)}$
(up to rewriting of kind $r=s \rightsquigarrow rs^{-1}$), say ${\rm \bf (Mi)}$.
The corresponding relation for ${\text{\rm St}}_{{\text{\rm Aut}}(F_n)}(X_k)$ is $\Phi(\gamma(f_1))\Phi(\gamma(f_2))\dots \Phi(\gamma(f_m))$, which is obviously
a word of length $m$ in elements of $\mathcal{M}_{n,k}\cup \mathcal{W}_{n,k}$. We claim that it also has the form ${\rm \bf (Mi)}$.

For instance, consider a 2-cell with boundary $f_1f_2f_3f_4^{-1}$, such that the equation $$\phi(f_1)\phi(f_2)\phi(f_3)= \phi(f_4)$$
has the form ${\rm \bf (M6)}$:
$\sigma^{-1}\cdot (A,a)\cdot\sigma= (A\sigma,a\sigma)$.

Again by the above remark, $f_2$ and $f_4$ are loops. Denoting $\tau=\phi(e_{i(f_1)}), \delta=\phi(e_{i(f_2)})$,
we can write  the corresponding relation
$$\Phi(\gamma(f_1))\cdot \Phi(\gamma(f_2))\cdot \Phi(\gamma(f_3))=\Phi(\gamma(f_4))$$
in ${\text{\rm St}}_{{\text{\rm Aut}}(F_n)}(X_k)$:
$$\tau\sigma^{-1}\delta^{-1}\cdot (\delta(A,a)\delta^{-1})\cdot\delta\sigma\tau^{-1}=\tau(A\sigma,a\sigma)\tau^{-1}.$$
This is exactly $$\tau\sigma^{-1}\delta^{-1}\cdot (A\delta^{-1},a\delta^{-1})\cdot\delta\sigma\tau^{-1}=(A\sigma\tau^{-1},a\sigma\tau^{-1})$$
and it has the form ${\rm \bf (M6)}$.
The other cases can be considered similarly. \hfill $\Box$

\begin{prop}\label{stabilizer} Let $X_n=\{x_1,\dots ,x_n\}$. For every $1\leqslant k\leqslant n$, a finite presentation for ${\text{\rm St}}_{{\text{\rm SAut}}(F_n)}(X_k)$ is given by generators
\begin{eqnarray}\label{E}
\{E_{ab}\,|\, a\in \{x_{k+1},\dots ,x_n\}^{\pm}, b\in X_n^{\pm}, a\neq b,b^{-1}\}
\end{eqnarray}
subject to those relations {\rm ($\mathcal{R}1$)} -- {\rm ($\mathcal{R}5$)} which contain only these
generators.
\end{prop}

{\it Proof.}
We deduce this proposition from Proposition~\ref{stabilizer1}. Using~(\ref{WN})
we conclude, that ${\text{\rm St}}_{{\text{\rm Aut}}(F_n)}(X_k)$ is generated by $\mathcal{M}_{n,k}$ and all $E_{ba}$
with \begin{eqnarray}\label{restrN}
b\in \{x_{k+1},\dots,x_n\}^{\pm}.
\end{eqnarray}
Now we can rewrite the relations in Proposition~\ref{stabilizer1} in terms of these generators.
We should do that exactly as Gersten in his proof of~\cite[Theorem~1.2]{G}.
As a result, we deduce that ${\text{\rm St}}_{{\text{\rm Aut}}(F_n)}(X_k)$ has the presentation as in Theorem~\ref{Gersten2} with the only restriction, that
all generators of type $E_{ba}$ satisfy Condition~(\ref{restrN}).
Finally, we apply the Reidemeister-Schreier method
to obtain a presentation for the subgroup
${\text{\rm St}}_{{\text{\rm SAut}}(F_n)}(X_k)$. Arguing as in the proof of \cite[Theorem~1.2]{G},
we complete our proof. \hfill $\Box$

\section{Appendix B: Subgroup structure of ${\text{\rm Aut}}(F(X))$.
Problems}\label{Appendix B}

\subsection{Some important subgroups of ${\text{\rm Aut}}(F(X))$}

An automorphism $\alpha$ of $F(X)$ is called {\it bounded}, if there exists a constant $C$ such that $|\alpha(x)|\leqslant C$ and $|\alpha^{-1}(x)|\leqslant C$ for every $x\in X$. The minimal natural number $C$ with
this property will be denoted by $||\alpha||$.
The group of all bounded automorphisms of $F(X)$  is denoted by $\mathcal{B}(X)$.
Clearly, $\langle\mathcal{E}(X), \mathcal{M}(X)\rangle \leqslant \mathcal{B}(X)$.

The following proposition shows that the group ${\text{\rm Aut}}(F(X))$ contains a lot of nonstandard automorphisms
if $X$ is infinite.


\begin{prop}
{\rm 1)} If $X$ is nonempty and finite, then $|{\text{\rm Aut}}(F(X)):\langle \mathcal{E}(X)\rangle|=2$.
Moreover, ${\text{\rm Aut}}(F(X))=\langle \mathcal{E}(X), \tau\rangle $, where
$\tau$ is a monomial automorphism which inverts one chosen letter of $X$ and fixes the others.

{\rm 2)}
If $X$ is infinite, then $|\langle\mathcal{E}(X)\rangle|=|X|$, $|\mathcal{M}(X)|=2^{|X|}$,  and $|{\text{\rm Aut}}(F(X))|=2^{|X|}$. Moreover, $|{\text{\rm Aut}}(F(X)):\mathcal{B}(X)|=|\mathcal{B}(X):\langle\mathcal{E}(X), \mathcal{M}(X)\rangle|=2^{|X|}$.
\end{prop}

\demo
1) This claim is due to Nielsen (see~\cite{LS}).

2) Let $X$ be infinite. The statement about cardinalities of
$\langle\mathcal{E}(X)\rangle$, $\mathcal{M}(X)$, and ${\text{\rm Aut}}(F(X))$
is obvious.
We prove that $|{\text{\rm Aut}}(F(X)):\mathcal{B}(X)|=2^{|X|}$. Consider the set
$$\frak{A}=\{S\subset X\,|\, |X\setminus S|\geqslant |S|=\aleph_0\}.$$ We introduce the equivalence relation $\sim$ on $\frak{A}$ by the following rule: For $S_1,S_2\in \frak{A}$ we write $S_1\sim S_2$, if the symmetric difference $\Delta(S_1,S_2)$ is finite. Clearly, $|(\frak{A}/\sim)| = |\frak{A}|=2^{|X|}$.
For every set $S\in \frak{A}$ we enumerate the elements of $S$ by natural numbers and define the automorphism $\alpha_S$ of $F(X)$ by the rule: $x\mapsto x$ for $x\in X\setminus S$,
and $x_n\mapsto x_2^{-1}(x_nx_1^{n-1})x_2$ for $n\in \mathbb{N}$. Then $S=\{x\in X\,|\, x\alpha_{S}\neq x\}$.

We show that the cosets $\alpha_{S_1}\mathcal{B}(X)$ and $\alpha_{S_2}\mathcal{B}(X)$ are different if $S_1\nsim S_2$. Suppose that $\alpha_{S_2}=\alpha_{S_1}\beta$ for some $\beta\in \mathcal{B}(X)$.
Then the set $$S_2\smallsetminus S_1=\{x\in S_2\,|\, x\alpha_{S_1}=x \}$$
is finite. Indeed, for $x$ from this set we have $|x\alpha_{S_2}|=|x\beta|\leqslant ||\beta||$, which can hold for at most $||\beta||$ elements of $S_2$.
Analogously $S_1\smallsetminus S_2$ is finite.
Hence, $S_1\sim S_2$ and we are done.




Now we prove that $|\mathcal{B}(X):\langle\mathcal{E}(X), \mathcal{M}(X)\rangle|=2^{|X|}$. Choose an element $x_0\in X$ and consider the set $\frak{B}=\{S\subset X\,|\, x_0\notin S\}$.
We introduce the equivalence relation $\sim$ on $\frak{B}$ by the following rule:
For $S_1,S_2\in \frak{B}$ we write $S_1\sim S_2$, if  $\Delta(S_1,S_2)$ is finite.
Clearly, $|\frak{B}/\sim|=|\frak{B}|=2^{|X|}$.
For every set $S\in \frak{B}$ we define the automorphism $\beta_S$ of $F(X)$ by the rule:
$x\mapsto x$ for $x\in X\setminus S$ and $x\mapsto xx_0$ for $x\in S$.

We show, that the cosets $\alpha_{S_1}\langle \mathcal{E}(X),\mathcal{M}(X)\rangle$ and $\alpha_{S_2}\langle \mathcal{E}(X),\mathcal{M}(X)\rangle$ are different if $S_1\nsim S_2$. Suppose that $\alpha_{S_2}=\alpha_{S_1}\beta$ for some $\beta\in \langle \mathcal{E}(X),\mathcal{M}(X)\rangle$.
Then the set $$S_2\setminus S_1=\{x\in S_2\,|\, x\alpha_{S_1}=x \}$$
is finite. Indeed, for $x$ from this set we have $2=|x\alpha_{S_2}|=|x\beta|$, which can hold for only finitely many $x\in X$.
Analogously $S_1\setminus S_2$ is finite.
Hence, $S_1\sim S_2$ and we are done.
\hfill $\Box$

\medskip



\begin{rmk} {\rm

1) Given an infinite $X$, does there exist a countable subset $C\subset {\text{\rm Aut}}(F(X))$, such that
${\text{\rm Aut}}(F(X))=\langle C, \mathcal{M}(X)\rangle$? In~\cite{BE}, Bryant and Evans proved that ${\text{\rm Aut}}(F_{\omega})$ has uncountable confinality (the {\it confinality} of a given group $G$ being the least cardinality of
a chain of proper subgroups whose union is $G$). This implies that if such $C$ exists, then $C$ can be chosen to be finite.

\medskip

2) In~\cite[Corollary 4.4]{T} Tolstykh proved that the group ${\text{\rm Aut}}(F_{\omega})$ has {\it universally finite width}
(recall that $G$ has u.f.w., if for every generating set $S$ of $G$
with $S^{-1}= S$ we have that $G = S^k$ for some natural number $k$).

\medskip

3) We would like to attract attention to the following conjecture of D.~Solitar: For infinite countable $X$,
 the group $\mathcal{B}(X)$ is generated by the set of {\it generalized} elementary automorphisms (for definitions see
 the paper of R.~Cohen~\cite{Cohen}).}

\end{rmk}

\subsection{Normal subgroups of ${\text{\rm Aut}}(F(X))$}

Let $G$ be a countable group. For any subset $Y$ of $G$ we denote by ${\text{\rm St}}(Y)$ the subgroup of ${\text{\rm Aut}}(G)$
consisting of all automorphisms which stabilize all elements of $Y$. Clearly, if $Y$ is finite, then ${\text{\rm St}}(Y)$ has countable index in ${\text{\rm Aut}}(G)$. The group $G$ has the {\it small index property} if the converse holds: for every subgroup $H$ of index less than $2^{\aleph_0}$ in ${\text{\rm Aut}}(G)$ there exists a finite subset $Y$ of $G$ such that ${\text {\rm St}}(Y)\leqslant H$.

The group $F(X_{\omega})$ is known to have the small index property~(\cite{BE}). Similarly
(\cite{New}), if $X$ is a countably infinite set and $H$ is a subgroup of index less than $2^{\aleph_0}$ in $\Sigma(X)$, there exists a finite subset $Y$ of $X$ such that $H$ contains the group $\{\sigma\in \Sigma(X)\,|\,
y\sigma=y\, \hspace*{2mm}{\text{\rm for all}} \hspace*{2mm} y\in Y\}$.
Finally~(\cite{E}), a corresponding property holds for subgroups of ${\text{\rm GL}}(V)$, where $V$ is a vector space of countably infinite dimension.


For background information about the small index property the reader should consult the introduction to~\cite{HHLS}.

\medskip

\begin{theorem}\label{proper normal}  Every proper normal subgroup of ${\text{\rm Aut}}(F(X_{\omega}))$ has index $2^{\aleph_0}$.
\end{theorem}

{\it Proof.} Let $H$ be a normal subgroup of index less than $2^{\aleph_0}$ in ${\text{\rm Aut}}(F(X_{\omega}))$. Since ${\text{\rm Aut}}(F(X_{\omega}))$ has
the small index property, there exists a finite subset $Y\subset X_{\omega}$, such that ${\text{\rm St}}(Y)\leqslant
H$. We will show that $H={\text{\rm Aut}}(F(X_{\omega}))$.

Let $\alpha\in {\text{\rm Aut}}(F(X_{\omega}))$. By the density lemma (Lemma~\ref{density}), there exists a finite
subset $\widetilde{Y}\subset X_{\omega}$, containing $Y$, and there exists an automorphism $\beta\in {\text{\rm Aut}}(F(X_{\omega}))$, such that $\beta\in {\text{\rm St}}(X_{\omega}\smallsetminus \widetilde{Y})$ and $\alpha|_{Y}=\beta|_{Y}$. Then $\alpha\beta^{-1}\in {\text{\rm St}}(Y)\leqslant H$. It is sufficient to prove that $\beta\in H$.
Let $\sigma$ be an automorphism of $F(X_{\omega})$, sending $Y$ to a subset of $X_{\omega}\smallsetminus {\widetilde{Y}}$.
Then $\sigma\beta\sigma^{-1}\in {\text{\rm St}}(Y)\leqslant H$ and since $H$ is normal, we have $\beta\in H$.\hfill $\Box$



\medskip

By~\cite{Cohen} (see also~\cite{Maz}), the natural homomorphism ${\text{\rm Aut}}(F(X))\rightarrow {\text{\rm GL}}(V)$, where $V$ is the free $\mathbb{Z}$-module of dimension $|X|$, is an epimorphism.
So we have the following corollary.


\begin{cor} Let $V$ be a free $\mathbb{Z}$-module of infinite countable dimension. Then
every proper normal subgroup of
${\text{\rm GL}}(V)$ has index $2^{\aleph_0}$.
\end{cor}

\begin{rmk} {\rm In~\cite{T}, Tolstykh proved, with the help of a nice result of~Swan, that the commutator subgroup of ${\text{\rm GL}}(V)$ coincides with ${\text{\rm GL}}(V)$; see Theorem~2.5~(ii) and Proposition~3.1~(i) there.}
\end{rmk}


Is it true that the commutator subgroup of ${\text{\rm Aut}}(F(X_{\omega}))$
coincides with ${\text{\rm Aut}}(F(X_{\omega}))$?
If this is not true, then the commutator subgroup has index $2^{\aleph_0}$.

\section{Acknowledgements} We thank M.~Bridson, A. Klyachko, D.~Segal and P.~Zalesski for helpful discussions. The first named author thanks the MPIM at Bonn for its support and excellent working conditions during the fall 2010,
while this research was finished.

\end{document}